%% file: ON_CURVES_AND_TUBES_IN_RN.tex
\documentclass[leqno,a4]{article}
\usepackage[utf8]{inputenc}
\usepackage{amssymb}
\usepackage{amsthm}
\usepackage{amsfonts}
\usepackage{amsmath}
\usepackage[normalem]{ulem}
\usepackage{tikz}
\usepackage{hyperref}
\usepackage[]{pstricks}
\usepackage{mathrsfs}
\usepackage[]{epsfig}
\usepackage{graphicx}
\usepackage{indentfirst}
\everymath{\displaystyle}
\numberwithin{equation}{section}

\usetikzlibrary{patterns}

\def\Tb{{\mathbf T}}
\def\Nb{{\mathbf N}}
\def\R{{\mathbb R}}
\def\Bb{{\mathbf B}}
\def\du#1{\,{\rm d}u_{#1}}
\def\dY#1{\,{\rm d}y_{#1}}
\def\parl{\left(}\def\parr{\right)}
\def\dfi#1#2{\dparcial{\fifi_{#1}}{{#2}}}
\def\puplacomdois#1#2{(#1_1,\ldots,#1_{#2})}
\def\omegar{\Omega_I}
\def\bgamma{\partial_l\Gamma}\def\esfera#1{{\mathbb S}^{#1}}
\def\desfera#1{{\rm d}{\mathbb S}^{#1}}\def\fifi{\varphi}
\def\dbgamma{{\rm d}(\partial_l\Gamma)}
\def\dvol#1{\dx1\dx2\ldots\dx {#1}}
\def\dupla#1{x_#1,y_#1}\def\elipse#1{{\cal E}_#1}\def\Ti{T_I}
\makeatletter
\newcount\cla\cla=1
\def\claim{\edef\@currentlabel{\the\cla}{\bf Claim \the\cla}\global\advance\cla by 1}
\makeatother
\def\refeq#1{(\ref{#1})}
\usepackage{tikz-cd}
\input matrizes
\def\dx#1{\,{\rm d}x_#1}
\def\vol{{\rm vol}}
\def\area{{\rm area}}
\def\dt{\,{\rm d}t}
\def\dv{\,{\rm d}v}
\def\eps{\epsilon}
\def\dxx{\,{\rm d}x}
\def\dy{\,{\rm d}y}
\def\dz{\,{\rm d}z}
\def\dparcial#1#2{\frac{\partial #1}{\partial #2}}
\def\sen{\sin}
\def\inv #1{#1^{-1}}
\def\Mfr{\cal M}
\def\frsys{{\cal {FS}}}
\def\puplalinhacomt#1#2{(#1_1(t)\;#1_2(t)\;\ldots\; #1_{#2}(t))}
\def\vetlinha#1#2{\trans {\puplalinha {#1}{ #2}}}

\def\vetcol #1#2{\arraycolsep=0pt\def\arraystretch{1.2}
	\left(\begin{array}{c}
		#1_1\\
		#1_2\\
		\vdots\\
		#1_{#2}
	\end{array}\right)}
\def\vetcolcomt#1#2{\arraycolsep=0pt\def\arraystretch{1.2}
	\left({\begin{array}{c}
			#1_1(t)\\
			#1_2(t)\\
			\vdots\\
			#1_{#2}(t)
	\end{array}}\right)}
\def\jdeuman#1#2#3{#2\le #1\le#3}
\def\derivada#1#2{#1^{(#2)}}
\def\V#1{\,V_{#1}}
\def\Vf#1{\,V_{f#1}}
\def\alfa{\alpha}
\def\kapa#1{\, \kappa_{#1}}
\def\kapaf#1{\, \kappa_{f#1}}
\def\vetores#1#2{$#1'(t)$, $#1''(t)$, \ldots,  $#1^{(#2)}(t)$}
\def\funcoes#1{$#1_1(t)$, $#1_2(t)$, \ldots,  $#1_{n-1}(t)$}
\def\funcoessd#1{#1_1(t), #1_2(t), \ldots,  #1_{n-1}(t)}
\def\funcoesst#1{#1_1, #1_2, \ldots,  #1_{n-1}}
\def\frame#1{#1_1(t), #1_2(t), \ldots,  #1_{n}(t)}
\def\framest#1{#1_1, #1_2, \ldots,  #1_{n}}
\def\cal #1{{\mathcal #1}}
\def\nug{\,\til{\nu}}
\def\nuh{\,\barra{\nu}}
\def\dF{{\rm d}F}
\def\trans#1{\vphantom{#1}^{\mbox{\rm\scriptsize T}}#1}
\def\Vg#1{\til{V}_{#1}}
\def\Vh#1{\barra{V}_{#1}}
\def\kapag#1{\,\til{\kappa}_{#1}}
\def\kapah#1{\,\barra{\kappa}_{#1}}
\def\kapatruque#1(#2){#2\,\kapa #1}
\def\**{\displaystyle}
\def\pupla#1#2{(#1_1,#1_2,\ldots, #1_{#2})}
\def\puplasp#1#2{#1_1,#1_2,\ldots, #1_{#2}}
\def\puplaspdet#1#2#3{#1_1(#3),#1_2(#3),\ldots, #1_{#2}(#3)}
\def\puplalinha#1#2{(#1_1\;#1_2\;\ldots\; #1_{#2})}
\def\puplalinhacomtsp#1#2{#1_1(t)\;#1_2(t)\;\ldots\; #1_{#2}(t)}
\def\jesima#1{#1^{\mbox{\scriptsize th}}}
\def\reffrenetf{\V1,\V 2,\ldots,\V n}
\def\reffrenetg{\Vg 1,\Vg 2,\ldots,\Vg n}
\def\reffreneth{\Vh 1,\Vh 2,\ldots,\Vh n}
\def\reffrenetnew#1{\V 1(t),\V2 (t),\ldots,\V{#1}(t)}
\def\reffrenetfs{S\,\V 1,S\,\V  2,\ldots,S\,\V  n}
\def\reffrenetfcomt{\V 1(t),\V2 (t),\ldots,\V n(t)}
\def\kapasf{\kapa 1,\kapa 2,\ldots,\kapa {n-1}}
\def\kapasg{\kapag 1,\kapag 2,\ldots,\kapag {n-1}}
\def\kapash{\kapah 1,\kapah 2,\ldots,\kapah {n-1}}
\def\norma#1{\left\|#1\right\|}
\def\pvetor#1#2{#1_1\wedge#1_2\wedge\cdots\wedge#1_{#2}}
\def\pvetorf#1#2{#1'\wedge #1''\wedge\cdots\wedge #1^{(#2)}}
\def\crossvetor#1#2{#1_1\times#1_2\times\cdots\times #1_{#2}}
\def\crossvetorf#1#2{#1'\times#1''\times\cdots\times #1^{(#2)}}
\def\hzero#1{\hbox to 0pt{\hss #1\hss}}
\usepackage{graphicx}
\usepackage{hyperref}
\hypersetup{colorlinks,
	debug=false,
	linkcolor=blue,  
	citecolor=blue,  
	urlcolor=blue,   
	bookmarksopen=true,
	pdftitle={...},
	pdfauthor={Jadonai},
	pdfsubject={LaTeX},
	pdfkeywords={LaTeX}
	pdfpagemode=FullScreen
}
\usepackage[utf8]{inputenc}
\usepackage{enumerate}
\usepackage{tasks}
\DeclareInstance {tasks} {myenumerate} {default}
{ label =(\roman *)}
\DeclareInstance {tasks} {myenumerateexe} {default}
{ label =[\arabic *]\kern5pt}
\usepackage{amsmath}
\def\i{\mathbf{i}\,}
\long\def\nada++#1++{}
\newtheorem{teo}{\sc Theorem}[section]
\newtheorem{coro}[teo]{\sc Corollary}
\newtheorem{prop}[teo]{\sc Proposition}
\newtheorem{lema}[teo]{\sc Lemma}
\newtheorem{exe}[teo]{\sc Example}
\newtheorem{case}{\sc Case}
\newtheorem{fact}[teo]{\sc Fact}
\newtheorem{defi}[teo]{\sc Definition}
\newtheorem{rem}[teo]{\sc Remark}
\def\beq{\begin{equation}}
\def\eeq{\end{equation}}
\def\benum{\vskip 10pt\begin{enumerate}[(i)]\rm}
	\def\eenum{\end{enumerate}\vskip 5pt}
\def\bteo{\vskip 10pt\begin{teo}\sl}
	\def\eteo{\end{teo}\vskip 10pt\par}
\def\bcoro{\vskip 10pt\begin{coro}\sl}
	\def\ecoro{\end{coro}\vskip 10pt\par}
\def\bprop{\vskip 10pt\begin{prop}\sl}
	\def\eprop{\end{prop}\vskip 10pt\par}
\def\blema{\vskip 10pt\begin{lema}\sl}
	\def\elema{\end{lema}\vskip 10pt\par}
\def\bexe{\begin{exe}\rm}
	\def\bcase{\begin{case}\rm}
		\def\ecase{\end{case}}
	\def\eexe{\end{exe}\vskip 10pt\par}
\def\brem{\vskip 10pt\begin{rem}\rm}
	\def\erem{\end{rem}\vskip 10pt\par}
\def\bdefi{\begin{defi}\rm}
	\def\edefi{\end{defi}\par}
\def\bfact{\vskip 10pt\begin{fact}\sl}
	\def\efact{\end{fact}\vskip 10pt\par}

\def\bsec#1{\section{#1}}

\def\funcao#1#2#3{\hbox{$#1:#2 \longrightarrow #3$}}
\def\R#1{\hbox{${\mathbb R}^{#1}$}}
\def\C#1{\hbox{${\mathbb C}^{#1}$}}

\def\N{\hbox{${\mathbb N}$}}

\def\til#1{\widetilde{#1}}

\def\barra#1{\mbox{$\overline{#1}$}}

\begin{document}
\title{On Curves and Tubes in \R n}

\author{J. Adonai P. Seixas and Isnaldo Isaac Barbosa\\ \footnotesize{Instituto de Matemática}\\
\footnotesize{Universidade Federal de Alagoas}\\
 \footnotesize{Maceió, Alagoas}\\
\footnotesize{adonai@mat.ufal.br and isnaldo@pos.mat.ufal.br}}

\maketitle

\begin{abstract}
In this paper, we present a new explicit formula for the curvatures of a regular curve with an arbitrary parameter in the Euclidean space $\R n$, $n\geq 2$, expressed only in terms of its derivatives. We introduce also the notion of tube with arbitrary cross sections around a curve for which we calculate the volume and give a generalization for the second theorem of Pappus. The first theorem of Pappus is obtained for sphere tubes in arbitrary dimension. 
\end{abstract}


\bsec{Introduction}

In his book, Elementary Differential
Geometry (\cite{Oneilgd}), Barrett O'Neill, makes the following remark:
$$\vbox{\hsize=4.25in 
	... However, for explicit numerical computations--and occasionally for the
	theory as well--this transference is impractical, since it is rarely possible to
	find explicit formulas for. (For example, try to find a unit-speed parametrization
	for the curve $\alfa(t) = (t, t^2, t^3)$.}$$
Such an observation concerns the explicit calculation of the Frenet apparatus of a curve $\alfa$, without using its arc length parametrization, that is, using an arbitrary parameter, as in the cited example. One of our goals here is to generalize this calculation for curves in $\R n$.

Next, we use the material contained mainly in~\cite{Gluck1},~\cite{Kuhnel} and~\cite{Rolfpp}. For this,
let \funcao{f}{I}{\R n}, $n>1$, a parametrized $(n-1)$-regular\footnote{$f$ is said to be $k$-regular if \vetores fk are linearly independent.} curve, that is, the derivatives \vetores f{n-1} are linearly independent, for all $t$ in the interval $I$. In this case, looking carefully  in the above reference, we get $(n-1)$ real functions, defined in $I$, \funcoes{\kappa},
$ \kapa j>0$, $j<n-1$, and a positively oriented orthonormal  frame field  along $f$, ${\mathcal F}_t=\{\frame V\}$, the {\sl Frenet frame of $f$}. The functions $ \kapa j$, $1\le j\le n-1$, will be called of {\sl curvatures} of $f$. When $n=3$, the last curvature,  $ \kapa {2}$, is  called {\sl torsion} and is indicate by $\tau$. The set
$$\cal A_f(t)=\{\funcoessd \kappa,\frame V\}$$
is called {\sl Frenet apparatus} of $f$.  The elements of this set satisfy the following equations, known  {\sl Frenet equations}, in which, $\nu(t)=\norma{f'(t)}$ denotes the speed of $f$ and $\V1(t)=\frac{f'(t)}{\nu(t)}$ the unit tangent field. For simplicity, henceforth,  we omit  the parameter $t$. 
\begin{equation}\label{freneteq}
\left\{
\begin{array}{l}
\V1'= \kapa 1 \nu  \V2   \\
\V j'=- \kapa {j-1}  \nu  \V{j-1}+\kapa {j}  \nu   \V{j+1},\quad 2\le j\le n-1\\
\V{n}'=- \kapa {n-1}  \nu  \V{n-1}.\\
\end{array}\right.
\end{equation}
Furthermore, 
\begin{equation}\label{vndef}
\V{n}=\V1 \times \V2 \times\ldots  \times \V {n-1},
\end{equation}
and, for each $j$, $1\le j\le n-1$, and $t\in I$, the space generated by $\{f',f'',\ldots,f^{(j)}\}$ coincides with that generated by $\{\puplasp V j\}$.
Note that this last fact can be rewritten by using  multivector ($j$-vector) objects, as
$$\pvetorf f j=\lambda \,  \pvetor V j,$$
where $\lambda(t)\neq 0$, for all $t\in I$. In Theorem~\ref{wedgeteomain}, we make explicit  $\lambda$ in terms of the curvatures of $f$.

When  $n=3$, we use the classic notation: 
$\kapa1=\kappa$, $\kapa2=\tau$, $V_1=\Tb$, the unitary tangent vector, $V_2=\Nb$,  the principal normal vector, and $V_3=\Bb=\Tb\times \Nb$,  the binormal vector. Thus, in this case, we have the classic Frenet equations:
$$\left\{
\begin{array}{l}
\Tb'=\kappa\,  \nu\,   \Nb   \\
\Nb'=-\kappa\,   \nu\,   \Tb+\tau\,   \nu \,   \Bb\\
\Bb'=-\tau\,   \nu\,   \Nb.\\
\end{array}\right.,$$
or, in the matrix form,
$$\left(\begin{array}{c}
\Tb' \\
\Nb' \\
\Bb' \\
\end{array}\right)= \nu 
\left(
\begin{array}{ccc}
0&\kappa &0\\
-\kappa &0&\tau \\
0&-\tau &0 \\
\end{array}
\right)
=
\left(\begin{array}{c}
\Tb \\
\Nb \\
\Bb \\
\end{array}
\right).$$
Now, for $n=5$, we get
\beq\label{freneteq5}\left(\begin{array}{c}
	\V1' \\
	\V2' \\
	\V3' \\
	\V4' \\
	\V5' 
\end{array}\right)= \nu 
\left(\begin{array}{ccccc}
	0& \kapa 1 &0&0&0\\
	- \kapa 1 &0& \kapa 2 &0&0\\
	0&- \kapa 2 &0& \kapa 3 &0\\
	0&0&- \kapa 3 &0&\kapa 4 \\
	0&0&0&-\kapa 4 &0
\end{array}\right)
\left(\begin{array}{c}
	\V1 \\
	\V2 \\
	\V3 \\
	\V4 \\
	\V5  \end{array}\right).\eeq
\bsec{Basic Facts on Wedge Product}

\def\pespaco#1#2{\bigwedge^#1(#2)}
Given the $p$-vectors $\pvetor vp,\pvetor wp\in\pespaco p{\R n}$, suppose that
$$\puplalinha vp=\puplalinha wp A,$$
for some $p\times p$ matrix $A=(a_{ij})$ (see Fact~\ref{factum}). This means that each $ v_j$ is a linear combination of the vectors $w_1$, $w_2$, \ldots, $w_p$. More precisely,
$$v_j=\sum_{i=1}^pa_{ij}w_i.$$ It is easy to see that
\begin{equation}\label{eq1}
\pvetor vp=(\det A) \pvetor wp.
\end{equation}
In particular, if $A$ is a triangular matriz, 
\begin{equation}\label{eq2}\pvetor vp=a_{11}a_{22}\ldots a_{pp}\, \pvetor wp.
\end{equation}
Another well known fact is that
\beq\label{pipvetor} (\pvetor vp)\cdot(\pvetor wp)=\det(v_i\cdot w_j)\eeq
defines an intern product in $\pespaco p{\R n}$ and, thus, we have the a well defined norm:
\begin{equation}\label{pvetornorma}
\norma{\pvetor vp}=\sqrt{\det(v_i\cdot v_j)},
\end{equation}
which is also known as {\sl area of parallelogram generated by  the vectors}  $v_1$, $v_2$, \ldots, $v_p$. It is convenient to remark that the $(p\times p)$ matrix $(v_i\cdot v_j)$ equals $\trans{\puplalinha vp}\puplalinha vp$, where $\trans M$ denotes the transpose of $M$ and $\puplalinha vp$ is the $(n\times p)$ matrix whose columns are the column vectors $\puplasp vp$. Thus, when $p=n$, we get a useful identity, namely
\beq\label{normanvetor} \norma{\pvetor vn}=|\det\puplalinha vn|.\eeq
In fact, $$(\det\puplalinha vn)^2=\det(\trans{\puplalinha vn}\puplalinha vn)=\det(v_i\cdot v_j).$$
Now, we consider a $(n-1)$-vector in $\pespaco p{\R n}$, say $w=\pvetor w{n-1}$. In this case, we also have the {\sl cross product} $\til w=\crossvetor w{n-1}$, which belongs to $\R n$ and, for all $X\in\R n$, 
$$\til w\cdot X=\det(\puplasp w{n-1},X)$$
holds.  This implies the $$\{\puplasp w{n-1},\crossvetor w{n-1}\}$$ is a positively oriented basis, whenever $\{\puplasp w{n-1}\}$ is linearly independent.  In fact, $\til w\neq 0$ and 
$$0<\til w\cdot\til w=\det(\puplasp w{n-1},\crossvetor w{n-1}).$$

Of course we cannot compare $w$ with $\til w$, however two $(n-1)$-vectors are equal if, and only if, the corresponding cross products are equal, that is,
\beq\label{pvetecross}\pvetor v{n-1}=\pvetor w{n-1}\Longleftrightarrow\crossvetor v{n-1}=\crossvetor w{n-1}.\eeq
This comes from the fact that the coordinates of $\pvetor v{n-1}$, in the canonical basis of $\pespaco p{\R n}$, coincide, up to signal, with those of the $\crossvetor v{n-1}$.  Indeed,  $$\crossvetor v{n-1}=\ast (\pvetor v{n-1}),$$
where $\ast$ is the Hodge star operator.

\section{The Main Results}
Fixed $t\in I$ and given $X\in\R n$, $1\leq j\le n$, and  we indicate by $[X]_j$ the $\jesima j$ coordinate of $X$ in the Frenet frame $$\cal F(t)=\{\frame V\},$$ 
that is
\def\nabase#1{[X]_1 \V1+[X]_2 \V2+[X]_3 \V3+\cdots+ [X]_n \V n}
$$X=\nabase X.$$
\vskip10pt

Using the above notation, we have the following result.
\bteo{} \label{wedgeteomain}Let \funcao{f}{I}{\R n} be a $(n-1)$-regular parametrized curve in $\R n$ with speed $ \nu=\norma{f'}$  and Frenet apparatus  $$\cal A=\{\funcoesst \kappa,\framest V\}.$$ 
Then 
\benum
\item $[f']_1=\nu$,\quad  $ [f^{(m)}]_m=\kapa1\kapa2\ldots\kapa{m-1} \nu^m,\quad 2\le m\le n;$\\[2pt]
\item $f'=\nu  \V1,\quad$ $f' \wedge f'' \wedge\cdots \wedge  f^{(m)} =
\nu^{\frac{ m(m+1)}2} \kapa 1^{m-1} \kapa 2^{m-2}
\cdots \kapa {m-1} \V1 \wedge  \V2 \wedge\cdots \wedge \V m,\quad  2\le m\le n$.
\eenum
\eteo
\def\span{{\rm span}}
\proof We proceed by induction on $m$. 
We have $[f']_1=f'\cdot \V1= \nu  \V1\cdot \V1=\nu$. Since $ f^{(m-1)}$ belongs to the space generated by 
$\{\puplasp V{m-1}\}$, we have that $f^{(m-1)}\cdot \V {m}=0$, and thus
$$\arraycolsep=0pt\def\arraystretch{1.5}
\begin{array}{rcl}
[f^{(m)}]_m&\;=\;& f^{(m)}\cdot \V m=- f^{(m-1)}\cdot \V {m}'=-\kapa1\kapa2\ldots\kapa{m-2} \nu^{m-1} \V {m-1}\cdot (-\kapa{m-1} \nu \V {m-1})\\
&\;=\;&\kapa1\kapa2\ldots\kapa{m-1} \nu^m,
\end{array}$$
which proofs (i). Now,
$$\arraycolsep=0pt\def\arraystretch{1.5}
\begin{array}{rcl}
f' \wedge f'' \wedge\cdots \wedge  f^{(m)} &\;=\;&(f' \wedge f'' \wedge\cdots \wedge  f^{(m-1)})\wedge f^{(m)}\\
&=& ( \nu^{(1+2+\cdots (m-1))} \kapa 1^{m-2}
\kapa 2^{m-3}
\ldots \kapa {m-2} \V1 \wedge  \V2 \wedge\cdots \wedge \V {m-1})\wedge ( [f^{(m)}]_m\V m)\\
&\;=\;& \nu^{(1+2+\cdots m)} \kapa 1^{m-1} \kapa 2^{m-2}
\ldots \kapa {m-1} \V1 \wedge  \V2 \wedge\cdots \wedge \V m,
\end{array}$$
where, in the last step, we use (i).\qed
\def\normacomp(#1){\norma{#1}}
\vskip 20pt
\def\somanmenosum{\frac{n(n-1)}2}
\def\soman{\frac{n(n+1)}2}
The following corollary establishes an initial algorithm to calculate all curvatures of $f$,  using only the its derivatives. 
\bcoro \label{coro2}\mbox{}
\beq\label{firsteq}\kapa {1}=\frac{\norma{f' \wedge f''}}
{ \nu^{3}},\quad \kapa {m}=\frac{\norma{f' \wedge f'' \wedge\cdots \wedge  f^{(m+1)}}}
{ \nu^{m+1} \kapa 1 \kapa 2\ldots \kapa {m-1} \norma{f' \wedge f'' \wedge\cdots \wedge  f^{(m)}}},\quad 2\le m\le n-2\eeq
and, the last curvature,  which has a sign, is
\beq\label{secteq}\kapa {n-1}=\frac{(f' \times f'' \times\cdots \times  f^{(n-1)})\cdot  f^{(n)}}
{ \nu^{n} \kapa 1\kapa 2\ldots \kapa {n-2} \norma{f' \wedge f'' \wedge\cdots \wedge  f^{(n-1)}}}.\eeq
Furthermore,
\beq\label{tereq} \V {n}=\frac{f' \times f'' \times\cdots \times  f^{(n-1)}}
{\norma{f' \wedge f'' \wedge\cdots \wedge  f^{(n-1)}}}.
\eeq
\ecoro
\proof Taking the norm in  Theorem~\ref{wedgeteomain}-(ii), and noting $\pvetor V{m}$ and $\pvetor V{m-1}$ are unit vectors, produces
$$\norma{f' \wedge f'' \wedge\cdots \wedge  f^{(m+1)}}= \nu^{\frac{(m+1)(m+2)}2} \kapa 1^{m} \kapa 2^{m-1}
\cdots \kapa {m-1}^2 \kapa {m}$$
and
$$\norma{f' \wedge f'' \wedge\cdots \wedge  f^{(m)}}= \nu^{\frac{ m(m+1)}2} \kapa 1^{m-1} \kapa 2^{m-2}
\cdots \kapa {m-1}.$$
Dividing these equations,  we get~(\ref{firsteq}).
These arguments show that~(\ref{firsteq}) also applies to $m = n-1$,  but in this way we lost the signal of $\kapa{n-1}$, that is,  we only obtain  $|\kapa{n-1}|$. For this reason, we rewrite Theorem~\ref{wedgeteomain}-(ii)  using  the cross product, according to equation~(\ref{pvetecross}), in our introduction:
$$\crossvetorf f{n-1} =
\nu^{\somanmenosum} \kapa 1^{n-2} \kapa 2^{n-3}\ldots \kapa {n-2}  \V{n},$$
for $\V n= \crossvetor V{n-1}$.
Since $\kapa j>0$,  $1\le j\le n-2$, this last equation implies (\ref{tereq}) and
$$\arraycolsep=0pt\def\arraystretch{2}
\begin{array}{rcl}
(\crossvetorf f{n-1})\cdot f^{(n)}&\;=\;&\**( \nu^{\somanmenosum} \kapa 1^{n-2} \kapa 2^{n-3}\ldots \kapa {n-2}  \V{n})\cdot( [f]_n\V n)\\
&\;=\;&\** \nu^{{\soman}} \kapa 1^{n-1} \kapa 2^{n-2}\ldots \kapa {n-1}\\
&=& \nu^{n} \kapa 1 \kapa 2\ldots \kapa {n-2} \norma{f' \wedge f'' \wedge\cdots \wedge  f^{(n-1)}} \kapa {n-1},
\end{array}$$
which implies~(\ref{secteq}).\qed
\vskip 20pt
The next corollary certainly gives us the more efficient way to get the curvatures in terms of the derivatives only. Before, we will do a cancellation lemma.
\blema \label{lema3} Let \funcao p{\N}{\R{}} be a positive function with $p(0)=1$ and  $p(1)=c$. Now define $q(n)=\frac{p(n-1)p(n+1)}{c \,p^2(n)}$. Then $q(n)$ equals  $$\frac{1}{c^{n+1}q(1)q(2)\ldots q(n-1)} \frac{p(n+1)}{p(n)}.$$
\elema 
\proof Using the definition of $q$, it easy to check that the factors in 
$q(1)q(2)\ldots q(n-1)$ cancel  nicely resulting
$\frac{p(n)}{ c^{n}p(n-1)}$. (Here the reader could use induction on $n$, observing that $q(1)=p(2)/c^3$.) Hence
$$\frac{1}{c^{n+1}q(1)q(2)\ldots q(n-1)} \frac{p(n+1)}{p(n)}=\frac1{c^{n+1}}\frac{c^{n}p(n-1)}{p(n)}\frac{p(n+1)}{p(n)}=q(n),$$
and the proof is complete.\qed
\bcoro  \label{coro4}Given $2\le m\leq n-2$,
\beq
\kapa {m}=\frac{\norma{f' \wedge f'' \wedge\cdots \wedge  f^{(m-1)}} \norma{f' \wedge f'' \wedge\cdots \wedge  f^{(m+1)}}}
{ \nu  \norma{f' \wedge f'' \wedge\cdots \wedge  f^{(m)}}^2}
\eeq
and
\beq
\kapa {n-1}=\frac{\norma{f' \wedge f'' \wedge\cdots \wedge  f^{(n-2)}}\, ((\crossvetorf f{n-1})\cdot f^{(n)})}
{ \nu  \norma{f' \wedge f'' \wedge\cdots \wedge  f^{(n-1)}}^2}.
\eeq
\proof Just set $p(n)=\norma{\pvetorf f n}$ and use Lemma~\ref{lema3} together with Corollary~\ref{coro2}.\qed
\ecoro
\vskip 10pt
The following corollary establishes an algorithm that calculates,  using only the derivatives, part of  the Frenet apparatus (up to $\kapa 4$ and $\V3$) of $f$, for any $n$.
\bcoro\mbox{}\label{coro5}
\vskip10pt
\begin{tasks}[style=myenumerate,label-format={\normalfont},label-width={5mm}, item-indent={20pt}, item-format={\sl},label-align=right, after-item-skip={4pt}, column-sep={80pt}](2)
	\task \ \ \ \ \ \ \ $\V 1= \frac{f'}{\nu}$;\\[2pt]
	\task \ \ \ \ \ \ \ $\kapa1=\frac{\normacomp(f'\wedge f'')}{\nu^3}$;\\[2pt]
	\task*(2)\ \ \ \ \ \ \ $\V2
	=\frac{\nu f''-f' \nu'}{\nu^3 \kapa 1}=\frac{\nu  f''-f' \nu'}{\normacomp(f'\wedge f'')}=\frac{-(f' \cdot f'') f' +\nu^2 f''}{\nu \normacomp(f'\wedge f'')}=-\left(\frac{ f' \cdot f''}{\nu \normacomp(f'\wedge f'')}\right)f'+\frac{\nu}{\normacomp(f'\wedge f'')} f''
	$;\label{coro5i3}\\[2pt]
	\task \ \ \ \ \ \ \ $\kapa 2=\frac{\normacomp(f'\wedge f''\wedge f''')}{\normacomp(f'\wedge f'')^2}$;\\[2pt]
	\task \ \ \ \ \ \ \ $k_3= \frac{\normacomp(f'\wedge f'') \, \normacomp(f'\wedge f''\wedge f'''\wedge f'''')}{\nu \normacomp(f'\wedge f''\wedge f''')^2}$;\\[2pt]
	\task \ \ \ \ \ \ \  $k_4= \frac{\normacomp(f'\wedge f'')\,  \normacomp(f'\wedge f''\wedge f'''\wedge f''''\wedge f''''')}{\nu \normacomp(f'\wedge f''\wedge f'''\wedge f'''')^2}$;\\[5pt]
	or, if $n=5$,\\[5pt]
	$k_4=\frac{ \normacomp(f'\wedge f''\wedge f''') \,((f'\times f''\times f'''\times f'''') \cdot f''''')}{\nu \normacomp(f'\wedge
		f''\wedge f'''\wedge f'''')^2}$;\\[2pt]
	\task  \ \ \ \ \ \ \ $\kapa 1'=
	\frac{\nu^3 f''\cdot f'''-\nu' (2 \normacomp(f'\wedge f'')^2+\nu^3 \nu'')}{\nu^4 \normacomp(f'\wedge f'')}$;\\[2pt]
	\task*(2)  \ \ \ \ \ \ \ $\V3=\frac{\nu (f''' \nu^2 \kapa 1-\nu f'' (3 \kapa 1 \nu'+\nu \kapa 1')+f' (-\nu \kapa 1 \nu''+\nu \nu' \kapa 1'+3 \kapa 1 \nu'^2+\nu^4\kapa 1^3))}{\normacomp(f'\wedge f''\wedge f''')}$.
\end{tasks}
\def\and{\quad{\rm and}\quad}
\ecoro
\proof From the definition of $\V1$, it follows (i) and, thus,
$$f''=\V 1 \nu'+\nu^2 \kapa 1 \V 2\quad{\rm and}\quad f'''= \left(\nu''-\nu^3 \kapa 1^2\right)\V 1+\left(3 \nu \kapa 1 \nu'+\nu^2 \kapa 1'\right)\V 2 +\nu^3 \kapa 1 \kapa 2 \V 3,$$
whence, we deduce (iii) and
$$f'' \cdot f'''=\nu' \left(\nu''-\nu^3 \kapa 1^2\right)+\nu^2 \kapa 1 \left(3 \nu \kapa 1 \nu'+\nu^2 \kapa 1'\right)=2 \nu^3 \kapa 1^2 \nu'+\nu^4 \kapa 1 \kapa 1'+\nu' \nu'',$$
that gives us 
$$\kapa 1'= \frac{-2 \nu^3 \kapa 1^2 \nu'+f'' \cdot f'''-\nu' \nu''}{\nu^4 \kapa 1},$$
hence (vii),
and the triangular linear system
$$ \left(
\begin{array}{ccc}
\nu & 0 & 0 \\
\nu' & \nu^2 \kapa 1 & 0 \\
\nu''-\nu^3 \kapa 1^2 & \kapa 1' \nu^2+3 \kapa 1 \nu' \nu & \nu^3 \kapa 1 \kapa 2 \\
\end{array}
\right)\left(
\begin{array}{c}
V_1 \\
V_2 \\
V_3 \\
\end{array}
\right)=\left(
\begin{array}{c}
f' \\
f'' \\
f''' \\
\end{array}
\right),$$
whose solution yields, in particular,
$$\arraycolsep=0pt\def\arraystretch{1.2}
\begin{array}{rcl}
\V3&\;=\;&\**\frac{f''' \nu^2 \kapa 1-\nu f'' (3 \kapa 1 \nu'+\nu \kapa 1')+f' (-\nu \kapa 1 \nu''+\nu \nu' \kapa 1'+3 \kapa 1 \nu'^2+\nu^4\kapa 1^3)}{\nu^5\kapa 1^2 \kapa 2}\\[10pt]
&\;=\;&\**\frac{\nu (f''' \nu^2 \kapa 1-\nu f'' (3 \kapa 1 \nu'+\nu \kapa 1')+f' (-\nu \kapa 1 \nu''+\nu \nu' \kapa 1'+3 \kapa 1 \nu'^2+\nu^4\kapa 1^3))}{\normacomp(f'\wedge f''\wedge f''')},
\end{array}$$
as we wanted.\qed
\brem Patiently replacing (or using a symbolic language software for to to this) the formulas (ii) and (vii) given above in (viii), we get a formula  in an arbitrary parameter
for $\V3$\footnote{\tiny\  $\arraycolsep=0pt\def\arraystretch{2}
	\begin{array}[t]{rcl}
	\V3&\;=\;&\left(-\frac{\nu'' \normacomp(f'\wedge f'')}{\nu \normacomp(f'\wedge f''\wedge f''')}+\frac{\nu'^2 \normacomp(f'\wedge f'')}{\nu^2 \normacomp(f'\wedge f''\wedge f''')}+\frac{\nu \nu' f''\cdot f'''}{\normacomp(f'\wedge f''\wedge f''') \normacomp(f'\wedge f'')}-\frac{\nu \nu'^2 \nu''}{\normacomp(f'\wedge f''\wedge f''') \normacomp(f'\wedge f'')}+\frac{\normacomp(f'\wedge f'')^3}{\nu^4 \normacomp(f'\wedge f''\wedge f''')}\right)f'+{}\\[2pt]
	&&\quad\quad\quad{}+ \left(-\frac{\nu' \normacomp(f'\wedge f'')}{\nu \normacomp(f'\wedge f''\wedge f''')}+\frac{\nu^2 \nu' \nu''}{\normacomp(f'\wedge f'') \normacomp(f'\wedge f''\wedge f''')}-\frac{\nu^2 f''\cdot f'''}{\normacomp(f'\wedge f'') \normacomp(f'\wedge f''\wedge f''')}\right)f''+\frac{ \normacomp(f'\wedge f'')}{\normacomp(f'\wedge f''\wedge f''')}f'''.\\
	\end{array}$}.
\erem
\bcoro If $\nu=1$, that is, $t$ is the arc length parameter, then\vskip10pt
\begin{tasks}[style=myenumerate,label-format={\normalfont},label-width={5mm}, item-indent={20pt}, item-format={\sl},label-align=right, after-item-skip={4pt}, column-sep={80pt}](3)
	\task  \ \ \ \ \ \ \ $\V 1=f'$;\\[2pt]
	\task \ \ \ \ \ \ \  $\kapa1=\normacomp(f'')$;\\[2pt]
	\task \ \ \ \ \ \ \  $\V2= \frac{ f''}{\normacomp(f'')}$;\\[2pt]
	\task \ \ \ \ \ \ \  $\kapa 2=\frac{\normacomp(f'\wedge f''\wedge f''')}{\normacomp(f'')^2}$;\\[2pt]
	\task*(2) \ \ \ \ \ \ \  $k_3= \frac{ \normacomp( f'')\,  \normacomp(f'\wedge f''\wedge f'''\wedge f'''')}{\normacomp(f'\wedge f''\wedge f''')^2}$;\\[2pt]
	\task*(2) \ \ \ \ \ \ \  $k_4= \frac{\normacomp(f'\wedge f'') \, \normacomp(f'\wedge f''\wedge f'''\wedge f''''\wedge f''''')}{\normacomp(f'\wedge f''\wedge f'''\wedge f'''')^2}$\\[5pt]
	or, if $n=5$,\\[5pt]
	$k_4=\frac{ \normacomp(f'\wedge f''\wedge f''')\, ((f'\times f''\times f'''\times f'''') \cdot f''''')}{ \normacomp(f'\wedge
		f''\wedge f'''\wedge f'''')^2}$;\\[2pt]
	\task  \ \ \ \ \ \ \ $\kapa 1'=\frac{f''\cdot f'''}{\normacomp(f'\wedge f'')}$;\\[2pt]
	\task*(3) \ \ \ \ \ \ \  $V_3=\frac{\normacomp(f'')^3}{\normacomp(f'\wedge f''\wedge f''')}f' -\frac{ f''\cdot f'''}{\normacomp(f'\wedge f''\wedge f''')  \normacomp(f'')}f''+\frac{ \normacomp(f'')}{\normacomp(f'\wedge f''\wedge f''')}f'''$.
\end{tasks}
\mbox{}\\
In particular, for $n=5$, we get more attractive formulas for the results of~\cite{Yilmaz5}.
\ecoro
\proof If follows immediately from the anterior corollary by using  $\norma{f'\wedge f''}=\norma{f''}$, which is true, since  $f'$ and $f''$ are orthogonal vectors.\qed 

\bcoro {\rm (\cite{Oneilgd}--4.3 Theorem)} If \funcao f{I}{\R 3} is 2-regular, then
\vskip10pt
\begin{tasks}[style=myenumerate,label-format={\normalfont},label-width={5mm}, item-indent={20pt}, item-format={\sl},label-align=right, after-item-skip={4pt}, column-sep={45pt}](3)
	\task  \ \ \ \ \ \ \  $\Tb=\V1=\frac{f'}{\nu}$;
	\task \ \ \ \ \ \ \   $\kappa=\kapa1=\frac{\normacomp(f'\wedge f'')}{\nu^3}$;
	\task \ \ \ \ \ \ \   $\tau=\kapa2=\frac{(f'\wedge f'')\cdot f'''}{\normacomp(f'\wedge f'')^2}$;
	\task \ \ \ \ \ \ \   $\Bb=\V3=\frac{f'\times f''}{\normacomp(f'\times f'')}$;
	\task  \ \ \ \ \ \ \  $\Nb=\V2=\Bb\times \Tb$.
\end{tasks}
\ecoro
\proof The proof is very simple. We just observe that (v) follows from the positive orientation of the Frenet frame $\{\Tb,\Nb,\Bb\}$.\qed

\bexe \label{exe1}Consider $f(t)=(t,t^2,t^3,t^4)$, $t\in\R{}$. A direct computation givens the matrix of the derivatives of $f$:
$$(f'\;f''\;f'''\;f'''')=\left(
\begin{array}{cccc}
1 & 0 & 0 & 0 \\
2 t & 2 & 0 & 0 \\
3 t^2 & 6 t & 6 & 0 \\
4 t^3 & 12 t^2 & 24 t & 24 \\
\end{array}\right)$$
whose determinant is 288.
Hence $f$ is 4-regular. The basic objects for the computation of the Frenet apparatus will be obtained below.\vskip10pt
\begin{tasks}[style=myenumerateexe,label-format={\normalfont\bf},label-width={10mm}, item-indent={25pt}, item-format={\sl},label-align=right, after-item-skip={10pt}, column-sep={0pt}](2)
	\task \ \ \ \ \ \ \ $\nu=\norma {f'}=\sqrt{1+4 t^2+9 t^4+16 t^6}$;
	\task \ \ \ \ \ \ \ $\nu'
	=\frac{2 t \left(2+9 t^2+24 t^4\right)}{\sqrt{1+4 t^2+9 t^4+16 t^6}}$;
	\task*(2) \ \ \ \ \ \ \
	$\arraycolsep=0pt\def\arraystretch{1.2}
	\begin{array}[t]{rcl}
	\norma{f'\wedge f''}^2&\;=\;&\det\left(
	\begin{array}{cc}
	f' \cdot f'\; &\; f' \cdot f'' \\
	f'' \cdot f' \;&\; f'' \cdot f''
	\end{array}
	\right)\\[20pt]
	&\;=\;&\arraycolsep=5pt\det
	\left(
	\begin{array}{lr}
	1+4 t^2+9 t^4+16 t^6 & 2 t \left(2+9 t^2+24 t^4\right) \\
	2 t \left(2+9 t^2+24 t^4\right) & 4 \left(1+9 t^2+36 t^4\right) \\
	\end{array}
	\right)\\[20pt]
	&\;=\;&4 (1+9 t^2+45 t^4+64 t^6+36 t^8);
	\end{array}$
	\task*(2) \ \ \ \ \ \ \ $ f'\times f''\times f'''=\det\left (
	\begin{array}{cccc}
	1 & 0 & 0 & e_1 \\
	2 t & 2 & 0 & e_2 \\
	3 t^2 & 6 t & 6 & e_3 \\
	4 t^3 & 12 t^2 & 24 t & e_4 \\
	\end{array}
	\right)=(-48 t^3,72 t^2,-48 t,12)$;
	\task \ \ \ \ \ \ \ $( f'\times f''\times f''')\cdot f''''=12\cdot 24=288$;
	\task \ \ \ \ \ \ \ $\normacomp(f'\wedge f''\wedge f''')^2=144 (1+16 t^2+36 t^4+16 t^6)$;
\end{tasks}
\vskip10pt
Now, the  Frenet apparatus:\vskip10pt
\begin{itemize}
	\item[\textbf{a)}]   $\kapa 1=\frac{\normacomp(f'\wedge f'')}{\nu^3}=\frac{2 \sqrt{1+9 t^2+45 t^4+64 t^6+36 t^8}}{\left(1+4 t^2+9 t^4+16
		t^6\right)^{3/2}}$;	
	\item[\textbf{b)}]  $ \kapa 2=\frac{\normacomp(f'\wedge f''\wedge f''')}{\normacomp(f'\wedge f'')^2}=\frac{3 \sqrt{1+16 t^2+36 t^4+16 t^6}}{1+9 t^2+45 t^4+64 t^6+36 t^8}$;
	\item[\textbf{c)}]  $\kapa 3=\frac{ \normacomp(f'\wedge f'') (f'\wedge f''\wedge f''') \cdot f''''}{\nu \normacomp(f'\wedge f''\wedge f''')^2}=\frac{4 \sqrt{1+9 t^2+45 t^4+64 t^6+36 t^8}}{\sqrt{1+4 t^2+9 t^4+16 t^6} \left(1+16
		t^2+36 t^4+16 t^6\right)}$;
	\item[\textbf{d)}]  $\arraycolsep=0pt\V 1=\frac{f'}{\nu}=\frac{1}{\sqrt{1+4 t^2+9 t^4+16 t^6}}
	( 1,  2 t,  3 t^2 ,  4 t^3)$;
	\item[\textbf{e)}]  $\arraycolsep=0pt\V 2
	=\frac{1}{\sqrt{1+4 t^2+9 t^4+16 t^6} \sqrt{1+9 t^2+45 t^4+64 t^6+36 t^8}}\left(
	\begin{array}{c}
	-t \left(2+9 t^2+24 t^4\right) \\
	1-9 t^4-32 t^6 \\
	3 t+6 t^3-24 t^7 \\
	2 t^2 \left(3+8 t^2+9 t^4\right) \\
	\end{array}
	\right)$;
	\item[\textbf{f)}]  $\arraycolsep=2pt \V4=\frac{f'\times f''\times f'''}{\normacomp(f'\times f''\times f''')}=
	\frac{1}{12\sqrt{1+16 t^2+36 t^4+16 t^6}} \left(
	\begin{array}{cccc}
	1 & 0 & 0 & e_1 \\
	2 t & 2 & 0 & e_2 \\
	3 t^2 & 6 t & 6 &e_3 \\
	4 t^3 & 12 t^2 & 24 t &e_4\\
	\end{array}\right)=
	\frac{1}{\sqrt{1+16 t^2+36 t^4+16 t^6}} \left(
	\begin{array}{c}
	-4 t^3\\ 6 t^2\\ -4 t\\ 1 
	\end{array}
	\right)$.
\end{itemize}
\eexe
\noindent In the calculations above, mainly if done by hand,  that of $\V2$ requires a little more work, but by using the indicated formula (see~\ref{coro5}-\ref{coro5i3}, for other options), it follows nicely. It remains to calculate $\V 3$. For this, we use the positive orientation of  the Frenet frame $\{\V1,\V2,\V3,\V4\}$. We have,
$$\V3=-\V1\times\V2\times\V4,$$
since $\det(\V1,\V2,\V4,\V3)=-\det(\V1,\V2,\V3,\V4)=-1$. We omit the explicit result. In particular, the Frenet apparatus at $t=0$ is given by
$$\cal A(0)=\{{2,3,4,(1,0,0,0),(0,1,0,0),(0,0,1,0),(0,0,0,1)}\}.$$

\section{On the Fundamental Theorem for Curves}
In~\cite{Kuhnel} and~\cite{Rolfpp},  we found the following theorem, known as {\sl Fundamental Theorem of the Local Theory of Curves}. Its statement could be posed in this way:\vskip5pt
$$\parbox{4.25in}{{\bf Theorem.} \sl  Let $k_i(s)$, $i = 1,\ldots  n-1$, $s\in  I$, be smooth functions with $k_i(s)>0$, $i = 1,\ldots  n-2$ . For a fixed parameter $s_0\in I$, suppose we have been given a point $q_0\in\R n$ as well as an $n$ frame $e_1(0)$, \ldots, $e_n(0)$.  Then there is a unique curve \funcao c {I}{\R n} parametrized by arc length and satisfying the following three conditions:\vskip-5pt
	\benum
	\item[1.] $c(s_0)=q_0$,
	\item[2.]  $e_1(0)$, \ldots, $e_n(0)$ is the Frenet frame of $c$ at $q_0$,
	\item[3.] $k_i(s)$, $i = 1,\ldots  n-1$, are the curvatures of $c$.
	\eenum}$$
A similar theorem appears in~\cite{Manfra}, on tridimensional case. The proofs,  in both cases, use the general existence and uniqueness theorem for systems of linear differential equations. In that follows, we extend, the cited theorem to arbitrary speed curves in $\R  n$. A complete proof is presented. Our theorem is as follows.
\bteo \label{fundteo}Let $\nu(t)$ and $k_j(t)$, $j \in\{1,2,\ldots  n-1\}$, $t\in  I$, be smooth functions such that, for all $t\in I$, 
$\nu(t)>0$ and $k_j(t)>0$, $j\leq n-2$. Then there exists \funcao{f}{I}{\R n} with speed $\nu(t)$ and curvatures  $k_j(t)$, $j \in\{1,2,\ldots  n-1\}$. Furthermore, if  \funcao{g}{I}{\R n} is another curve satisfying these conditions, there exists an orientation preserving isometry
\funcao F{\R n}{\R n} such that $g=F\circ f$, that is, $f$ and $g$ are congruent curves.
\eteo

Before the proof of this theorem, we obtain some preliminary results. Initially, since we are going to work with two curves  $f$ and $g$, we will make a suitable adjustment in the notation: we will put a  tilde over each object associated to the curve $g$.
For example, $\nug$  and $\kapag 1$ indicate the speed and  the first curvature of $g$ whereas $\Vg 1$ indicates the first vector field of the Frenet frame of  $g$. Thus, given the curves $f$ and $g$, both defined in $I$, we have the Frenet apparatus: 
$$\cal A_f=\{\kapasf,\reffrenetf\}$$
and 
$$\cal A_g=\{\kapasg,\reffrenetg\}.$$
Now, consider \funcao f{I}{\R n} a parametrized curve in $\R n$ e let \funcao{F}{\R n}{\R n} be an orientation preserving isometry, that is, $F(X)=S X+X_0$, where $S$ is a orthogonal matrix of  determinant~1. Using~$F$, we construct a new curve, namely, $g=F\circ f$, or
$g(t)=F(f(t))$, $t\in I$. The next result establishes that in a certain way $g$ inherits the Frenet apparatus of $f$.
\bprop The speed of $g$ equals speed of $f$ and the Frenet apparatus of $g=F\circ f$ is given by
$$\cal A_g=\{\kapasf f,\reffrenetfs\}$$
In other words, $\kapag j=\kapa j$, $1\le j\le n-1$,  and $\Vg j=S \V j$, $1\le j\le n$.
\eprop
\proof Using the chain rule, we get $g'(t)=\dF_{f(t)}(f'(t))=S f'(t)$, Hence $\nug^2=(S f')\cdot(S f')=\nu^2$, because $S$ preserves the  inner  product, which proof the claim on the speeds. Firstly, we  study the curvatures. Using  Corollary~\ref{coro4} together the norm of a multivector given in~(\ref{pvetornorma}), we see that given $\puplasp v j\in\R n$, the following holds
$$\norma{\pvetor {S v}j}=\sqrt{\det ((S v_i)\cdot(S v_j))}=\sqrt{\det (v_i\cdot v_j)}=\norma{\pvetor vj},$$
again because $S$ preserves the  inner product. Thus, for $1\le m\le n-2$,
$$\arraycolsep=0pt\def\arraystretch{3}
\begin{array}{rcl}
\**\kapag m&\;=\;&\**\frac{\norma{\pvetorf g {m-1}}\, \norma{\pvetorf g {m+1}}}{\nug \norma{\pvetorf g {m}}^2}\\
&\;=\;&\** \frac{\norma{\pvetorf{S f}{m-1}}\, \norma{\pvetorf{S f}{m+1}}}{\nu \norma{\pvetorf{S f}{m}}^2}\\
&\;=\;&\** \frac{\norma{\pvetorf f {m-1}}\, \norma{\pvetorf f {m+1}}}{\nu \norma{\pvetorf f {m}}^2}=\kapa m,
\end{array}$$
since $g^{(j)}=S f^{(j)}$, for all $j$, by the chain rule. It remains to see the claim about the Frenet frames and the last curvature. From the  definitions, we get  
$$\Vg 1=\frac{g'}{\nug}=\frac{S\,f'}{\nu}=S\,\frac{ f'}\nu=S \V1.$$ 
By induction, assuming that $\Vg j=S \V j$, $1\leq j\le n-1$,  we show that  $\Vg{j+1}=S \V{j+1}$. In fact,
first 
$$\Vg{j}'=(S \V j)'=S(-\nu \kapa {j-1}\V{j-1}+\nu \kapa{j} \V{j+1})=
-\nu \kapa {j-1}\Vg{j-1}+\nu \kapa{j} S \V{j+1}.$$
On the other hand,
$$\Vg{j}'=-\nug \kapag {j-1}\Vg{j-1}+\nug \kapag{j} \Vg{j+1}=-\nu \kapa {j-1}\Vg{j-1}+\nu \kapa{j}\Vg{j+1}.$$
Hence $\nu \kapa{j} S \V{j+1}=\nu \kapa{j}\Vg{j+1}$ and thus $\Vg{j+1}=S \V{j+1}$.  Note that the preceding arguments could be used for the curvatures. We will make it so for $\kapag {n-1}$. Since that $\{\reffrenetf\}$ is a positively oriented frame and $\det S=1$, it is not hard to see that 
$$S \V n=S(\crossvetor V{n-1})=\crossvetor{S\, V}{n-1}=\crossvetor{\til V}{n-1}=\Vg n.$$
Indeed,  it is sufficient to note that $S(\crossvetor V{n-1})\cdot \Vg j=0$, for all $1\le j\le n-1$.  Differentiating  $\Vg n=S  \V n$ yields $\nug \kapag {n-1}\Vg{n-1}=\nu \kapa {n-1} S \V{n-1},$
from which it follows that $\kapag {n-1}=\kapa {n-1}$ and the proof is complete.\qed
\brem The proof above shows that if $F$ is an  orientation reversing isometry, then all works well, except that $\Vg n=-S \V n$ and $\kapag {n-1}=-\kapa {n-1}$.
\erem
The converse of the proposition above
is true. Its statement  is as follows. In it $f$ is as before, having speed $\nu$ and Frenet apparatus
$$\cal A_f=\{\kapasf,\reffrenetf\}.$$

\bprop \label{propo4}Let \funcao{h}{I}{\R n} be a parametrized curve with speed $\nuh$ and Frenet apparatus
$$\cal A_h=\{\kapash,\reffreneth\}.$$
If  $\nuh=\nu$  and  $\kapah j=\kapa j$, $1\le j\le n-1$, then there exists an 
preserving orientation isometry $F$ of $\R n$ such that $h=F\circ f$.
\eprop
\proof For simplicity, suppose that $0\in I$ and  $h(0)=f(0)=(0,0,\ldots, 0)$. Now, let $S$ be the orthogonal transformation that sends the Frenet frame $\cal F_f(0)$ to the  Frenet frame $\cal F_h(0)$, that is, $S \V j(0)=\Vh j(0)$, $j=1,2,\ldots n$. Of course that $\det S=1$, for these frames are positively oriented. In that follows, we use the ideas of O'Neill~\cite{Oneilgd}, in his proof for $n=3$. Consider $g=S\circ f$ and the real function 
$$r(t)=\sum_{j=1}^n \Vg j(t)\cdot\Vh j(t),\quad t\in I,$$
where, as before, $\cal F_g=\{\Vg j(t)=S \V j(t),\quad 1\le j\le n\}$ is the Frenet frame field of $g$. Of course that $g(0)=h(0)$, $\cal F_g(0)=\cal F_h(0)$, $\nug=\nuh=\nu$ and $\kapah j=\kapag j=\kapa j$, for $1\le j\le n-1$. The main idea now is to show that $g$ matches $h$. We start observing that $r(0)=n$, since $\cal F_g(0)=\cal F_h(0)$. 
We have that
$$r'=\sum_{j=1}^n( \Vg j'\cdot\Vh j+ \Vg j\cdot\Vh j').$$
Claim: $r'=0$. In fact, firstly, the first and second summands of $r'$ are
$$\Vg 1'\cdot\Vh 1+ \Vg 1\cdot\Vh 1'=\nu  \kapa 1   \Vg 2\cdot \Vh 1+\nu  \kapa 1   \Vg 1\cdot \Vh 2$$
and
$$\Vg 2'\cdot\Vh 2+ \Vg 2\cdot\Vh 2'=-\nu \kapa 1   \Vg 2\cdot\Vh 1-\nu \kapa 1   \Vg 1\cdot \Vh 2+\nu \kapa 2   \Vg 3\cdot \Vh 2+\nu \kapa 2   \Vg 2\cdot \Vh 3.$$
Hence
$$\sum_{j=1}^2( \Vg j'\cdot\Vh j+ \Vg j\cdot\Vh j')=\nu \kapa 2   \Vg 3 \cdot\Vh 2+\nu \kapa 2   \Vg 2 \cdot\Vh 3.$$
By induction on $m$, $m\le n-1$, we get
$$\sum_{j=1}^m( \Vg j'\cdot\Vh j+ \Vg j\cdot\Vh j')=\nu \kapa m   \Vg {m+1} \cdot\Vh m+\nu \kapa m   \Vg m \cdot\Vh {m+1}.$$
Thus
$$\arraycolsep=0pt\def\arraystretch{1.2}
\begin{array}{rcl}
\**r' &\;=\;&\**\nu \kapa {n-1}   \Vg {n} \cdot\Vh { n-1}+\nu \kapa { n-1}   \Vg { n-1} \cdot\Vh {n}+ (\Vg n'\cdot\Vh n+ \Vg n\cdot\Vh n')\\[10pt]
\**&\;=\;&\**\nu \kapa {n-1}   \Vg {n} \cdot\Vh { n-1}+\nu \kapa { n-1}   \Vg { n-1} \cdot\Vh {n}-\nu \kapa {n-1}   \Vg {n} \cdot\Vh {n-1}-\nu \kapa {n-1}   \Vg {n-1} \cdot\Vh {n}=0
\end{array}$$
which proofs the claim. 
So $r$ is a constant function. From $r(0)=n$ we get $r(t)=n$, $t\in I$. Since each summand of $r$ is at most 1, we get that all of them are equal to 1.  In particular, $\Vg 1=\Vh 1$ which implies that $g'=h'$ and, by integration on $[0,t]$, $g=h$, that is, $h=S\circ f$. The general case is obtained by considering $h(t)-h(0)$ and $f(t)-f(0)$ instead of $h$ and $f$, respectively. From this we conclude that 
$S(f(t)-f(0))=h(t)-h(0)$. Hence $g=F\circ h$, where $F(X)=S X+X_0$, $X_0=h(0)-S(f(0))$.\qed
\vskip10pt\def\deumam#1#2#3{#2\le #1\le#3}

We are almost ready for proving the fundamental theorem (Theorem~\ref{fundteo}). Before, let us summarize  some facts on vector matrices.

\bfact \rm\label{factum}Given the vectors $\puplasp V m$ in $\R n$, we indicate by $$V=\vetcol V m=\vetlinha Vm,$$ where $\trans M$ indicates the transpose of $M$, the $n\times 1$ vector matrix with elements $\V j$. $V$ is called a $(m\times 1)$-{\sl vector matrix}. Note that by stacking all coordinates of the vectors $\V j$ (viewed as column vectors), we obtain an $(mn)\times 1$ real matrix or an $mn$-column vector. 

Given a  $p \times q$ real matrix $A$ and a $(q\times 1)$-vector matrix $V=\vetlinha Vq$, the  {\sl product} $W=A\, V$ is defined to be the $(q\times 1)$-vector column matrix $W=\vetlinha Wq$ such that
$$W_i=\sum_{j=1}^m a_{ij}V_j,\; a_{ij}\in \R{},$$
or
$$\vetcol W p=\arraycolsep=5pt\left(
\begin{array}{cccc}
a_{11} & a_{12} & \ldots & a_q \\
a_{21} & a_{22} & \ldots & a_{2 q} \\
\vdots & \vdots & \ddots & \vdots \\
a_{\text{p1}} & a_{p2} & a_{p3} & a_{\text{pq}} \\
\end{array}
\right)\vetcol V q.$$
Of course, if $V=\puplalinha Vp$ is a $(1\times p)$-vector matrix, the product $VA$ is a well defined $(1\times q)$-vector matrix. The following properties hold true. In them, $a\in \R{}$, $C$ is a $r\times p$ real matrix, $A$ and $B$ are $p\times q$ real matrices and $V$ and $W$ are $(q\times 1)$-vector matrices.
\\[5pt]
\begin{tasks}[style=myenumerate,label-format={\normalfont},label-width={5mm}, item-indent={20pt}, item-format={\sl},label-align=right, after-item-skip={4pt}, column-sep={40pt}](2)
	\task $A\,(a V)=(aA)\,V=a (A\,V)$;
	\task $(A+B)\,V=A\,V+B\,V$;
	\task $A(V+W)=A\,V+A\,W$;
	\task $(CA)V=C(A\,V)$.
\end{tasks}
\vskip10pt

Now, consider a linear system of the type $AV=W$, where $A$, $V$ and $W$ are as above. From this,  an interesting exercise arrives: find a usual linear system equivalent to it. It is easy. Just stack the coordinates of the elements of the matrices $V$ and $W$, obtaining $\til V$ and $\til W$ in $\R{qn}$ and replace $A$ by $\til A$, of order $pn\times qn$, where $\til a_{ij}$ is the block $\til a_{ij}=a_{ij}{\rm Id}_n$ and ${\rm Id}_n$ is the $n\times n$ identity matrix. In the Table~1 below, we see an example for $p=2$, $q=3$ e $n=4$, where $\V i=(\V{i1},\V{i2},\V{i3},\V{i4})\in\R4$.

Another useful product of vector matrices is what uses the inner product in its construction. We will deal  only with vector matrices of order either $n\times 1$ or $1\times n$. For this, let $V=\vetlinha V n$ and $W=\puplalinha W n$. The {\sl dot  product} of $V$ by $W$ is the $n\times n$ real matrix 
$$V\cdot W=\vetcol Vn\cdot\puplalinha W n=\arraycolsep=5pt\left(
\begin{array}{cccc}
V_1 \cdot W_1 & V_1 \cdot W_2 & \ldots & V_1 \cdot W_n \\
V_2 \cdot W_1 & V_2 \cdot W_2 & \ldots & V_2 \cdot W_n \\
\vdots & \vdots & \ddots & \vdots \\
V_n \cdot W_1 & V_n \cdot W_2 & \ldots & V_n \cdot W_n \\
\end{array}
\right).$$
\brem 
With this new notation, the inner product in~(\ref{pipvetor}) becomes  
$$(\pvetor vn)\cdot(\pvetor wn)=\det (v\cdot w).$$
\erem
\fbox{\vbox{{\scriptsize
			$$\arraycolsep=2pt
			\begin{array}{rcl}
			\AA\Vs&=&\Mav\\[30pt]
			&\;=\;&\Mbidexp\VV\\[30pt]
			&\;=\;&\Mbid\VV=\WW
			\end{array}$$}
		\mbox{}\\[5pt]
		\hbox to\textwidth{\bf \hss {\sc Table 1}: Expanding a vector linear system to an usual linear system\hss}}
}
\\[10pt]

Given $n\times n$ real matrices $A$ and $B$, a $(n\times 1)$-vector matrix $V$ and a $(1\times n)$-vector matrix $W$, the following  hold:
\begin{tasks}[style=myenumerate,label-format={\normalfont},label-width={5mm}, item-indent={20pt}, item-format={\sl},label-align=right, after-item-skip={4pt}, column-sep={40pt}](2)
	\task $A\,(V\cdot W)=(A\,V)\cdot W$;
	\task $(V\cdot W)A=V\cdot (WA)$.
\end{tasks}
Of course, we have distributive properties for appropriate choices of the vector matrices. 
Moreover, if $V$ e $W$ depend differentiably on $t\in I$, then $$(V\cdot W)'=V'\cdot W+V\cdot W'.$$
\efact
\vskip -10pt 
Finally, using the facts above, we will prove Theorem~\ref{fundteo}. The proof involves four steps, namely:\\
\begin{tasks}[style=myenumerate,label-format={\normalfont},label-width={15mm}, item-indent={35pt}, item-format={\rm},label-align=right, after-item-skip={10pt}, column-sep={40pt}](1)
	\task[(Step 1)] \hspace{1cm} From the given functions $\nu$ and  $\kapa j$, $\jdeuman j1{n-1}$, we construct, based on the Frenet equations, a system of $n^2$ first order linear differential equations, which we refer as $(\frsys)$.
	\task[(Step 2)] \hspace{1cm}  We apply the general existence and uniqueness theorem for systems of linear differential equations to $(\frsys)$ and take the solution that satisfies a certain initial condition. Such  a solution is an $n$-vector column matrix $V=\vetlinha Vn$ that depends on $t$.
	\task[(Step 3)] \hspace{1cm} The existence assertion of the theorem: we verify that $\{\V1(t),\V2(t),\ldots,\V n(t)\}$ is actually a orthonormal frame field and, from the vector function $\V 1(t)$,  we construct a curve $f$ with speed $\nu$, curvatures $\kapa j$, $\jdeuman j1{n-1}$, and Frenet frame field $\cal F_f=\{\reffrenetf\}$.
	\task[(Step 4)] \hspace{1cm} The uniqueness assertion of the theorem:  given a curve $g$ with speed $\nu$ and curvatures $\kapa j$, $\jdeuman j1{n-1}$ there exists an isometry of $\R n$ such that $g=F\circ f$. This step we have already seen in the Proposition~\ref{propo4}.
\end{tasks}
We start writing the Frenet equations~(\ref{freneteq}) of a given $f$ in a matrix form (see~(\ref{freneteq5}), for $n=5$). They become
$$\arraycolsep=2pt\left(
\begin{array}{c}
\V1'  \\
\V2' \\
\V3'  \\
\V4' \\
\vdots\\
\V{n-2}' \\
\V{n-1}' \\
\V{n}' \\
\end{array}
\right)=\left(
\begin{array}{ccccccc}
0 & \kapatruque 1  (\nu) & 0 & 0 & \ldots & 0 & 0 \\
\kapatruque 1  (-\nu) & 0 & \kapatruque 2  (\nu) & 0 & \ldots & 0 & 0 \\
0 & \kapatruque 2  (-\nu) & 0 & \kapatruque 3  (\nu) & \ldots & 0 & 0 \\
0 & 0 & \kapatruque 3  (-\nu) & 0 & \ldots  & 0 & 0 \\
\vdots & \vdots &\vdots & \ddots & \vdots &\vdots & \vdots \\
0 &0 & 0 & \ldots &0 & \kapatruque {n-2}  (\nu) & 0 \\
0 & 0 & 0 &\ldots & \kapatruque {n-2} (-\nu) & 0 & \kapatruque {n-1} (\nu) \\
0 & 0 & 0 & \ldots & 0 & \kapatruque {n-1} (-\nu) & 0 \\
\end{array}
\right)
\left(
\begin{array}{c}
\V1  \\
\V2 \\
\V3  \\
\V4  \\
\vdots\\
\V{n-2} \\
\V{n-1} \\
\V{n} \\
\end{array}
\right),$$
or simply
$$V'=\Mfr\,  V,\leqno{(\frsys)} $$
where $\V 1=\frac{f'}{\nu}$, $V=\vetlinha V n$ and $V'=\vetlinha {V'} n$ are $n$-vector column matrices constructed from the Frenet frame of $f$, and $\Mfr$, called the {\sl Frenet matrix} of $f$, is the  $n\times n$ skew symmetric matrix such that 
$$\arraycolsep=0pt \Mfr_{ij}=\V i'\cdot\V j=V j\cdot(-\nu \kapa{i-1}\V{i-1}+\nu \kapa{i}\V {i+1})=\left\{
\begin{array}{ll}
-\nu \kapa{i-1},&\quad j=i-1\\
\nu \kapa{i},&\quad j=i+1\\
0,&\quad j\not\in\{i-i,i+1\}.
\end{array}\right.$$
Of course even without knowing $f$, we can consider $(\frsys)$ on the interval $I\ni 0$, because we can construct~$\Mfr$. Hence, we get a first order system of differential equations for the unknown vector functions $\reffrenetf$, which, according to the Fact~\ref{factum}, can be viewed as a usual system of $n^2$ first order linear differential equations. Thus we can apply the general existence and uniqueness theorem for systems of linear differential equations to it and achieve a unique set of functions $\{\reffrenetfcomt\}$ that satisfies $(\frsys)$ and the initial conditions $\V j(0)=e_j$, $\jdeuman j1n$, where $\{\puplasp en\}$ is the canonical basis of $\R n$. In this way, we go through the steps 1 e 2. 

For the step 3, we consider the $n\times n$ matrix function $A(t)=(\V i(t)\cdot\V j(t))$, or
$$A(t)=\vetcolcomt Vn\cdot\puplalinhacomt Vn=V(t)\cdot \trans {V(t)}.$$
It is convenient to remark that $A(0)={\rm Id}_n$, where ${\rm Id}_n$ is the $n\times n$ identity matrix.
Differentiating~$A$, we get
$$\arraycolsep=0pt\def\arraystretch{1.5}
\begin{array}{rcl}
\**A' &\;=\;&\**V'\cdot\trans V+V\cdot \trans V'=(\Mfr\, V)\cdot \trans V+ V\cdot\trans{(\Mfr\,V)}\\
\**&\;=\;&\**\Mfr( V\cdot \trans V)+(V\cdot \trans V)(\trans \Mfr)=\Mfr\,A+A\,(\trans \Mfr)
\end{array}$$
Now, note that the expression $\Mfr\,A+A\,(\trans \Mfr)$ is linear as function of $A$. Thus, we have that $A$ is a solution of the linear matrix differential equation  $$X'=\Mfr\,X+\trans \Mfr\, X$$ with the initial condition $X(0)={\rm Id}_n$. By using vectorization (that is, by  stacking columns) of  matrices, it is not hard to show that this equation reduces to a system of $n^2$ first order linear differential equations (For $n=3$, see the Table 2 below). Hence $A$ is the unique solution of  $$X'=\Mfr\,X+\trans \Mfr\, X, \quad X(0)={\rm Id}_n.$$ This fact, together with skew symmetry of $\Mfr$, implies $A(t)={\rm Id}_n$, for all $t\in I$, since $X={\rm Id}_n$ 
satisfies trivially the equation. In fact, ${\rm Id}_n'=0$ and $\Mfr\,{\rm Id}_n+\trans \Mfr\, {\rm Id}_n=\Mfr+\trans \Mfr=0.$ We conclude that $\{\reffrenetfcomt\}$ is actually a orthonormal frame field.  In reality, it is a positively oriented orthonormal frame field, since it coincides with the canonical basis at $t=0$ (the $\det$  is a continuous function). It remains only to construct a curve $f$ for attaining the step 3.

In this point, we have a positively oriented orthonormal frame field $\{\reffrenetfcomt\}$, $t\in I$, that satisfies $(\frsys)$. Since $\V 1$ must be the unit tangent vector of $f$, there is a natural way to choose the curve 
$f$, namely,  $$f(t)=\int_0^t\nu(u)\V 1(u)\,{\rm d}u.$$ From this, we get $f'=\nu\V1$, or $\V1=f'/\nu$. For the moment, we suppose that $f$ is a $ (n-1)$-regular curve with speed $\nu_f$, curvatures $\kapaf j$, ${\jdeuman j1{n-1}}$, and Frenet frame $\{\Vf j,\jdeuman j1n\}$. Hence 
\begin{tasks}[style=myenumerate,label-format={\normalfont},label-width={35mm}, item-indent={20pt}, item-format={\sl},label-align=right, after-item-skip={5pt}, column-sep={30pt}](2)
	\task $\nu_f=\nu$;
	\task $\Vf 1=\V1$;
	\task $f''=\nu' \V 1 +\nu^2 \kapa 1  \V 2$;
	\task $f'\wedge f''=\nu^3 \kapa 1  \V 1\wedge \V 2$.
\end{tasks}
The Corolary~\ref{coro4}, together with (ii),  yields $\kapaf 1=\frac{\normacomp(f'\wedge f'')}{\nu_f^3}=\kapa 1$. Using this and (ii), it follows that $\Vf2=\V 2$. In fact,  the differentiation of (ii) gives $\kapaf 1\nu_f \Vf 2=\kapa 1 \nu\V2$.
Now, a direct  calculation shows that
$$f'\wedge f''\wedge f'''=\nu^6 \kapa 1^2 \kapa 2  \V 1\wedge \V 2\wedge \V 3=\nu_f^6 \kapaf 1^2 \kapa 2  \V 1\wedge \V 2\wedge \V 3$$
which, together with the Corolary~\ref{coro4} and the derivative of  $\Vf2=\V 2$, implies $\kapaf2=\kapa 2$ and $\Vf3=\V3$. By repeating this process inductively (as in the Theorem~\ref{wedgeteomain}), we conclude 
$$f' \wedge f'' \wedge\cdots \wedge  f^{(m)} =
\nu^{\frac{ m(m+1)}2} \kapa 1^{m-1} \kapa 2^{m-2}
\cdots \kapa {m-1} \V1 \wedge  \V2 \wedge\cdots \wedge \V m,\quad  2\le m\le n,$$
that, in particular,  shows the regularity of $f$,
$\kapaf j=\kapa j$, $\jdeuman j1{n-1}$, and $\Vf j=\V j$, for $\jdeuman j1n$. We are done.\qed\\[10pt]
\fbox{\vbox{{\scriptsize
			$$\arraycolsep=2pt
			\usualsode$$}
		\mbox{}\\[5pt]
		\hbox to\textwidth{\hss \sc Table 2\hss}\\[2pt]
		\hbox to\textwidth{\bf \hss  $X'=\Mfr\,X+\trans \Mfr\, X$ as an usual  system of  linear differential equations, $\mathbf n=3$\hss
		}\\[2pt]
		\hbox to\textwidth{\bf \hss  $\mathbf X=\trans{(1\;0\;0\;0\;1\;0\;0\;0\;1)}$ is the solution that corresponds to $\mathbf{\rm Id}_3$\hss}}
}
\\[10pt] 
\bexe We will consider the system $(\frsys)$ with $n=3$, $\nu=\sqrt 2$, $\kapa 1=\kapa 2=\frac1{\sqrt 2}$ and initial condition 
\beq\label{condini47} V(0)=\parl \parl 0,\frac{1}{\sqrt{2}},\frac{1}{\sqrt{2}}\parr,\parl -1,0,0\parr,\parl 0,-\frac{1}{\sqrt{2}},\frac{1}{\sqrt{2}}\parr\parr.\eeq
So, we have
$$\cal M\,V=V':\quad\quad\left(
\begin{array}{ccc}
0 & \frac{1}{\sqrt{2}} & 0 \\
-\frac{1}{\sqrt{2}} & 0 & \frac{1}{\sqrt{2}} \\
0 & -\frac{1}{\sqrt{2}} & 0 \\
\end{array}
\right) \left(
\begin{array}{c}
V_1 \\
V_2 \\
V_3 \\
\end{array}
\right)=
\left(
\begin{array}{c}
V_1' \\
V_2' \\
V_3'\\
\end{array}
\right),\leqno(3-\frsys)$$
subject to~\refeq{condini47}

\noindent
We will indicate two ways for solving $(3-\frsys)$: one for illustrating the conversion to an usual system of linear differential equations as in the Table~1 and the other going hand in hand with an algorithm that  plays a key rule in the classification of the curves of constant curvatures. Thus, at the end, we will get a curve with speed $\sqrt2 $, curvature and torsion equal to $\frac12$.\\[5pt]
\underline{Solution 1.} 
Let $\V i(t)=(V_{i1}(t),V_{i2}(t),V_{i3}(t))$, $i=1,2,3$. Hence $V=(\V1,\V2,\V3)$ is the unknown of the system, which converted to its usual form becomes
$$\til{\cal M}\,\til V=\til V':\quad\quad {\scriptsize\usual3}$$
with
$$\til V(0)=\parl 0,\frac{1}{\sqrt{2}},\frac{1}{\sqrt{2}},-1,0,0,0,-\frac{1}{\sqrt{2}},\frac{1}{\sqrt{2}}\parr.$$
Remember that the important piece of the solution $V$ is the vector $\V1=(V_{11},V_{12},V_{13})$, from which, by integration of $\nu \V1$, we construct the curve $f$.

It is well known fact that $$\til V(t)=(e^{t\til\Mfr})\,\til V(0)=\left(\sum_{j=0}^\infty \frac {t^j}{j!}\til\Mfr^j\right)\til V(0)$$ is the unique solution of~$(3-\frsys)$).  An elementary procedure to calculate this solution is indicated next. Patiently, only looking at with the first three coordinates of $\til V(t)$ and using induction on $j$, it is possible to obtain that $(V_{11},V_{12},V_{13})$ equals
$$\frac 1{\sqrt2}\left(-t+\frac{t^3}{3!}-\frac{t^5}{5!}+\cdots,\frac{(-1)^{m-1} t^{2 m+1}}{(2 m+1)!}+\cdots,1-\frac{t^2}{2!}+\frac{t^4}{4!}+\cdots,\frac{(-1)^m t^{2 m}}{(2 m+1)!}+\cdots,1\right).$$
Hence $$\V 1(t)=\frac 1{\sqrt2}(-\sin t,\cos t,1).$$
Now, by integrating $\nu\V1$,
$$\int_0^t\nu\V1(u)\,{\rm d}u=\int_0^t\sqrt 2\V1(u)\,{\rm d}u=\int_0^t(-\sin u,\cos u,1)\,{\rm d}u=(\cos t-1,\sin t,t),$$
which, after a translation, yields $f(t)=(\cos t,\sin t, t)$, the very well known circular helix of $\R3$, as expected.\\[5pt]
\underline{Solution 2.} Here, we come back to consider the original system $(3-\frsys)$: $\Mfr\, V=V'$, where \mbox{$V=(\V1,\V2,\V3)$}, together with 
$$V(0)=(\V1(0),\V2(0),\V3(0))=\parl \parl 0,\frac{1}{\sqrt{2}},\frac{1}{\sqrt{2}}\parr,\parl -1,0,0\parr,\parl 0,-\frac{1}{\sqrt{2}},\frac{1}{\sqrt{2}}\parr\parr.$$
Again, the solution is $V(t)=(e^{t\Mfr})\, V(0)$. We are going to calculate $V(t)$ and then the curve $f(t)$.\\
\begin{tasks}[style=myenumerate,label-format={\normalfont},label-width={15mm}, item-indent={35pt}, item-format={\rm},label-align=right, after-item-skip={4pt}, column-sep={40pt}](1)
	\task[(Step 1)] \hspace{1cm} To reduce $\Mfr$ to a simpler form, which is possible because it is a skew symmetric matrix. Actually, there exists a orthogonal matrix $Q$ such that $\cal N=\inv Q\Mfr Q$ is a block matrix of the kind
	$$\cal N=\left(
	\begin{array}{ccc}
	0 & -a & 0 \\
	a & 0 & 0 \\
	0 & 0 & 0 
	\end{array}
	\right)=\left(
	\begin{array}{cc}
	A & 0 \\
	0 & 0 
	\end{array}
	\right),$$
	for some $a\in \R{}$.
	This comes from the {\sl  normal-form of a skew symmetric matrix theorem}, which can be found in~\cite{Greub}.
	\task[(Step 2)] \hspace{1cm}  To calculate $e^{t\cal N}$, This is easy, because the powers of $\cal N$ reduce to those of the block $A$, whose exponential is not hard to get. We have that
	$$e^{t\,\cal N}=
	\left(\begin{array}{cc}
	e^{t A} & 0 \\
	0 & 1 
	\end{array}
	\right)=\left(
	\begin{array}{ccc}
	\cos a t & -\sin a t & 0 \\
	\sin a t & \cos a t & 0 \\
	0 & 0 & 1 
	\end{array}
	\right).$$\\[-10pt]
	\task[(Step 3)] \hspace{1cm}  We solve the new system $\cal N\, W=W'$, with $W(0)=\inv Q\,V(0)$.  Of course the solution of this equation is
	$$W=e^{t\,\cal N}(\inv Q\,V(0)).$$
	\task[(Step 4)] \hspace{1cm}  Finally, we obtain the desired solution of $\Mfr\,V=V'$, namely, $V=Q\,W$. In fact, from the previous step, we get $Q\,W(0)=V(0)$ and 
	$$V'=Q\,W'=Q\,\cal N\,W =Q\,(\inv Q\,\Mfr Q)(\inv Q\,V)=\Mfr\, V. $$
\end{tasks}

Now the execution of the steps. In the Step 1, according to the proof of the cited theorem, we must consider $M_s=\Mfr^2$ that a symmetric matrix and thus there exist a basis of eigenvectors that diagonalizes it. From this basis we construct the matrix $Q$ in the Step 1. A direct calculation shows that $M_s$ has 
an double eigenvalue $\lambda_1=-1$ and the other one  $\lambda_2=0$. The  vectors $v_1=(-\frac{1}{\sqrt{2}},0, \frac{1}{\sqrt{2}})$, $v_2=(0,1,0)$ and $v_3=(\frac{1}{\sqrt{2}},0, \frac{1}{\sqrt{2}} )$ are unit eigenvectors of $M_s$, $v_1$ and $v_2$ are associated to~$\lambda_1$. (Remember that if $P=(v_1\;v_2\;v_3)$ is the matrix whose columns are the column vectors $v_1$, $v_2$ and $v_3$ then
$\inv P\,M_s\,P={\rm diag}\,(-1,-1,0)$.) 
Looking carefully at the proof of the normal-form of a skew symmetric matrix theorem that we cite above, we construct the matrix $Q$. Its columns are $q_1=v_1$, $q_2=\frac{1}{\sqrt{-\lambda 1}}\Mfr\, v_1$ and $q_3=v_3$. More precisely,
$$Q=\left(
\begin{array}{ccc}
-\frac{1}{\sqrt{2}} & 0 & \frac{1}{\sqrt{2}} \\
0 & 1 & 0 \\
\frac{1}{\sqrt{2}} & 0 & \frac{1}{\sqrt{2}} 
\end{array}
\right).$$
Hence
$$\cal N=\inv Q\,\Mfr\,Q=\left(
\begin{array}{ccc}
0 & -1 & 0 \\
1 & 0 & 0 \\
0 & 0 & 0
\end{array}
\right).$$
Now, we execute the Step 2 and get
$$e^{t\,\cal N}=\left(
\begin{array}{ccc}
\cos t & -\sin  t & 0 \\
\sin  t & \cos  t & 0 \\
0 & 0 & 1 
\end{array}
\right),$$\\[-10pt]
Then the Step 3 gives us
$$\arraycolsep=0pt
\begin{array}{rcl}
W&\;=\;&{\arraycolsep=5pt e^{t\,\cal N}(\inv Q\,V(0))=\left(
	\begin{array}{ccc}
	\cos t & -\sin  t & 0 \\
	\sin  t & \cos  t & 0 \\
	0 & 0 & 1 
	\end{array}
	\right)
	\left(
	\begin{array}{ccc}
	-\frac{1}{\sqrt{2}} & 0 & \frac{1}{\sqrt{2}} \\
	0 & 1 & 0 \\
	\frac{1}{\sqrt{2}} & 0 & \frac{1}{\sqrt{2}}
	\end{array}
	\right)\left(
	\begin{array}{c}
	\text{V}_1(0) \\
	\text{V}_2(0) \\
	\text{V}_3(0) \\
	\end{array}
	\right)}\\[30pt]
&=&\**\left(
\begin{array}{c}
\frac{1}{2} \left(\sqrt{2}\, V_2(0) \sin t+\left(V_1(0)-V_3(0)\right) \cos t+V_1(0)+V_3(0)\right) \\[10pt]
\frac{\left(V_3(0)-V_1(0)\right) \sin t}{\sqrt{2}}+V_2(0) \cos t \\[10pt]
\frac{1}{2} \left(-\sqrt{2}\, V_2(0) \sin t+\left(V_3(0)-V_1(0)\right) \cos t+V_1(0)+V_3(0)\right)
\end{array}
\right).\\
&&
\end{array}$$
By substituting $\V1(0)$, $\V2(0)$ and $\V3(0)$ by their values (as row vectors), we obtain
$$\arraycolsep=5pt
W=
\left(
\begin{array}{ccc}
\sin  t & -\cos  t & 0 \\[5pt]
-\cos  t & -\sin t & 0 \\
0 & 0 & 1 \\
\end{array}
\right).$$
From this, it comes the solution $V$ of $\Mfr V=V'$:
$$V=Q W=\left(
\begin{array}{ccc}
-\frac{\sin  t}{\sqrt{2}} & \frac{\cos  t}{\sqrt{2}} & \frac{1}{\sqrt{2}} \\[5pt]
-\cos (t) & -\sin t & 0 \\
\frac{\sin  t}{\sqrt{2}} & -\frac{\cos (t)}{\sqrt{2}} & \frac{1}{\sqrt{2}} \\
\end{array}
\right).$$
The Steps are all done. By using the first row of $V$, we get the curve $f$, exactly as at the end of  the Solution~1, that is, $$f(t)=(1,0,0)+\int_0^t\sqrt 2\left (-\frac{\sin  u}{\sqrt{2}},\frac{\cos u }{\sqrt{2}} , \frac{1}{\sqrt{2}}\right)\,{\rm d}u=(\cos t,\sin t,t).$$
I hope this example, mainly its second solution, helps you in the proof of the theorem on the classification of the curves of constant curvatures that we make in the next section.
\eexe
\section{On Curves of Constant Curvatures}
In this section we will use the flowing notation, in order to simplify some calculations involving  inner product in $\R{2m}$.\def\hatx#1{#1_{\tiny\C{}}}

Given $X = (a_1, b_1, \ldots, a_{2m}, b_{2m}) \in\R{2m}$, the {\sl complex representation of} $X$ will be indicate by 
$$\hatx{X}= (a_1 +\i b_1, ...a_{2m}+\i b_{2m}) \in \C m,\quad \i=\sqrt{-1}.$$ 
With this notation, together with the real part of $z\in\C{}$, $\Re\, z$, we get the following useful properties:\\
\begin{tasks}[style=myenumerate,label-format={\normalfont},label-width={5mm}, item-indent={20pt}, item-format={\sl},label-align=right, after-item-skip={4pt}, column-sep={40pt}](1)
	\task $X\cdot Y=\Re{(\hatx X\cdot\hatx{Y}})$, 
	\task  $X\cdot X=\hatx X\cdot\hatx{X}$, that is, $\norma X=\norma{ \hatx X}$,
\end{tasks}\mbox{}\\
where the dot, $\cdot$ , indicates the usual inner products in $\R {2m}$ as well as that\footnote{$Z=\pupla z m\in\C m$ and $W=\pupla w m\in\C m$, $Z\cdot W=\sum_{j=1}^m z_j\barra{w_j}$.} in $\C{m}$.

\subsection{The Even Dimension Case}

We start with an example in $\R4$. Consider
\beq \label{fdoexe}f(t)=(r_1 \cos \left(a_1 t\right),r_1 \sin \left(a_1 t\right),r_2 \cos \left(a_2 t\right),r_2 \sin \left(a_2 t\right)),\quad t\in \R{},\eeq
where $a_1$, $a_2$, $r_1$ and $r_2$ positive numbers with $a_1\neq a_2$. The complex representation of $f$ in $\C2$ is
$$g(t)=\hatx {(f(t))}=(r_1e^{\i t a_1},r_2e^{\i t a_2}).$$
Note that $g'(t)=(\i r_1a_1e^{\i t a_1},\i r_2a_2e^{\i t a_2}),$   $g''(t)=(-r_1a_2^2e^{\i t a_1},- r_2a_1^2e^{\i t a_2}).$
More generally, $$g^{(j)}(t)=(\i^jr_1a_1^je^{\i t a_1},\i^j r_2a_2^je^{\i t a_2}).$$
So, we obtain  fast
$$\nu^2=\norma{f'}^2=\norma{g'}^2=\i r_1a_1e^{\i t a_1}(-\i r_1a_1e^{-\i t a_1})+\i r_2a_2e^{\i t a_2}(-\i r_2a_2e^{-\i t a_2})=r_1^2a_1^2+r_2^2a_2^2,$$
$$\norma{f''}^2=\norma{g''}^2=(a_1^2 r_1 \left(-e^{\i a_1 t}\right))(a_1^2 r_1 \left(-e^{-\i a_1 t}\right))+(a_2^2 r_2 \left(-e^{\i a_2 t}\right))(a_2^2 r_2 \left(-e^{-\i a_2 t}\right)=a_1^4 r_1^2+a_2^4 r_2^2$$
and 
$$f'\cdot f''=\Re{(g'\cdot g'')}=\Re{(-\i \left(a_1^3 r_1^2+a_2^3 r_2^2\right))}=0.$$
The notable here is that, given any $j$ and  $k$ in $\N$,  the inner product
$$\arraycolsep=0pt\def\arraystretch{1.2}
\begin{array}{rcl}
\** g^{(j)}(t)\cdot g^{(k)}(t)&\;=\;&\**(r_1 \left(\i a_1\right)^j e^{\i a_1 t})r_1 \left(-\i a_1\right)^k e^{-\i a_1 t}+r_2 \left(\i a_2\right)^j e^{\i a_2 t})r_2 \left(-\i a_2\right)^k e^{-\i a_2 t})\\[5pt]
\**&\;=\;&\**\i^j (-\i)^k \left(r_1^2 a_1^{j+k}+r_2^2 a_2^{j+k}\right)=(-1)^k\i^{j+k} \left(r_1^2 a_1^{j+k}+r_2^2 a_2^{j+k}\right)
\end{array}$$ does not depend on $t$ and thus  $f^{(j)}(t)\cdot f^{(k)}(t)=\Re{(g^{(j)}(t)\cdot g^{(k)}(t))}$ as well. Thus for any $j\in \N$
$\norma{\pvetorf f j}$ does not depend on $t$ and then all the curvatures of $f$ must be constant, according to the Corollary~\ref{coro4}. Of course it remains to verify that $f$ is at least $3$-regular, because  without this information, it makes no sense to calculate its curvatures. For this we return to $\R4$ and write
\beq\label{oneparam}\arraycolsep=2pt f(t)=\Sdet\left(
\begin{array}{c}
	r_1 \\
	0 \\
	r_2 \\
	0 \\
\end{array}
\right)\eeq
or $f(t)=M(t)f(0)$. Note that $M(t)$ is an one-parameter family of orthogonal matrices of determinant~1. Moreover, we have that  $f^{(j)}(t)=M(t)f^{(j)}(0)$, for any $j\in\N$. Collecting this information,  for $j\in\{1,2,3,4\}$, in an matrix form, we obtain
$$(f'(t)\;f''(t)\;f'''(t)\;f^{(4)}(t))=M(t)\,(f'(0)\;f''(0)\;f'''(0)\;f^{(4)}(t)).$$
In other words, the matrix of the derivatives  $(f'(t)\;f''(t)\;f'''(t)\;f^{(4)}(t))$ equals 
$$\arraycolsep=2pt \Sdet\donga4.$$
Thus the rank of $(f'(t)\;f''(t)\;f'''(t)\;f^{(4)}(t))$ is equal to the rank of $(f'(0)\;f''(0)\;f'''(0)\;f^{(4)}(0))$ that is equal to 4, which can be calculated directly or by using the next Lemma. Now we can calculate the curvatures, by using the Corollary~\ref{coro4}.
\vskip10pt
\begin{tasks}[style=myenumerateexe,label-format={\normalfont},label-width={2mm}, item-indent={20pt}, item-format={\sl},label-align=right, after-item-skip={4pt}, column-sep={40pt}](3)
	\task  $\kapa 1=\frac{\norma{f''}}{\nu^2}=\frac{\sqrt{a_1^4 r_1^2+a_2^4 r_2^2}}{a_1^2 r_1^2+a_2^2 r_2^2}$;
	\task $\kapa 2=\frac{a_2 \left(a_1^3-a_1 a_2^2\right) r_1 r_2}{\left(a_1^2 r_1^2+a_2^2 r_2^2\right) \sqrt{a_1^4 r_1^2+a_2^4 r_2^2}}$;
	\task $\kapa 3=\frac{a_1 a_2}{\sqrt{a_1^4 r_1^2+a_2^4 r_2^2}}$.
\end{tasks}
\vskip10pt
In the Theorem~\ref{teo52} below, we get, in particular,  the converse of this example: If a 3-regular curve in $\R 4$ has constant speed and curvatures, then it is as that in~(\ref{fdoexe}), up to an isometry of $\R4$. To finish, a remark on the angles $a_1$ and $a_2$. When they are equal, the trace of $f$ is contained in the intersection of the hyperplanes $\frac {x_3}{r_4}-\frac {x_1}{r_1}=0$ and  $\frac {x_4}{r_4}-\frac {x_2}{r_1}=0$, which has dimension 2. This would imply that $f',f'',f'''$ are linearly dependent, $k_2=0$ and $\kapa 3$ is not defined. Thus the condition $a_1\neq a_2$ guarantee that the curve $f$ is $4$-regular. 
\blema \label{dongalema}Given $m$ distinct real numbers $\puplasp am$ and any other real number $b$,  define
{\scriptsize
	$$\def\arraystretch{1.25}\arraycolsep=2pt D\pupla am= \dongam.$$
}and
{\scriptsize
	$$\def\arraystretch{1.25}\arraycolsep=2pt \til D(\puplasp am,b)= \dongamum.$$
}Then 
\benum
\item $\det D\pupla am=\left(\prod _{j=1}^m a_j^3\right) \prod _{i<j}^m \left(a_i^2-a_j^2\right)^2$.\\[5pt] \item $\det \til D(\puplasp am,b)=\left(\prod _{j=1}^m a_j^5\right) \prod _{i<j}^m \left(a_i^2-a_j^2\right)^2$.
\eenum
\elema
\proof We use induction on $m$. The idea of this proof is similar to that used in the calculation of the Vandermonde determinant. Differently,  here, the induction hypotheses is attained after two steps. In fact, we start writing $D= D\pupla am=\pupla {D}{2m}$, where $D_j$ denotes the $\jesima j$ row of $D$ and then by replacing  $a_1^2$ by $\lambda$ in the first row of $ D$ to construct the polynomial $$p(\lambda)=\det D((0,-\lambda,0,\lambda^2,0,\ldots,(-1)^m\lambda^m),D_2,\ldots,D_{2m}).$$  Hence, $p$ either is zero or has degree at most $m$. Since a matrix with either one zero row or two equal rows has zero determinant, we get that $p$  vanishes at $0$ and $a_{j}^2$, $2\leq j\le m$. We claim that the degree of $p$ is exactly $m$. In effect, the coefficient of $\lambda^m$ is  $c_m=(-1)^{1+2m}(-1)^m\,M_{1(2m)}=(-1)^{m+1}\,M_{1(2m)}$, where $M_{1(2m)}$ is the $(1,2m)$-minor of $D$. By factoring  out $a_1$, we get
$$M_{1(2m)}=a_1\det\dongaminor.$$
Denoting the $(2m-1)\times(2m-1)$  matrix above by $B=\pupla B {2m-1}$, we define $q(\lambda)$ to be $$\det((1,0,-\lambda,0,\ldots,(-1)^{m-1}\lambda^{m-1}),B_2,\ldots,B_{2m-1}).$$
Note that we have substituted $a_1^2$ by $\lambda$ in the first row of $B$. The  induction hypotheses guarantees that the polynomial $q$ has exactly degree $m-1$ because the  coefficient of $\lambda^{m-1}$ equals $$d_{m-1}=(-1)^{1+2m-1}(-1)^{m-1}\det D(a_2,a_3,\ldots,a_m)=(-1)^{m-1}\det D(a_2,a_3,\ldots,a_m)\neq 0.$$
Looking back in the polynomial  $p$, we see that $c_m=(-1)^{m+1}\,a_1\,q(a_1^2)\neq 0$, because $\partial q=m-1$ and $q$ vanishes at $a_2^2$, $a_3^2$, \ldots, $a_m^2$ that are $m-1$ distinct numbers. Now, we can factor $p$ as 
$$p(\lambda)=c_m\,\lambda(\lambda-a_2^2)(\lambda-a_3^2)\cdots(\lambda-a_m^2)=(-1)^{m+1}\,a_1\,q(a_1^2)\,\lambda(\lambda-a_2^2)(\lambda-a_3^2)\cdots(\lambda-a_m^2),$$
which implies
\beq \label{detdonga}\det D=p(a_1^2)=(-1)^{m+1}\,a_1\,q(a_1)\,a_1^2\,\prod_{j=2}^m(a_1^2-a_j^2)=(-1)^{m+1}\,a_1^3\,q(a_1^2)\,\prod_{j=2}^m(a_1^2-a_j^2).\eeq
Again from the induction hypothesis, we obtain that
$$\arraycolsep=0pt\def\arraystretch{1.5}
\begin{array}{rcl}
\**q(\lambda) &\;=\;&\**d_{m-1}(\lambda-a_2^2)(\lambda-a_3^2)\cdots(\lambda-a_m^2)\\
\**&\;=\;&\**(-1)^{m+1}\det D(a_2,a_3,\ldots,a_m)\,\prod_{j=2}^m(\lambda-a_j^2)\\
\** &\;=\;&\**(-1)^{m+1}\det D(a_2,a_3,\ldots,a_m)\,\prod_{j=2}^m(\lambda-a_j^2) \\[10pt]
\** &\;=\;&\**(-1)^{m+1}\left(\prod _{j=2}^m a_j^3\right)\left( \prod _{2\le i<j}^m \left(a_i^2-a_j^2\right)^2\right)\,\prod_{j=2}^m(\lambda-a_j^2),
\end{array}$$
whence
$$q(a_1^2)=(-1)^{m+1}\left(\prod _{j=2}^m a_j^3\right)\left( \prod _{2\le i<j}^m \left(a_i^2-a_j^2\right)^2\right)\,\prod_{j=2}^m(a_1^2-a_j^2)$$
which substituted in~(\ref{detdonga}) yields finally
$$\det D=\left(\prod _{j=1}^m a_j^3\right) \prod _{i<j}^m \left(a_i^2-a_j^2\right)^2.$$
The result in (ii) follows easily from~(i), 
completing the proof.\qed
\brem The element $d_{ij}$ of the matrix $D\pupla am$ equals $\derivada fj(0)\cdot e_i$, where $e_i$ is the $\jesima i$ vector of the canonical basis of the $\R {2m}$ and $f$ is as in the following classification  theorem.
\erem 
\bteo \label{teo52}Let $\nu>0$ and $\kappa_j$, $1\le j\le n-1$, be constants such that $\kappa_j>0$, $1\le j\le n-2$, and $\kappa_{n-1}\neq 0$.  Let $f$ be a curve with speed $\nu$ and  curvatures  $\kappa_j$, $1\le j\le n-1$. Suppose that $n$ is even. Then there exist positives real numbers $a_j$ and $r_j$, $1\le j\le m$, $m=\frac n2$, such that, up to an isometry,
$$f(t)=(r_1e^{\i t\,a_1},r_2e^{\i t\,a_2},\ldots,r_{m}e^{\i t\,a_m}),$$
where $e^{\i t\,a_j}=\cos ta_j+\i \sin t a_j$.
\eteo 
\proof Some of the ideas in this proof are suggested in~\cite{Kuhnel}  ({\bf 2.16-Remark}). We give a full proof for $n=4$ and believe that its extension to the general case is very easy. 

Let $\til f$ be the $3$-regular curve obtained from the unique solution $V=(\V1,\V2,\V3,\V4)$ of the linear system $\Mfr\, V=V'$ or, more explicitly,
\def\quadra #1{(#1_1,#1_2,#1_3,#1_4)}\def\quadrasp #1{#1_1,#1_2,#1_3,#1_4}
$$\sys4,$$
with the initial condition $V(0)=\quadra Q$, where $\{\quadrasp Q\}$  will be an orthonormal basis defined as follows. The existence of the curve $\til f$ is guaranteed by the Theorem~\ref{fundteo}. Remember, $\til f(t)=\int_0^t \nu\V1(u) {\rm d}u$.
Since $\Mfr$ is skew symmetric, we obtain from the normal-form of a skew symmetric matrix theorem an orthonormal basis $\{\quadrasp Q\}$  that reduces  $\Mfr$ to\def\Nfr{\cal N}
$$\cal N=\left(
\begin{array}{cccc}
0 & a_1 & 0 & 0 \\
-a_1 & 0 & 0 & 0 \\
0 & 0 & 0 & a_2 \\
0 & 0 & -a_2 & 0 \\
\end{array}
\right),
$$
where $0<a_1$ and  $0<a_2$.
Actually,
$\cal N=\trans Q\,\Mfr\,Q$, where $Q$ is the $4\times 4$ matrix whose columns are the vectors $\quadrasp Q$ viewed as column vectors. By rewriting the equation $\Mfr\, V=V'$, with $V(0)=\quadra Q$, in terms of $\cal N$, we get
$(Q\,\Nfr\,\trans Q)\,V=V'$, which is the same as
$$\Nfr\,\trans (Q\, V)=\trans Q\,V'=(\trans Q\, V)',\quad (\trans Q\,V)(0)= E=\quadra e,$$
where $\{\quadrasp e\}$ is the canonical basis of $\R4$. By the uniqueness of solution of the equation above, we must have $Q\,V=e^{t \Nfr}\, E$ or $V=\trans Q\,e^{t \Nfr}\,E$. In other words,
$$\arraycolsep=2pt V(t)=\transmatrizQ\,\sencos\,\colcan.$$
The fours first rows of the product above is
$$\arraycolsep=0pt\def\arraystretch{1.5}\V1(t)=\left(
\begin{array}{c}
q_{11} \cos \left(a_1 t\right)-q_{21} \sin \left(a_1 t\right) \\
q_{11} \sin \left(a_1 t\right)+q_{21} \cos \left(a_1 t\right) \\
q_{31} \cos \left(a_2 t\right)-q_{41} \sin \left(a_2 t\right) \\
q_{31} \sin \left(a_2 t\right)+q_{41} \cos \left(a_2 t\right) \\
\end{array}
\right)=\left(
\begin{array}{c}
\sqrt{q_{11}^2+q_{21}^2} \left(-\sin \left(a_1 t-c_1\right)\right) \\
\sqrt{q_{11}^2+q_{21}^2} \cos \left(a_1 t-c_1\right) \\
\sqrt{q_{31}^2+q_{41}^2} \left(-\sin \left(a_2 t-c_2\right)\right) \\
\sqrt{q_{31}^2+q_{41}^2} \cos \left(a_2 t-c_2\right) \\
\end{array}
\right),$$
where $c_1$ and $c_2$ are such that
$$\arraycolsep=0pt\sin c_1=\frac{q_{11}}{\sqrt{q_{11}^2+q_{21}^2}},\
\cos c_1  =\frac{q_{21}}{\sqrt{q_{11}^2+q_{21}^2}},\ \sin c_2  =\frac{q_{31}}{\sqrt{q_{31}^2+q_{41}^2}},\ \cos c_2  =\frac{q_{41}}{\sqrt{q_{31}^2+q_{41}^2}}.$$
Note that $\V1$ is in fact an unit vector field.
Integration of $\nu \V1$ yields
$$\arraycolsep=0pt \til f(t)=\int\nu\V1(t){\rm d}t=\nu\left(
\begin{array}{c}
\frac{\sqrt{q_{11}^2+q_{21}^2} }{a_1}\, \cos \left(a_1 t-c_1\right)\\
\frac{\sqrt{q_{11}^2+q_{21}^2} }{a_1}\, \sin \left(a_1 t-c_1\right)\\
\frac{\sqrt{q_{31}^2+q_{41}^2}}{a_2} \, \cos \left(a_2 t-c_2\right)\\
\frac{\sqrt{q_{31}^2+q_{41}^2} }{a_2}\, \sin \left(a_2 t-c_2\right))
\end{array}
\right)
=\left(
\begin{array}{c}
r_1\, \cos \left(a_1 t-c_1\right)\\
r_1\, \sin \left(a_1 t-c_1\right)\\
r_2\, \cos \left(a_2 t-c_2\right)\\
r_2\, \sin \left(a_2 t-c_2\right))
\end{array}
\right).$$
Since we know that $\til f$ is 3-regular, it follows that $a_1\neq a_2$, according to the remark that we did at the end of the example above. Furthermore, $\nu_{\til f}=\nu$ and $\kapag j=\kapa j$, $1\le j\le 3$. 
A simple computation shows that the curve $f$ obtained from $\til f$ by taking $c_1=c_2=0$ has also the speed $\nu$ and  curvatures $\kapa j$, $1\le j\le 3$, for the inner products involving its derivatives are exactly the same as those of $\til f$. We are done.\qed
\brem In~\cite{Lucas} and ~\cite{Rolfpp},  we find an approach involving one-parameter subgroup of isometries to perform the calculations of the curvatures as well as to prove the constant curvature classification theorem. Here, in~(\ref{oneparam}), we have one example of such a subgroup, namely, $M(t)$, since ${M(t_1+t_2)=M(t_1)M(t_2)}$.  
\erem
\subsection{The Odd Dimension Case}
Consider the following generalization of the circular helix
\beq \label{fdoexe5}f(t)=(r_1 \cos \left(a_1 t\right),r_1 \sin \left(a_1 t\right),r_2 \cos \left(a_2 t\right),r_2 \sin \left(a_2 t\right),b t),\quad t\in \R{},\eeq
where $a_1$, $a_2$, $r_1$ and $r_2$  positive numbers with $a_1\neq a_2$ and $b\neq 0$. The complex representation of $f$ in $\C3$ is
$$g(t)=\hatx {(f(t))}=(r_1e^{\i t a_1},r_2e^{\i t a_2},bt),$$
Note that $g'(t)=(\i r_1a_1e^{\i t a_1},\i r_2a_2e^{\i t a_2},b),$   $g''(t)=(-r_1a_1^2e^{\i t a_1},- r_2a_1^2e^{\i t a_2},0).$
More generally, $$g^{(j)}(t)=(\i^jr_1a_1^je^{\i t a_1},\i^j r_2a_1^je^{\i t a_2},0),$$
whenever $j>1$.
So
$$\nu^2=\norma{f'}^2=\norma{g'}^2=\i r_1a_1e^{\i t a_1}(-\i r_1a_1e^{-\i t a_1})+\i r_2a_2e^{\i t a_2}(-\i r_2a_2e^{-\i t a_2})+b^2=r_1^2a_1^2+r_2^2a_2^2+b^2,$$
$$\norma{f''}^2=\norma{g''}^2=(a_1^2 r_1 \left(-e^{\i a_1 t}\right))(a_1^2 r_1 \left(-e^{-\i a_1 t}\right))+(a_2^2 r_2 \left(-e^{\i a_2 t}\right))(a_2^2 r_2 \left(-e^{-\i a_2 t}\right)=a_1^4 r_1^2+a_2^4 r_2^2$$
and 
$$f'\cdot f''=\Re{(g'\cdot g'')}=\Re{(-\i \left(a_1^3 r_1^2+a_2^3 r_2^2\right))}=0.$$
As before, we get that 
$f^{(j)}(t)\cdot f^{(k)}(t)=\Re{(g^{(j)}(t)\cdot g^{(k)}(t))}$  does not depend on $t$ and  thus all of the curvatures of $f$ must be constant. Of course it remains to verify that $f$ is at least $3$-regular, because  without this information, it makes no sense to calculate its curvatures. For this we return to $\R4$ and write
\beq\label{oneparamb}\arraycolsep=2pt f(t)=\Sdetcinco\left(
\begin{array}{c}
	r_1 \\
	0 \\
	r_2 \\
	0 \\
	0
\end{array}
\right)+\left(\begin{array}{c}
	0 \\
	0 \\
	0 \\
	0 \\
	b\,t
\end{array}
\right)\eeq
or $f(t)=N(t)f(0)+(0,0,0,0,b\,t)$. Note that $N(t)$ is also a family of one-parameter of orthogonal matrix of determinant~1. By collecting the derivatives $\derivada fj(t)$,  for $j\in\{1,2,3,4,5\}$, in an matrix, we obtain that
$(f'(t)\;f''(t)\;f'''(t)\;f^{(4)}(t)\;f^{(5)}(t))$ equals 
$$\arraycolsep=2pt \Sdetcinco\dongacinco.$$
The Lemma~\ref{dongalema} implies that the rank of $(f'(t)\;f''(t)\;f'''(t)\;f^{(4)}(t)\;\derivada f5(t))$ is equal to 5, which is its rank at $t=0$. Hence $f$ is 5-regular and then  we can calculate its curvatures,  by using the Corollary~\ref{coro4}.
\vskip10pt
\begin{tasks}[style=myenumerateexe,label-format={\normalfont},label-width={5mm}, item-indent={20pt}, item-format={\sl},label-align=right, after-item-skip={4pt}, column-sep={40pt}](2)
	\task  $\kapa 1=\frac{\sqrt{a_1^4 r_1^2+a_2^4 r_2^2}}{a_1^2 r_1^2+a_2^2 r_2^2+b^2}$;
	\task $\kapa 2=\frac{\sqrt{a_1^6 b^2 r_1^2+a_2^2 r_2^2 \left(a_2^4 b^2+a_1^2 \left(a_1^2-a_2^2\right)^2 r_1^2\right)}}{\sqrt{a_1^4 r_1^2+a_2^4 r_2^2} \left(a_1^2 r_1^2+a_2^2 r_2^2+b^2\right)}$;
	\task*(2) $\kapa 3=\frac{a_1^2 a_2^2 \left(a_1^2-a_2^2\right) r_1 r_2}{\sqrt{a_1^4 r_1^2+a_2^4 r_2^2} \sqrt{a_1^6 b^2 r_1^2+a_2^2 r_2^2 \left(a_2^4 b^2+a_1^2 \left(a_1^2-a_2^2\right)^2 r_1^2\right)}}$;
	\task $\kapa 4=\frac{a_1 a_2 b \sqrt{a_1^4 r_1^2+a_2^4 r_2^2}}{\sqrt{a_1^2 r_1^2+a_2^2 r_2^2+b^2} \sqrt{a_1^6 b^2 r_1^2+a_2^2 r_2^2 \left(a_2^4 b^2+a_1^2 \left(a_1^2-a_2^2\right)^2 r_1^2\right)}}$.
\end{tasks}
\vskip10pt

Now the theorem for the odd dimension case.
\bteo \label{teo53}Let $\nu>0$ and $\kappa_j$, $1\le j\le n-1$, be constants such that $\kappa_j>0$, $1\le j\le n-2$, and $\kappa_{n-1}\neq 0$.  Let $f$ be a curve with speed $\nu$ and  curvatures  $\kappa_j$, $1\le j\le n-1$. Suppose that $n$ is odd. Then there exist positives real numbers $a_j$ and $r_j$, $1\le j\le m$, $m=\frac {n-1}2$, and $b\neq 0$, such that, up to an isometry,
$$f(t)=(r_1e^{\i t\,a_1},r_2e^{\i t\,a_2},\ldots,r_{m}e^{\i t\,a_m},b\,t).$$
\eteo 
\proof It is similar to that of the Theorem~\ref{teo52}. The only difference is that the matrix $\Nfr$ becomes
\beq\label{Ngood}\arraycolsep=2pt\Nfr=\left(
\begin{array}{ccccc}
	0 & a_1 & 0 & 0&0 \\
	-a_1 & 0 & 0 & 0&0 \\
	0 & 0 & 0 & a_2&0 \\
	0 & 0 & -a_2 & 0 &0\\
	0 & 0 &0 & 0 &0
\end{array}
\right),
\eeq
whence, 
$$\arraycolsep=2pte^{t \Nfr}=\left(
\begin{array}{ccccc}
\cos \left(a_1 t\right) & \sin \left(a_1 t\right) & 0 & 0 & 0 \\
-\sin \left(a_1 t\right) & \cos \left(a_1 t\right) & 0 & 0 & 0 \\
0 & 0 & \cos \left(a_2 t\right) & \sin \left(a_2 t\right) & 0 \\
0 & 0 & -\sin \left(a_2 t\right) & \cos \left(a_2 t\right) & 0 \\
0 & 0 & 0 & 0 & 1 
\end{array}
\right).$$
\vskip10pt
\noindent
This leads to 
$$f(t)=(r_1 \cos \left(a_1 t\right),r_1 \sin \left(a_1 t\right),r_2 \cos \left(a_2 t\right),r_2 \sin \left(a_2 t\right),b t).$$
\vskip10pt
\noindent Of course, there exist other possibilities for $\Nfr$, for instance, $\Nfr$ could be
$$\arraycolsep=2pt\Nfr=\left(
\begin{array}{ccccc}
0 & a_1 & 0 & 0&0 \\
-a_1 & 0 & 0 & 0&0 \\
0 & 0 & 0 &0&0 \\
0 & 0 &0 & 0 &0\\
0 & 0 &0 & 0 &0
\end{array}
\right),
$$
which would lead to a curve of the kind
$$f(t)=(r_1 \cos \left(a_1 t\right),r_1 \sin \left(a_1 t\right),b_1\,t,b_2\,t,b_3\,t),$$
which has $\kapa 3=0$ and undefined $\kapa 4$. Similarly, any  possibility other than $\Nfr$ in (\ref{Ngood}) is discarded. 
Of course, we have considered  $n=5$. The general case is almost as this.\qed


\section{Tubes in \R{n+1}}
Probably, the torus of revolution $T(a,b)\subset\R3$, generate by revolving the circle $${(x-b)^2+z^2=a^2},\quad y=0,$$
around the $z$-axis,  is the best known example of tube around a curve $f$ in \R3. This curve, which is called  the axis of $T(a,b)$, is  given by $f(t)=(b\cos t,b\sin t,0)$, $0\le t\le 2\pi$, whose trace is the circle $x^2+y^2=b^2$ in the $xy$-plane.  $T(a,b)$ is a regular surface, whenever $0<a<b$ and 
$$X(v,t)=( (b+a \cos v)\cos t, (a\cos v+b)\sin t,a \sin v),\quad (v,t)\in[0,2\pi]\times[0,2\pi],$$
is onto $T(a,b)$. Of course the restriction of $X$ to the open rectangle  $(0,2\pi)\times(0,2\pi)$ is a parametrization of $T(a,b)$. A good way of rewriting this map consists in replacing $(a\cos v,a\sin v)$ by $(x_1,x_2)$  to obtain a new map
$$G(x_1,x_2,t)=( (b+ x_1)\cos t, (x_1+b)\sin t,x_2),\quad (x_1,x_2,t)\in S^1(a)\times[0,2\pi] \subset\R3,$$
where $S^1(a)$ is the circle $x_1^2+x_2^2=a^2$ of the $x_1x_2$-plane of the tridimensional space \R3 with euclidean coordinates $(x_1,x_2,t)$. This map  sends the right circular cylinder over $S^1$ having height~$2\pi$ onto~$T(a,b)$, transforming  each copy of $S^1$ in the cylinder onto a copy of the generator circle of~$T(a,b)$. By composing $G$ with $h(v,t)=(a\,\cos v,a\,\sin v,t)$ we can recover $X$.   At this moment, it is convenient to observe that $\til G$, the extension of $G$ to the solid cylinder $ B^{[2]}[a]\times{[0,2\pi]}$, fills the solid enclosed by $T(a,b)$, where $B^{[2]}[a]$ denotes the compact disk of radius $a$ of the $x_1x_2$-plane. $\til G$ becomes
\beq \label{gtil}\til G(x_1,x_2,t)=( (b+ x_1)\cos t, (x_1+b)\sin t,x_2),\quad (x_1,x_2,t)\in B^{[2]}[a]\times[0,2\pi]\subset\R3\eeq
which is onto the solid torus  $\barra {T(a,b)}$. The next picture shows  $X$,  $G$ and the torus $T(a,b)$ together with the Frenet frame $\{{\mathbf T},{\mathbf N},{\mathbf B}\}$ of its axis $f(t)=(b\cos t,b\sen t,0)$. It is not hard to see that 
$${\mathbf T}(t)=(-\sin t,\cos t,0),\quad {\mathbf N}(t)= (-\cos t,-\sin t,0),\quad {\mathbf B}(t)=(0,0,1).$$

\includegraphics[scale=0.4]{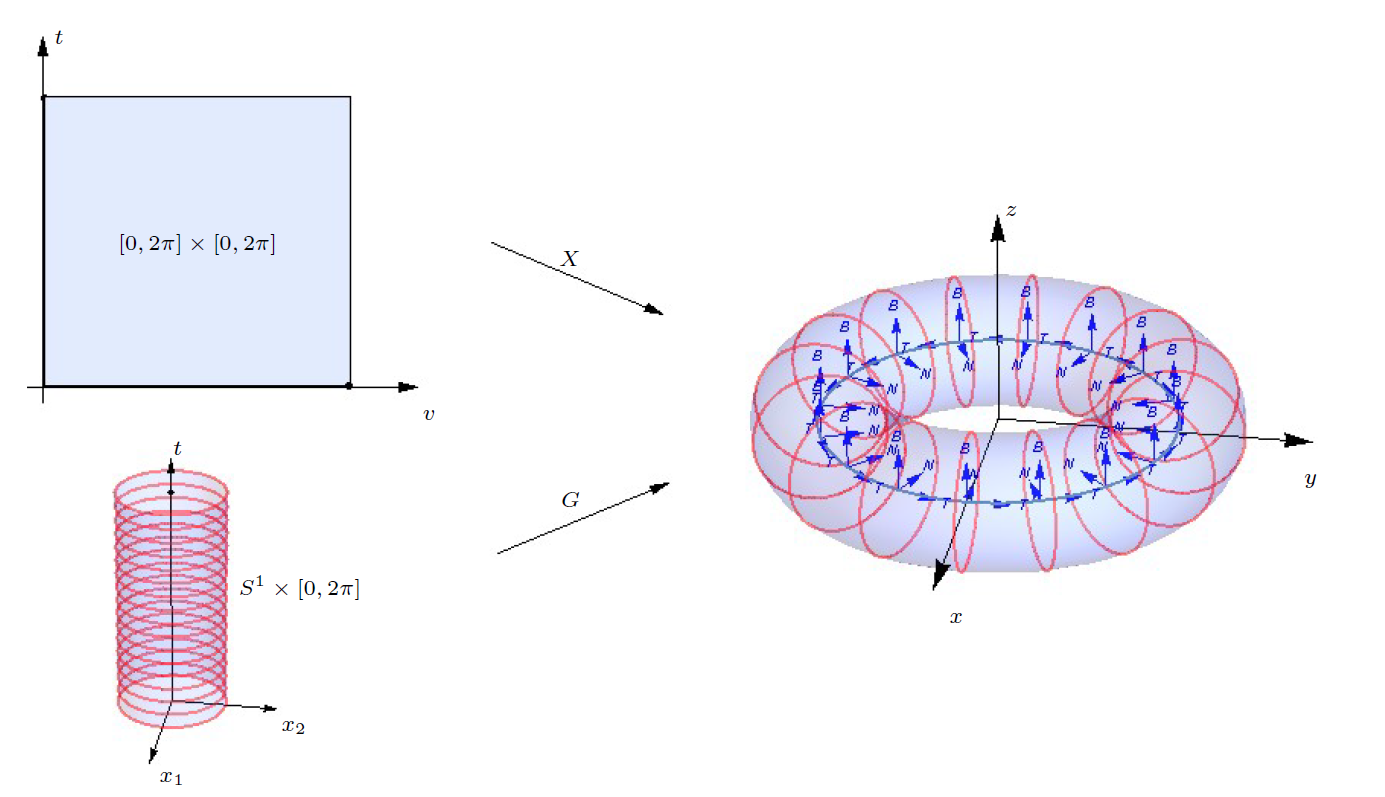}
Furthermore, $\kappa=1/b$ and $\tau=0$ are the curvature and torsion of $f$. In order to generalize the ideas in the exposition above, we verify easily that $G$ ($\til G$) coincides with 
$$f(t)-x_1\,{\mathbf N}(t)+ x_2\,\Bb(t).$$
The minus signal was included because $\{-\Nb,\Bb\}$ has the same orientation as $\{e_2,e_3\}$  in the revolving $x_1x_2$-plane. However, we can choose with no additional difficulties
$$f(t)+x_1\,\Nb(t)+ x_2\,\Bb(t)=(\left(b- x_1\right)\cos t,  \left(b-x_1\right)\sin t, x_2)$$
to generate the same torus. From this model, we will construct tubes of either spheres or disks or any other nice region around a $n$-regular curve $f$ in \R{n+1}. Going back to~(\ref{gtil}) and taking the rule of $\til G$, we define the following map
$$H(x_1,x_2,t)=( (b+ x_1)\cos t, (x_1+b)\sin t,x_2),\quad (x_1,x_2,t)\in\Omega,$$
where $\Omega$ is the open set $ B^{[2]}(\til a)\times(0,2\pi)$,  $B^{[2]}(\til a)$ is the open disk of radius $\til a$, for some $\til a$ such that  $b>\til a>a$. It is not hard to see that $H$ is injective and the absolute values of the its jacobian determinant is 
$$|\det JH|=\left|\det\left ( \dparcial H {x_1},\dparcial H {x_2},\dparcial H t\right)\right|=\left|\frac{\partial(x,y,z)}{\partial(x_1,x_2,t)}\right|=\norma{\dparcial H {x_1}\wedge\dparcial H {x_2}\wedge\dparcial H t}=|b+x_1|=b+x_1> 0,$$
since $-b<-\til a<x_1<\til a<b$. Thus, for a sufficiently small $\eps>0$ we can calculate the volume of $\barra T_\eps\subset\barra{T(a,b)}$, the image under $H$ of the compact set $ B^{[2]}[a]\times[\eps,2\pi-\eps]$, 
which, by using the change of variable in multiple integrals, is 
$$\arraycolsep=0pt\def\arraystretch{2}
\begin{array}{rcl}
	\**\vol\;\barra T_\eps &\;=\;&\**\int\!\!\!\!\int\!\!\!\!\int_{\barra T_\eps}\dxx\dy\dz=\int_{\eps}^{2\pi-\eps}\left(\int\!\!\!\!\int_{B^{[2]}[a]}(b+x_1)\dx 1\dx 2\right)\dt\\
	\**&\;=\;&\**b\int_{\eps}^{2\pi-\eps}\left(\int\!\!\!\!\int_{B^{[2]}[a]}\dx 1\dx 2\right)\dt=2 \pi  a^2 b (\pi -\eps),
\end{array}$$

\includegraphics[scale=0.55]{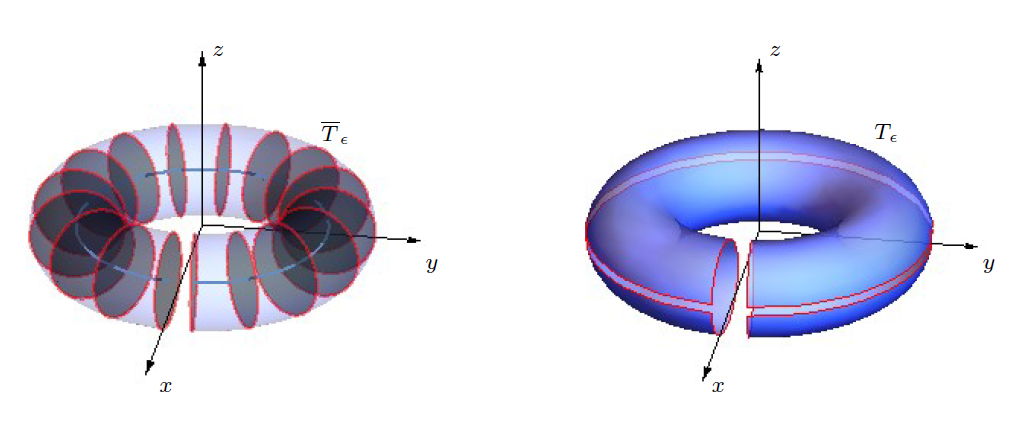}
where we have used $\int\!\!\int_{B^{[2]}[a]}x_1\dx 1\dx 2=0$, fact that is intimately connected to the barycenter\footnote{The barycenter (see also~\refeq{bary}, page~\pageref{bary}) of a compact region $K\subset\R 2$ with $\area>0$ is the point
	$$B_K=\frac 1{\area\, K}(\int\!\!\!\!\int_{K}x\dxx\dy,\int\!\!\!\!\int_{K}y\dxx\dy).$$}
of the disk $D[a]$, namely its center $(0,0)$. Hence $\vol\,\barra{T(a.b)}=\lim_{\eps\rightarrow 0}\vol\,\barra T_\eps=2\pi^2a^2b,$
which is a well known result and that can be obtained from the second theorem of Pappus applied to the torus: the volume of $\barra{T(a,b)}$ equals the area of $B^{[2]}[a]$ times the length of the curve described by the center of the revolving  disk.  By the same reasoning, we can find the  area of the $T(a,b)$, by working with the restriction of $X$ to the rectangle $R_\eps=[\eps, 2\pi-\eps]\times[\eps, 2\pi-\eps]$ whose image is $T\eps$, as in the picture above (see~\cite{Manfra}{\bf -2.5-Example 5}). The calculations are as follows.
$$\arraycolsep=0pt\def\arraystretch{2.5}
\begin{array}{rcl}
	\**\area\; T_\eps &\;=\;&\**\int\!\!\!\!\int_{R_\eps}\norma{\dparcial Xt\wedge\dparcial Xv}\dt\dv=
	\**\int\!\!\!\!\int_{R_\eps}\sqrt{\norma{\dparcial Xt}^2\norma{\dparcial Xv}^2-\left(\dparcial Xv\cdot\dparcial Xv\right)^2}\dt\dv\\
	\**&\;=\;&\**\int_{\eps}^{2\pi-\eps}\left(\int_{\eps}^{2\pi-\eps} (b+a \cos v)\dv\right)\dt=4 a^2 \eps \sin  \eps-4 \pi  a^2 \sin  \eps+4 a b \eps^2-8 \pi  a b \eps+4 \pi ^2 a b,
\end{array}$$
whose limit, as $\eps$ tends to zero,  is $4\pi^2ab$, which is the area of $T(a,b)$. Now, we have the length of $S^1(a)$ times the length of curve described by the center of the revolving circle, which comes from the first theorem of Pappus. The cited theorems of Pappus concerns either surfaces or solids of revolution generated by either a simple plane curve or the region enclosed by it. 
Next, we will deal with these theorems in a more general sense in any dimension. A generalization of them, in the tridimensional case, appears in~\cite{Good}. 

Given an interval $J$, let $\Omega_J$ be the set
$\Omega_J(S_t)=\cup_{t\in J}S_t\subset\R{n+1}$, where 
$$S_t=\{(\puplasp xn,t);\quad \pupla xn\in \til S_t\subset\R n\}$$ and $\til S_t$, the projection of $S_t$ in $\R n$, has positive volume. So, $\Omega_J$ is the solid whose intersection with the hyperplane $x_{n+1}=t$ is $S_t$,  for each $t\in J$. When all projections coincide with a certain region $S$, $\Omega_J(S)$ is a {\sl solid cylinder} with cross section $S$, case in which we write $S_t=(S,t)$ and $\vol\, \til S_t=\vol\, S>0$, of course measured in $\R n$.

\includegraphics[scale=0.52]{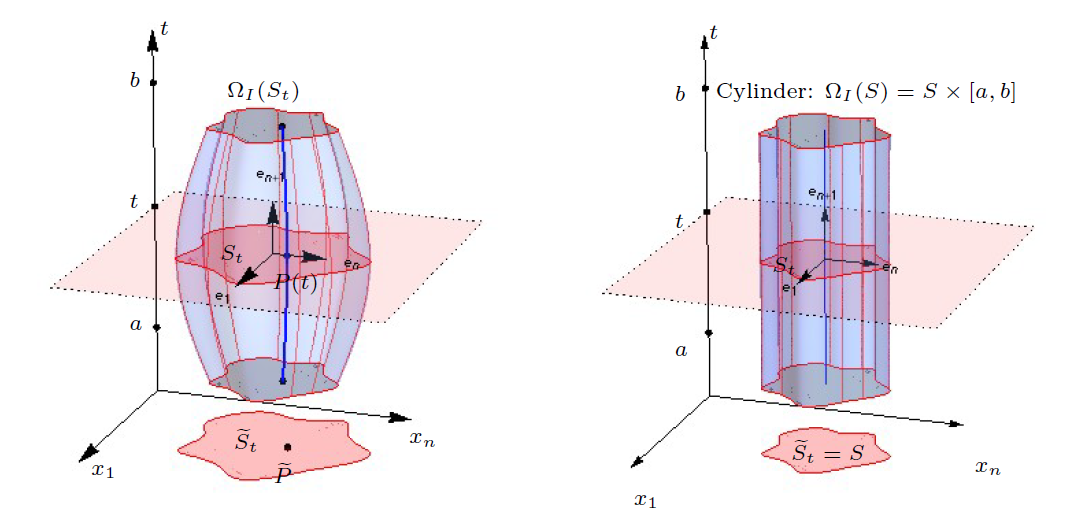} 
Now, given $I=[a,b]\subset J$ and $\funcao{F}{\Omega_J(S_t)}{\R{}}$ continuous, we suppose that the integral 
$$\int_{\Omega_I(S_t)}F\pupla x{n+1}\dx1\dx2\ldots\dx{n+1}=\int_{\Omega_I(S_t)}F(\puplasp x{n},t)\dx1\dx2\ldots\dt$$
always exists. Thus, by using the Fubini theorem, we get
$$\int_{\Omega_I(S_t)}F(\puplasp x{n},t)\dx1\dx2\ldots\dx n\dt=\int_a^b\left(\int_{\til S_t}F(\puplasp x{n},t)\dx1\dx2\ldots\dx n\right)\dt.$$
In particular,
$\vol\;\Omega_J(S_t)=\int_a^b\vol\; \til S_t\dt.$


Let \funcao{f}{J}{\R{n+1}} be a regular curve with Frenet referential $\reffrenetnew{n+1}$ and $\Omega_J(S_t)$ as above. We construct a tube around $f$, which we call  the {\sl tube generated by  $\Omega_J(S_t)$} as follows. Choose in each slice $\R n\times \{t\}$ one point, say $P=P(t)=(\puplaspdet pnt,t)$  (in general, we choose $P(t)$ in $S_t$) and put $S_t$ in the hyperplane orthogonal to the $f$ at $f(t)$, making the point $P(t)$ stay on $f(t)$ and each $e_j$ of the canonical basis of $\R n$ on $\V{j+1}$. This process ends with what we call {\sl the tube around $f$ generated by $\Omega_J(S_t)$ and centered at $P(t)$}, which will be indicated by $\Gamma_P(f, \Omega_J(S_t))$ or simply $\Gamma$ with little danger of confusion.  Intuitively, everything happens as if a region $S$ is moving along a curve by expanding or shrinking and always stuck to the curve by a point $P$. We also say  that the curve $f$ is the {\sl axis of the tube}. 

$\Gamma_P(f, \Omega_J(S_t))$  can be obtained as image of  $\Omega_J(S_t)$ under the map
\funcao{G}{\R{n}\times J}{\R{n+1}}  given  by 
\beq\label{mapG}
G(X,t)=G(\puplasp xn,t) =f(t)+\sum_{j=1}^n (x_j-p_j(t))\V{j+1}(t).
\eeq
Of course that this map changes with a change in $\til P(t)=\puplalinhacomt pn$.
Geometrically, $G$ sends each cross-sections $S_t$ onto a copy of it, also called {\sl cross section of the tube},  into the hyperplane passing through $f(t)$ and parallel to the subspace generated by $\V 2$, $\V3$, \ldots, $\V{n+1}$,  gluing $P(t)$ and $f(t)$ and laying $e_j$, the $j$-th element of the canonical basis,  on $\V{j+1}$, for $j$ running from 1 to $n$, as we see in the picture below.

\includegraphics[scale=0.5]{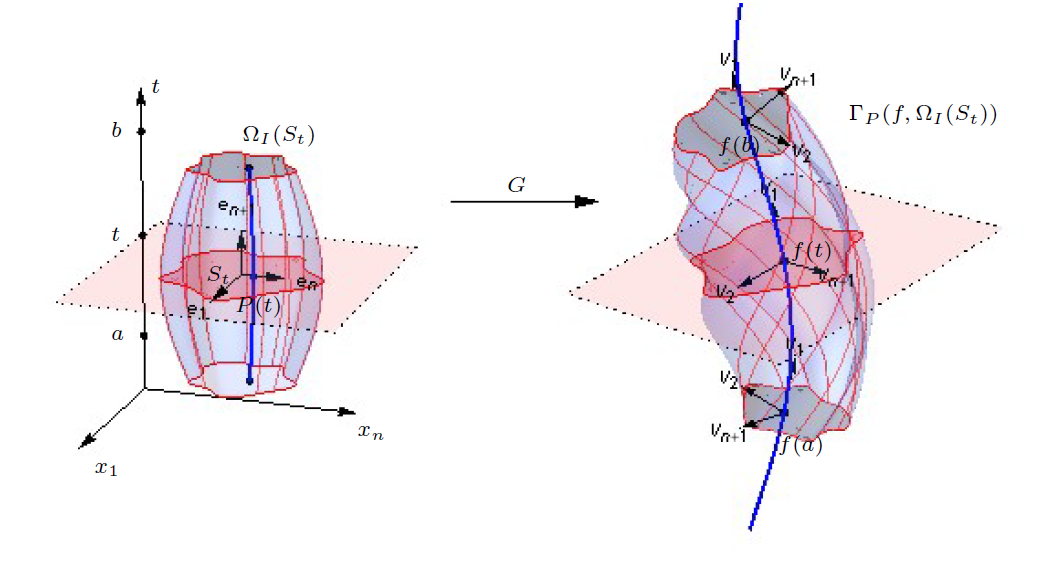}

\brem \label{remtubes} 
Actually, in (\ref{mapG}), we could use $\til F_t(\V{j+1}(t))$ instead $\V{j+1}(t)$, where $\til F_t$ is a differentiable family of linear isometries of $W_t={\rm span}\{\V j(t),\; 2\le j\le n+1\}$  induced by a differentiable family of linear isometry of $\R n$, $F_t$, producing so a {\sl twisted tube}. If $\funcao{L_t}{\R n}{W_t}$ is the linear isometry such that $L_t(e_j)=\V{j+1}(t)$, then $\til F_t$ is constructed as in the commutative  diagram below. Note that when
each cross section $S_t$ of $\Omega_J(S_t)$ is a disk centered at the $(0,0, ...,t)$, then the disk tube constructed from a family of isometries $F_t$ must coincide with that obtained from G in (\ref{mapG}), which is constructed by using the identity map. In fact, disks centered at origin are invariant under linear  isometries.
$$\begin{tikzcd}
\R n \arrow[r, "L_t"]
& W_t \arrow[d, "\til F_t"] \\
\R n  \arrow[u, "F_t"]
&\arrow[l, "\inv L_t"] W_t
\end{tikzcd}$$
\erem
Of course a tube  $\Gamma=\Gamma_P(f,\Omega(S_t))$ may have self-intersections, which depend on the curvatures of $f$ and the size of the $S_t$. If $\Gamma_P$ has no such a self-intersection, that is, if $G$ is injective, $\Gamma$ is called a {\sl genuine tube}. In the picture above, we have a {\sl genuine  tube}. In the cases where the cross section of $\Omega_J(S_t)$ are disks, we say that $\Gamma$ is a {\sl disk tube} around $f$. Below, we see two cases of self-intersections in {\sl disk tubes} around the same curve. In both cases, a reduction in the sizes of the cross sections would produce genuine disk tubes. The next theorem guarantees the existence of genuine disk tubes.

\includegraphics[scale=0.6]{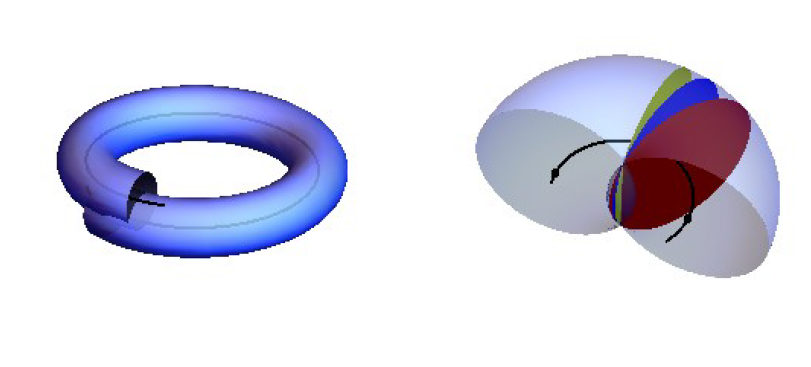}

\bteo{} \label{tuboteo}Let \funcao{f}{J}{\R{n+1}} be a $n$-regular curve with speed $\nu$, Frenet apparatus $$\cal A=\{\funcoesst \kappa,\framest V\},
$$ 
and consider a compact interval $I=[a,b]\subset J$. If $f$ is injective, then there exists $r>0$ such that the map $G$ in~(\ref{mapG}), with $p_j(t)=0$, $1\le j\le n$, carries $B^{[n]}[r]\times I$ onto the compact genuine disk tube $\Gamma=G(B^{[n]}[r]\times I)$, where $B^{[n]}[r]$ is the compact disk of radius $r$ in $\R n$. 
\eteo
\proof We start by calculating the absolute value of the jacobian determinant of $G$, that is given by
$$|\det JG|=\left|\det\left(\dparcial G {x_1}\;\dparcial G {x_2}\cdots\dparcial G {x_n}\; \dparcial G t\right)\right|=\norma{\dparcial G {x_1}\wedge\dparcial G {x_2}\wedge\cdots\wedge\dparcial G {x_n}\wedge\dparcial G t}.$$
From $G(\puplasp xn,t) =f(t)+\sum_{j=1}^n x_j\V{j+1}$, we get that 
$\dparcial G {x_j}=\V{j+1}$, for  $1\le j\le n$, and that 
$$\arraycolsep=0pt\def\arraystretch{1.2}
\begin{array}{rcl}
\** \dparcial G {t}&\;=\;&\**f'+\sum_{j=1}^n x_j\V{j+1}'=\nu\left(\V1+\sum_{j=1}^{n-1} x_j(-\kapa j\V j+\kapa {j+1}\V {j+2})- x_{n}\kapa n\V n\right)\\
\**&\;=\;&\**\nu(1-x_1\kapa1)\V1+\nu\sum_{j=2}^{n+1} z_j\V{j},
\end{array}$$
for some coefficients $z_j$, which depend on $x_k$, $1\le k\le n$,  and $\kapa m$, for $2\le m\le n$. Hence,
$$\norma{\dparcial G {x_1}\wedge\dparcial G {x_2}\wedge\cdots\wedge\dparcial G {x_n}\wedge\dparcial G t}=\nu|1-x_1\kapa1|\,\norma{\V2 \wedge\V3\wedge\cdots\wedge\V {n+1}\wedge\V1}=\nu|1-x_1\kapa1|$$
and, thus, the absolute value of the jacobian determinant of $G$ equals $\nu|1-x_1\kapa1|$. In particular, fixed $t\in J$, and letting $P= (0,0,\ldots,0,t)$, we obtain that  $|\det JG (P)|=\nu>0$. Since, $f$ is injective, it follows that $G$ is injective on the compact set $K=\{(0,0,\ldots,0)\}\times [a,b]$. Hence, by using the generalized version of the inverse function theorem, we get an open set $U\supset K$ and an open set  $V\supset G(K)=f([a,b])$ such that the restriction \funcao{G}{U}{V} is a diffeomorphism. Since $U$ is an open set, there is a family of small open cylinders $B^{[n]}(\eps_\lambda)\times (a_\lambda,b_\lambda)$ whose union covers $K$. From\break

\includegraphics[scale=0.6]{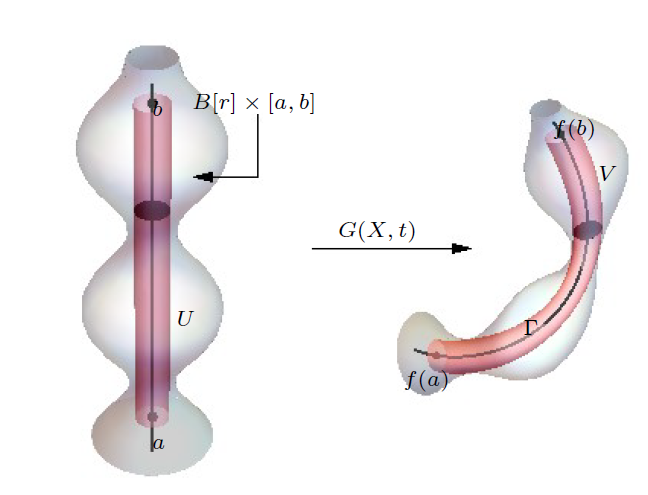}
this family, we get a finite number of such cylinder, say  $B^{[n]}(\eps_j)\times (a_j,b_j)$, $1\le j\le m$, such that $K\subset\cup_{j=1}^m\left(B^{[n]}(\eps_j)\times (a_j,b_j)\right)$.  Now, $\Gamma=G(B^{[n]}[r]\times[a,b])$, $r<\min\{\eps_j\}$, becomes a genuine compact disk tube.\qed 
\brem Many times, we have curves \funcao f{J}{\R n} such that its trace is contained in an affine subspace, say, $W=X_0+S$, where $X_0$ is a point and $S$ is a $m$-dimensional subspace of $\R n$, $m<n$ (think of $f$ as a plane curve in $\R3$). In many of these cases, when $f$ is $(m-1)$-regular, we may construct a partial Frenet apparatus for $f$, obtaining $$\cal A=\{\kapa1>0,\kapa2>0,\ldots,\kapa{m-1}>0,\kapa{m}=0,\V1,\V2,\ldots,\V m,\V {m+1}\},$$ which satisfies the Frenet equations and,  moreover, for each $t\in J$, the subspace generated by $\{\V1,\V2,\ldots,\V m\}$ equals $S$ and $\V{m+1}\in S^{\perp}$ is constant. Then, we consider an orthonormal basis for the orthogonal complement of $S$, say, $\{\V{m+1},\V{m+2},\ldots,\V{n}\}$ and define  $\kapa{m+1}$, \ldots, $\kapa {n-1}$ to be all zero. So, we have a full Frenet apparatus:
$$\cal A=\{\funcoesst \kappa,\framest V\},$$
where $\{\V{m+1},\V{m+2},\ldots,\V{n}\}$ is constant and $\kapa j=0$, $m\le j\le n-1$. Moreover, the Frenet equations remain true for this extended apparatus.

The apparatus of the straight line $f(t)=X_0+ t\, V_1$, $\norma{\V1}=1$, is now defined and is given by
$$\cal A=\{\funcoesst \kappa,\framest V\},$$
where the curvatures are zero and $\{\framest V\}$ is an orthonormal basis of  $\R n$. 
Also, we may consider the regular curves in $\R m$, $m<n$, as curves in $\R n$, by filling with zeros after the $m$-th coordinate: $f(t)=(\puplalinhacomtsp fm,0,0,\ldots,0)$, case in which the Frenet apparatus becomes
$$\cal A=\{\funcoesst \kappa,\framest V\},$$
where $\{\V{m+1},\V{m+2},\ldots,\V{n}\}$ is part of the canonical basis of $\R n$ and $\kapa j=0$, $m\le j\le n-1$. Of course $\V j$, $1\le j\le m$, are those vector fields of the original curve with zeros on the last $n-m$ coordinates. For example, the curve in $\R4$, $f(t)=(\cos t, \sin t, 0,0)$, obtained from $f(t)=(\cos t, \sin t)$, has its Frenet apparatus with curvatures $\kapa1=1$, $\kapa2=0$, $\kapa3=0$ and the Frenet frame
$$\cal F=\{\V1=(-\sin t,\cos t,0,0),\V2=(-\cos t,- \sin t, 0,0),\V3=(0,0,1,0),\V4=(0,0,0,1)\}.$$
Finally, note that the Theorem~\ref{tuboteo} may be applied to any regular curve that has a well defined Frenet apparatus. 
\erem

At the moment that  we have a new theorem, it is always good to seek an example to highlight its statement. We do not find in the literature any non-trivial example of a genuine tube, even of disks, other than the torus of revolution. For this reason, we will try to expose such an example including a significant level of details. Much of the calculations will be omitted, but they will be indicated and left as exercises for  the reader.  We will construct a maximal genuine disk tube around the circular helix $f(t)=(\cos t,\sin t,t)$, $t\in\R{}$.  The example will be decomposed in steps, some of them including some interesting lemmas of real analysis. Actually, we have here a surprisingly laborious and elegant example. 
\bexe{} \label{mainexe}Let \funcao{f}{\R{}}{\R 3}, $f(t)=(\cos t,\sin t,t)$. The apparatus of $f$ is 
\vskip10pt
\begin{tasks}[style=myenumerateexe,label-format={\normalfont\bf},label-width={5mm}, item-indent={25pt}, item-format={\sl},label-align=right, after-item-skip={10pt}, column-sep={0pt}](2)
	\task $\kapa {}=\kapa 1=\frac 12$;
	\task $\tau= \kapa 2=\frac 12$;
	\task $\Tb=V_1=\parl -\frac{\sin (t)}{\sqrt{2}},\frac{\cos (t)}{\sqrt{2}},\frac{1}{\sqrt{2}}\parr$;
	\task $\Nb=V_2=(-\cos (t),-\sin (t),0)$;
	\task $\Bb=V_3=\parl\frac{\sin (t)}{\sqrt{2}},-\frac{\cos (t)}{\sqrt{2}},\frac{1}{\sqrt{2}}\parr$.
\end{tasks}
\vskip10pt
We will study the disk tube of radius 2 around $f$ and centered at $P=(0,0,t)$, which will be indicated by $\Gamma^{\ast}=G^{\ast}(B^{[2]}[2]\times\R{})$, where $B^{[2]}[2]\subset\R2$ is the  compact disk of radius 2 centered at origin and  $G^{\ast}$ from $\R 3$ into $\R3$ is 
$$\arraycolsep=0pt\def\arraystretch{2}
\begin{array}{rcl}
\** G^{\ast}(X,t)&\;=\;&\**f(t)-x\, \Nb(t)+y\, \Bb(t)\\
\**&\;=\;&\**\parl(x+1) \cos (t)+\frac{y\, \sin (t)}{\sqrt{2}},(x+1) \sin (t)-\frac{y\, \cos (t)}{\sqrt{2}},t+\frac{y}{\sqrt{2}}\parr,
\end{array}$$
where $X=(x,y)$. Looking at the map $G$ in (\ref{mapG}), we see that in the rule above we consider $-\Nb(t)$ instead $\Nb(t)$ (compare to the Remark~\ref{remtubes}). The reason is that this simplifies the  handling of certain angles, mainly $\Ti$ and $T$ (see the next picture). The disk tubes obtained from $G^{\ast}$ and $G$ coincide, since $G^{\ast}(-x,y,t)$ equals $G(x,y,t)=f(t)+x \Nb(t)+y \Bb(t)$, which is the map in~(\ref{mapG}). Also, it is easy to see that the injectivity of one of them implies the injectivity of the other one. In what follows, we will conclude that $G^{\ast}$ (and thus $G$) is injective in $B^{[2]}[2]\times \R{}$ and that this fails in $B^{[2]}(R)\times\R{}$ for any open disk of radius $R>2$, $B^{[2]}(R)$. So,  suppose that $G^{\ast}(A,t_1)=G^{\ast}(B,t_2)$, where $A=(x_1,y_1)$ and $B=(x_2,y_2)$ lie in $B^{[2]}[2]$ and, without loss of generality, $t_1\le t_2$. We are going to our first claim.\\[10pt]
\claim: 
\vskip10pt
\begin{tasks}[style=myenumerateexe,label-format={\normalfont\bf},label-width={5mm}, item-indent={25pt}, item-format={\sl},label-align=right, after-item-skip={10pt}, column-sep={0pt}](2)
	\task $t_2-t_1=\frac{y_1}{\sqrt{2}}-\frac{y_2}{\sqrt{2}}$;
	\task $\left(x_1+1\right)^2+\frac{y_1^2}{2}=\left(x_2+1\right)^2+\frac{y_2^2}{2}$.
\end{tasks}
\vskip10pt
In fact, $G^{\ast}(x_1,y_1,t_1)=G^{\ast}(x_2,y_2,t_2)$ is the same as 
\beq \label{Gid}\def\arraystretch{1.5}
\left(
\begin{array}{c}
	\left(x_1+1\right) \cos \left(t_1\right)+\frac{y_1 \sin \left(t_1\right)}{\sqrt{2}} \\
	\left(x_1+1\right) \sin \left(t_1\right)-\frac{y_1 \cos \left(t_1\right)}{\sqrt{2}} \\
	t_1+\frac{y_1}{\sqrt{2}} \\
\end{array}
\right)=\left(
\begin{array}{c}
	\left(x_2+1\right) \cos \left(t_2\right)+\frac{y_2 \sin \left(t_2\right)}{\sqrt{2}} \\
	\left(x_2+1\right) \sin \left(t_2\right)-\frac{y_2 \cos \left(t_2\right)}{\sqrt{2}} \\
	t_2+\frac{y_2}{\sqrt{2}} \\
\end{array}
\right)
\eeq
whence we get easily~[1]. By summing the squares of the two first entries in the matrices above, we get~[2], which shows that $A$ and $B$ lie in the ellipse $\cal E_r$
$$\cal E_r=\{(x,y)\in\R2;\; H(x,y)=\left(x+1\right)^2+\frac{y^2}{2} =r^2\},$$ where $r^2$ equals the identical  real numbers in~[2]. The ellipse $\cal E_r$ plays a crucial role in the discussion that follows. Since  $0=H(-1,0)\le H(x,y)\le 9=H(2,0)$, for all $(x,y)\in B^{[2]}[2]$, and $(2,0)$ is the only maximum point of $H$ on $B^{[2]}[2]$, we see that $0\le r\le 3$ and, moreover, $r=0$ implies that $A=B=(-1,0)$ and $r=3$ yields $A=B=(2,0)$. In any case, [1] gives $t_1=t_2$ and thus $(\dupla1,t_1)=(\dupla2,t_2)$. This is what we want. Hence, we may work with $0<r<3$.  Some observations on the ellipse $\cal E_r$: its semi-minor axis is $r$ and its semi-major axis is $\sqrt 2\, r$; it is contained in the open disk $B^{[2]}(2)$, when $r<1$; when $r=1$, it is contained in the compact disk $B^{[2]}[2]$ and ${\cal E_r}\cap B^{[2]}[2]=\{(-2,0)\}$; if $r>1$, the intersection ${\cal E_r}\cap B^{[2]}[2]$ equals the arc of the $\elipse r$ joining the points $I_1$ and $I_2$, as in the picture below (the angle $\theta$ that appears in it will be considered later). Now, we parametrize $\elipse r$ by $$E(t)=(r\cos (t)-1,r\sqrt2 \sin (t)),\quad t\in[-\pi,\pi].$$

\includegraphics[scale=0.45]{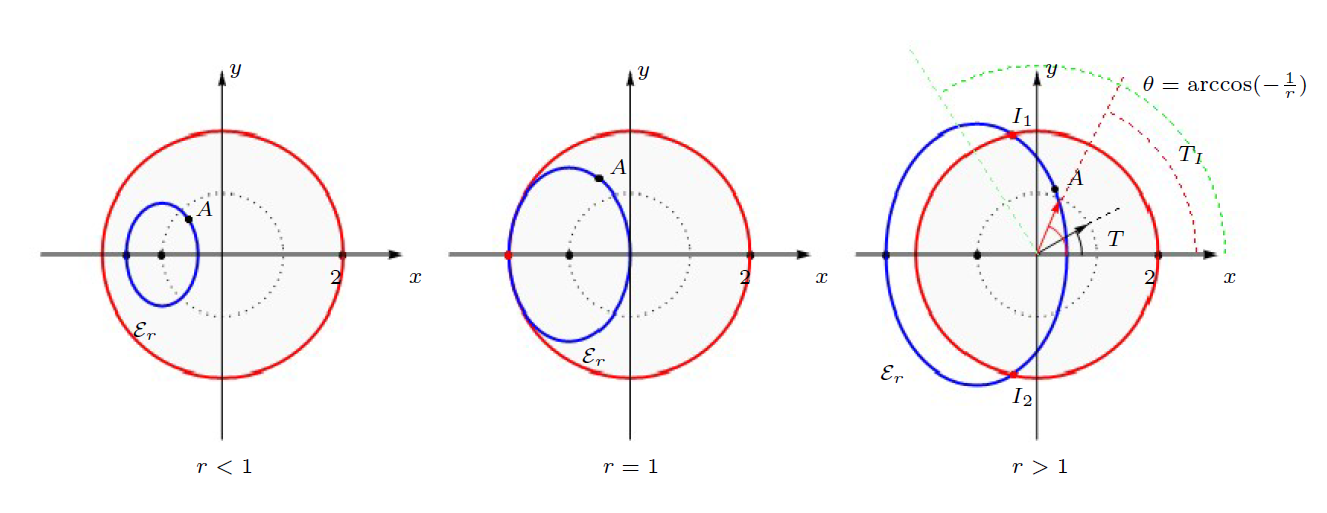}
For the moment, we consider $1<r<3$. A direct calculation shows that
$$I_1=(x_I,y_I)=\parl -2+\sqrt{2} \sqrt{r^2-1},\sqrt{2} \sqrt{1-r^2+2 \sqrt{2} \sqrt{r^2-1}}\parr,$$
and $I_2=(x_I,-y_I)$. We have unique angles $T_I\in (0,\pi)$ and $T\in (-\pi,\pi)$ such that $E(\Ti)=I_1$ and $E(T)=A$ (there exists an analogous angle for $B$ as well). Hence, 
\beq\label{sinnecosTs}
\arraycolsep=10pt\def\arraystretch{3}
\begin{array}{|l|l|}
	\hline
	\**\cos(\Ti)=\frac{1+x_I}{r}=\frac{\sqrt{2} \sqrt{r^2-1}-1}{r} &\**\sin(\Ti)=\frac{y_I}{r\sqrt 2}=\frac{\sqrt{1-r^2+2 \sqrt{2} \sqrt{r^2-1}}}{r}\\\hline
	\**\cos(T)=\frac{x_1+1}{r}&\**\sin(T)=\frac{y_1}{r\sqrt 2}\\
	\hline
\end{array}
\eeq
\vskip10pt\noindent
\claim: $$-\Ti=-\arccos\parl\frac{\sqrt{2} \sqrt{r^2-1}-1}{r}\parr\leq T\le \Ti=\arccos\parl\frac{\sqrt{2} \sqrt{r^2-1}-1}{r}\parr.$$
This claim follows from the discussion above.\\[10pt]
\claim: Without restrictions on $r$, we have that
\beq\label{clidfun}x_1 \sin \left(t_2-t_1\right)+\frac{y_1 \left(\cos \left(t_2-t_1\right)-1\right)}{\sqrt{2}}+t_2-t_1+\sin
\left(t_2-t_1\right)=0.\eeq
In fact, firstly
$$(f(t_2)-G^{\ast}(B,t_2))\cdot f'(t_2)=(-x_2\, \Nb(t_2)+y_2 \, \Bb(t_2))\cdot f'(t_2)=0.$$
On the other hand, since $G^{\ast}(B,t_2)=G^{\ast}(A,t_1)$,
$$\arraycolsep=0pt\def\arraystretch{2}
\begin{array}{rcl}
\**0= (f(t_2)-G^{\ast}(B,t_2))\cdot f'(t_2)
\**&\;=\;&\**(f(t_2)-G^{\ast}(A,t_1))\cdot f'(t_2)\\
\** &\;=\;&\**x_1 \sin \left(t_2-t_1\right)+\frac{y_1 \left(\cos \left(t_2-t_1\right)-1\right)}{\sqrt{2}}+t_2-t_1+\sin
\left(t_2-t_1\right),
\end{array}$$
after some simplifications. Now, we get the fundamental function for our example, obtained from the left side of the claim~3 by putting $s=t_2-t_1\geq 0$ and by using the angle $T$ in the claim~2:
$$g(s)=s-r \sin (T)+r \sin (s+T),\; s\geq 0.$$
In effect,
$$\arraycolsep=0pt\def\arraystretch{2}
\begin{array}{lcl}
\**x_1 \sin \left(t_2-t_1\right)+\frac{y_1 \left(\cos \left(t_2-t_1\right)-1\right)}{\sqrt{2}}+t_2-t_1+\sin
\left(t_2-t_1\right)=&&\\
\quad\quad\quad \**(r\cos(T)-1) \sin \left(s\right)+r\sin(T) \left(\cos \left(s\right)-1\right)+s+\sin
\left(s\right)&\;=\;&\**s-r \sin (T)+r \sin (s+T).
\end{array}$$
Note that $g(0)=0$. The main idea is to show that $s=0$ is the only zero of $g$. Suppose, for a moment, that this is indeed true.  The equation above, together with~(\ref{clidfun}), says that $t_2-t_1$ is a zero of $g$. Hence $t_1=t_2$. By using  [1] of the claim~1 and the equation (3), we deduce that
$$\left(
\begin{array}{c}
\left(x_1+1\right) \cos \left(t_1\right) \\
\left(x_1+1\right) \sin \left(t_1\right) 
\end{array}
\right)=\left(
\begin{array}{c}
\left(x_2+1\right) \cos \left(t_1\right) \\
\left(x_2+1\right) \sin \left(t_1\right)
\end{array}
\right),$$
or
$$((x_1-x_2)\cos(t_1),(x_1-x_2)\sin(t_1))=(0,0).$$
Therefore, $x_1-x_2=0$ and $(\dupla1,t_1)=(\dupla2,t_2)$, as we desire. In general, $g$ can have positive zeros, but with the constraints that come from our problem, this cannot happen and, in fact, $g$ must have only a zero, as we will see below. 

In the cases $0<r\le1$, the problem is very simple. Indeed, since 
$$g'(s)=r \cos (s+T)+1\ge 1-r\ge  0.$$ 
we get that $g$ is increasing and its  critical points are isolated, which implies that $g$ is a strictly increasing function. In particular, $g(s)>0$, for all $s>0$. Hence, $g$ has no positive zero and we are done: $(x_1,y_1,t_1)=(x_2,y_2,t_2)$. We are left with the cases $1<r<3$, which involves many technicalities. Firstly, since $A$ and $B$ belong to the disk $B^{[2]}[2]$, we have that the absolute values of their coordinates are lower than or equal to 2. Hence, by using [1] in the claim 1, we get $s=|s|=|t_2-t_1|\le 2\sqrt 2<\pi$. Therefore, we can consider $g$ defined in $[0,\pi]$:
$$g(s)=s-r\sin (T)+r\sin(s+T),\; 0\le s\le \pi,$$
subject to the constraints
$$-\pi<-\Ti\leq T\le \Ti<\pi,$$
where $\Ti=\arccos(h(r))$,
$$h(r)=\frac{\sqrt{2} \sqrt{r^2-1}-1}{r},$$
and always keeping in mind that 
$$\cos(T)=\frac{x_1+1}{r},\quad\sin(T)=\frac{y_1}{r\sqrt 2},$$
$A=(x_1,y_1)$ and, thus, $r$ and $T$ are all fixed.
Another important date in our discussion is the angle $\theta\in(\frac\pi2,\pi)$ (as shown in the picture above) given by $\theta=\arccos(-\frac 1r)$. Of course $\theta>\Ti$, because $h(r)+\frac1r>0$ and the $\arccos$ is a strictly decreasing function. Our strategy will be to describe the critical points of $g$ and then to show that $g$ is positive at each of these points, which implies that $g>0$ on $(0,\pi]$, according to the next lemma (Lemma~\ref{criticolema}). We are going to describe the critical points of $g$. For this purpose, let $s_c\in(0,\pi)$ be a critical point of $g$. \\[10pt]
\claim: The critical point $s_c$ satisfies only one of the following descriptions:
\vskip10pt\noindent
\begin{tasks}[style=myenumerateexe,label-format={\normalfont\bf},label-width={5mm}, item-indent={25pt}, item-format={\sl},label-align=right, after-item-skip={10pt}, column-sep={0pt}](3)
	\task $s_c+T=\theta$;
	\task $s_c+T=-\theta$;
	\task $s_c+T=-\theta+2\pi$.
\end{tasks}
\vskip10pt\noindent
Indeed, since $g'(s_c)=1+r\cos(s_c+T)=0,$ we get that either $s_c+T=\theta+2k\pi$ or $s_c+T=-\theta+2k\pi$, for some integer $k$. This, together with $\theta\in(\pi/2,\pi)$ and  $s_c+T\in(-\pi,2\pi)$, due to the constraints on $s_c$ and $T$, implies the three items above. \\[10pt]
\claim: $g(s_c)>0$.
\\[10pt]
We will study separately the three cases in the claim~4.
\vskip10pt
{\leftmargini=0pt
	\begin{enumerate}\parskip=10pt
		\item[{\rm\bf [1]}] Here, $s_c+T=\theta$, whence 
		$$g(s_c)=s_c-r \sin (T)+r \sin(\theta)=(\theta +r \sin (\theta ))-(r \sin (T)+T)=\psi(\theta)-\psi(T),$$
		where, 
		\beq\label{psi}\psi(u)=u +r \sin (u).\eeq
		Of course $g(s_c)>0$, if $T\le0$. When $0<T<\Ti<\theta<\pi$, we need a more delicate analysis. It is as follows. It is not hard to see that $\psi$, on the interval $[0,\pi]$, attains its strict global maximum at $u=\theta$, its only critical point. Hence, $g(s_c)=\psi(\theta)-\psi(T)>0$.
		\item[{\rm\bf [2]}]  Now, $s_c+T=-\theta$. Hence,
		$- \pi <s_c + T <-\frac{\pi}2$ and $\cos(s_c+T)=-1/r$. Suppose, by contradiction, that $g(s_c) \le  0$. Since $ g (0) = 0$ and 
		$$g'(0)=1+r\cos (|T|)\ge 1+r\cos(\Ti)=\sqrt{2} \sqrt{r^2-1}>0,$$
		we would have $g (s_1) = 0$, for some $0 <s_1\le s_c$. Therefore, we would have $g '(s_2) = 0$, for some $0 <s_2 <s_c$, hence $\cos(s_c + T) = \cos ( s_2 + T) = - 1 / r$ and $ -\pi <s_2 + T <s_c + T <-\pi / 2$. We have an absurdity, because the $\cos$  is injective in $(-\pi, -\pi / 2)$.
		\item[{\rm\bf [3]}] $s_c+T=-\theta+2\pi$. Actually, this is the harder part of our example. First, note that $T>0$. We have
		$$g\left(s_c\right)=g(-\theta -T+2 \pi )=-\theta-T+2\pi -r \sin
		(T)-r \sin (\theta ) =2 \pi -(\psi (\theta )+\psi (T)).$$
		Now, since $\psi$ (see (6)) is strictly increasing in $(0,\theta)$, we get $-\psi(T)>-\psi(\Ti)$ and thus
		$$g\left(s_c\right)=2 \pi -(\psi (\theta )+\psi (T))>2\pi-(\psi (\theta )+\psi (\Ti))=2\pi-\eta(r),$$
		where
		$$\eta(u)=\left(\arccos\left(-\frac{1}{u}\right)+\sqrt{u^2-1}\right)+\left(\arccos\left(\frac{\sqrt{2} \sqrt{u^2-1}-1}{u}\right)+\sqrt{1-u^2+2 \sqrt{2} \sqrt{u^2-1}}\right),$$
		$u\in[1,3]$, 
		which we can patiently verify the following:
		\begin{tasks}[style=myenumerate,label-format={\normalfont},label-width={5mm}, item-indent={20pt}, item-format={\rm},label-align=right, after-item-skip={4pt}, column-sep={0pt}](2)
			\task $\eta(1)=2\pi$;
			\task $\eta(3)=2 \sqrt{2}+\arccos\left(-\frac{1}{3}\right)<2\pi$;
			\task*(2) $\eta'(u)=0$ if and only if $u=\sqrt 3$;
			\task*(2) $\eta(\sqrt 3)=2 \sqrt{2}+\arccos\left(-\frac{1}{\sqrt{3}}\right)+\arccos\left(\frac{1}{\sqrt{3}}\right)=2 \sqrt{2}+\pi<2\pi,$ for $\arccos(-x)+\arccos(x)$ is constant.
		\end{tasks}
		These properties imply that $u=1$ is the global maximum point for $\eta$. Hence,  $\eta(u)<2\pi$ in the interval $(1,3)$. As a matter of fact, since $\eta(1)<\eta(\sqrt 3)<\eta(3)$, $\eta$ is strictly decreasing in $[1,3]$.  So $g(s_c)>0$. 
	\end{enumerate}
}
We are almost done, because we still need to verify that
$ g (\pi) = \pi -2 r \sin (T)> 0. $
The lemma that we will use requires this fact, which will be left as an exercise. We suggest to consider separately $1<r<\sqrt{\frac 32}$, case in which $\Ti>\frac\pi2$, and $\sqrt{\frac 32}\le r<3$, where $\Ti\le\frac\pi2$. Of course the result is trivial if $T\le 0$.
Putting this all together, we get that $g(0)=0$, $g>0$ at its critical points and $g(\pi)>0$. By using  lemma below, it follows that $g>0$ in $(0,\pi)$. Thus, $g$ vanishes only at $s=0$ and it follows that $(X,t_1)=(Y,t_2)$ (see the main idea that we had talked above), that is, $G^{\ast}$ and $G$ are both injective. 

It remains to show that $\Gamma^{\ast}$ is a maximal genuine tube, that is, if $R>2$ then $G^{\ast}$ is no longer injective. In fact, let $R>2$ and define $r=\frac R2>1$. Then let $s_0$ be a positive solution of $s-r\sin s=0$ (the left side is the function $g$ for $T=\pi$) that exists, according to the~Lemma~\ref{rootsen} bellow. Now put $P_1=(X_1,0)$, $X_1=(0,-1-r)$,  and $P_2=(Y_1,s_0)$, $Y_1=(-1-r\cos(s_0),\sqrt 2\,s_0)$. Of course $P_1\neq P_2$,
$G^{\ast}\left(P_1\right)=(-r,0,0)$
and
$$\arraycolsep=0pt\def\arraystretch{2}
\begin{array}{rcl}
\** G^{\ast}\left(P_2\right)&\;=\;&\**(-r \cos ^2\left(s_0\right)-s_0 \sin \left(s_0\right),\cos \left(s_0\right) \left(s_0-r \sin
\left(s_0\right)\right),0)\\
\**&\;=\;&\**(r \sin ^2\left(s_0\right)-r-s_0 \sin \left(s_0\right),s_0 \cos \left(s_0\right)-r \sin \left(s_0\right) \cos
\left(s_0\right),0)\\
\** &\;=\;&\**\parl r\, \frac{s_0^2}{r^2}-r-s_0 \, \frac{s_0}{r},s_0 \cos \left(s_0\right)-r  \frac{s_0}{r} \cos
\left(s_0\right),0\parr
=(-r,0,0)=G^{\ast}(P_1).
\end{array}$$
Moreover $\norma{X_1}=1+r<R$ and $\norma{Y_1}<R$. This last inequality results of the fact that  $Y_1$ belongs to the ellipse $(x+1)^2+\frac{y^2}2=r^2=\frac{R^2}{4}$ whose intersection with the circle $x^2+y^2=R^2$ is empty, given that the system formed by these last two equations has no solution in $\R2$. Indeed, by working a little harder, we can establish that $$G^{\ast}(0,-1-r,t)=G^{\ast}(-1-r\cos(s_0),\sqrt 2\,s_0,t+s_0)=(-r \cos (t),-r \sin (t),t),$$
for an arbitrary  $t$. Finally, we see below  a piece of the infinite disk tube $\Gamma^{\ast}$, say $\Gamma_1$, together with another genuine tube $\Gamma_2\subset\Gamma_1$, whose cross sections are copies of the region 
\beq\label{cross}S=\left\{(x,y);\; -1\le x\le 1,\;-\frac{x^2}{4}-1\le y\le\frac{x^2}{4}+1\right\}
\eeq

\includegraphics[scale=0.55]{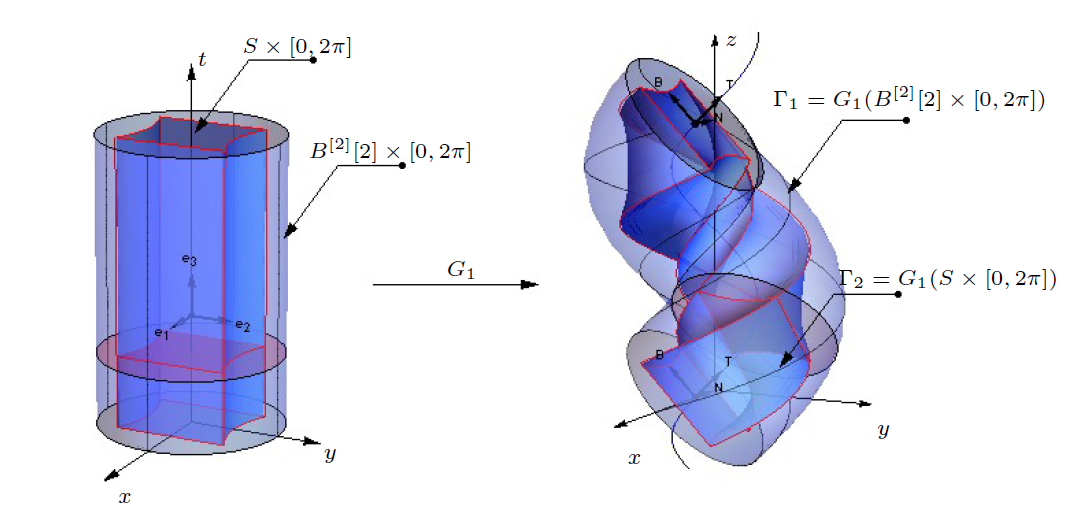}
The tubes $\Gamma_1$ and $\Gamma_2$ are obtained for $t$ running from $0$ to $2\pi$. In the Example~\ref{volumeexe}, we calculate their volumes.
\eexe
\blema \label{criticolema} Let \funcao{g}{ [0, b] }{\R{}} be a differentiable  function such that $g (0)\ge 0$ and $g (b)> 0$. If either $g$ has no critical point in $(0,b)$ or $g$ is positive at each of its critical points, then $g>0$ in $(0,b]$.
\elema
\proof Suppose, by contradiction, that $g(x_1)\le 0$ for some $x_1\in(0,b)$. Hence, there exists $x_0\in[x_1,b)$ such the $g(x_0)=0$, since $g(b)>0$.  By using the Rolle's theorem, if $g(0)=0$, we get $x_c\in(0,x_0)$ such that $g'(x_c)=0$. If $g(0)>0$, the minimum value of $g$ in $[0,b]$ is smaller than or equal to zero, because $g(x_0)=0$ and thus this minimum value must occur at some  $x_c\in (0,b)$, where $g'$ vanishes. Hence, $g$ must be positive in $(0,b]$, if $g$ has no critical point.  Now, we suppose that $g(x_c)>0$, whenever $g'(x_c)=0$. We have that  $g'(x_0)\neq 0$, by hypothesis. If $g'(x_0)<0$, we choose $x_1>x_0$ where $g$ is negative. This implies that the  minimum value of $g$ in $[x_0,b]$ is negative and occurs in $(x_0,b)$, say at $x_c$. Hence $g'(x_c)=0$ and $g(x_c)<0$, a contradiction. The case $g'(x_0)>0$ leads to a similar contradiction, now by taking  $x_1<x_0$ where $g$ is negative.\qed
\blema \label{rootsen} If $r>1$, then there exists $s_0\in(0,\pi)$ such that $s_0-r\sin s_0=0$.
\elema
\proof Define $g(s)=s-r\sen s$, $s\in \R{}$. If $s\neq 0$, $g(s)=s\parl 1-r\frac{\sin s}s\parr$ and,  since $$\lim_{s\rightarrow 0}\parl 1-r\frac{\sin s}s\parr=1-r<0,$$ we can choose $s_1>0$ near to zero such that $g(s_1)<0$. But $g(\pi)=\pi>0$. From the intermediate value theorem, we obtain $s_0$ between $s_1$ and $\pi$ such that $g(s_0)=0$, as we wanted.\qed
\subsection{On the Volume of Tubes}

The calculation of the volume of disk tubes in $\R {n+1}$ and in the euclidean sphere ${\mathbb S}^n$ appears in~\cite{Hot}.  Actually, the main goal of this paper is to work within the sphere. In~\cite{Good}, the author generalizes the  Pappus theorems for genuine tubes around an arbitrary curve in $\R 3$.  In the case of the second Pappus theorem, his proof makes use of the divergence theorem. Here, we will establish such a theorem for genuine tubes in $\R{n+1}$ by using a different strategy, namely, the formula for change variable in multiple integrals.

Next, we will consider an interval $I=[a,b]$ and the tube $\Gamma=\Gamma_P(f, \Omega_I(S))$  generated by the solid  $\Omega_I(S_t)$ around the curve \funcao f{J\supset I}{\R{n+1}} and centered at $P(t)=(\til P,t)$, where $S_t\subset\R n\times\{t\}$ has positive volume and $\til P=\pupla pn$ is constant. In what follows, we assume this to be the case. Thus, the cross sections of $\Gamma$ are stuck to $f$ by the same point. We have that $\Gamma$ is the image of $\Omega_I(S_t)$ under
\funcao{G}{\R n \times J}{\R{n+1}}  given  by 
\beq
G(X,t)=G(\puplasp xn,t) =f(t)+\sum_{j=1}^n (x_j-p_j)\V{j+1}(t).
\eeq

\includegraphics[scale=0.65]{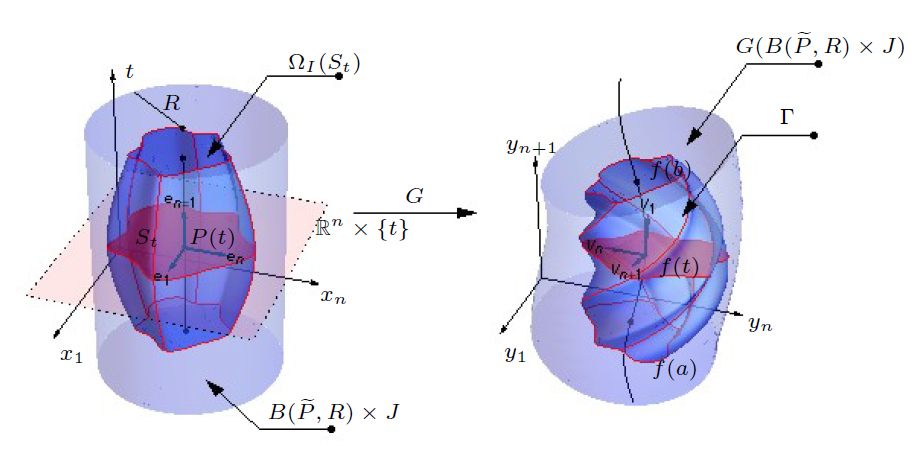}
Moreover, exactly as in the proof of the Theorem~\ref{tuboteo}, we have that 
$$|\det JG(\puplasp xn,t)|=\nu(t)|1-(x_1-p_1)\kapa 1(t)|.$$
Let $R\in \R{}$  be a positive number smaller than the {\sl radius of curvature} $\rho=\frac1{\sup_J \kapa 1(t)}$.  Let $B^{[n]}(\til P,R)$ be the open disk centered at $\til P$ of radius $R$ in $\R n$. Then, for all $X=\pupla xn\in  B^{[n]}(\til P,R)$ and $t\in J$, we have that
\beq \label{regularitycond}|\det JG(\puplasp xn,t)|=\nu(t)|1-(x_1-p_1)\kapa 1(t)|=\nu(t)(1-(x_1-p_1)\kapa 1(t))>0,
\eeq
because $$(x_1-p_1)\kapa1\le|x_1-p_1|\sup_J\kapa1\le\norma{X-\til P}\,\sup_J\kapa1<\rho\sup_J\kapa1=1.$$ 
We refre to the inequality  in \refeq{regularitycond} as  {\sl fundamental regularity condition}.
\brem
The formula above is obtained from the default map $G$ in \refeq{mapG}. For example, if we use the vector field $-\V2(t)$ instead $\V2(t)$ and consider the map
$$G^{\ast}(X,t)=f(t)-(x_1-p_1(t))\V2(t)+\sum_{j=2}^n (x_j-p_j(t))\V{j+1}(t),$$
the absolute value of its jacobian  becomes 
$$|\det JG^{\ast}(\puplasp xn,t)|=\nu(t)|1+(x_1-p_1)\kapa 1(t)|.$$
\erem
Now, suppose that $\Omega_I(S_t)\subset B^{[n]}(\til P,R)\times J$, and that $G$ is injective in $B^{[n]}(\til P,R)\times J$.  By using the  change of  variables  formula for multiple integrals and the Fubini theorem, it follows that
\beq \label{voltuboform}\arraycolsep=0pt\def\arraystretch{3}
\begin{array}{rcl}
	\** \vol\, \Gamma&\;=\;&\**\int_{\Gamma}\dY1\dY2\ldots\dY{ n+1}=\int_{G(\Omega_I(S_t))}\dY1\dY2\ldots\dY{ n+1}\\
	\**&\;=\;&\**\int_{\Omega_I(S_t)}|\det JG|\dx1\dx2\ldots\dx n\dt\\
	\**&\;=\;&\**\int_a^b\nu(t)\left(\int_{S_t}(1-(x_1-p_1)\kapa 1(t))\dx1\dx2\ldots\dx n\right)\dt\\
	\**&\;=\;&\**\int_a^b\nu(t)(\vol\,S_t)\dt-\int_a^b\nu(t)\kapa 1(t)\left(\int_{S_t}x_1\dx1\dx2\ldots\dx n\right)\dt+\\
	&&\**\quad\quad\quad\quad\quad\quad\quad\quad\quad{}+p_1\int_a^b\nu(t)\kapa 1(t)(\vol\,S_t)\dt,
\end{array}
\eeq
\mbox{}\\[5pt]
From this, we will get a series of interesting results on the volume of the tubes. For instance, under a reasonable condition on the injectivity of $G$,  $\vol\, \Gamma$  does not change if we replace the $(n-1)$ last coordinates of $\til P$ by any $(n-1)$-tuple. This means that we can shift the axis of the tube orthogonally to $\V2$ without altering the volume of the tube. 

Now, we give a definition.
Given a compact set $S\subset\R n$ with $\vol\, S>0$, we define the {\sl barycenter of $S$} to be the point ${\rm Center}(S)=\pupla c n$, where
\beq\label{bary}
c_j=\frac1{\int_S\dvol n}\int_Sx_j\dvol n=\frac1{\vol\, S}\int_Sx_j\dvol n,\quad 1\le j\le n.
\eeq
In some situations we use ${\rm Center}(S,j)$ to mean $c_j$. 

By using this definition,~(\ref{voltuboform}) becomes 
$$
\vol\,\Gamma=\int_a^b\nu(t)\,\vol\,S_t\dt-\int_a^b\nu(t)\kapa 1(t)\,(\vol\,S_t)\,{\rm Center}(S_t,1)\dt+p_1\int_a^b\nu(t)\kapa 1(t)\,\vol\,S_t\dt.
$$
In particular, if $S_t=S$ is constant,
\beq\label{voltuboform2}
\vol\,\Gamma=(\vol\,S)\,l(f,I)+(\vol\,S)\left(p_1-{\rm Center}(S,1)\right)\int_a^b\nu(t)\kapa 1(t)\dt,
\eeq
where $l(f,I)$ is the length of arc of the curve $f$  from $t=a$ to $t=b$.
\bprop \label{barydisk}Given $A=\pupla an$ and $R>0$,  the barycenter of the compact disk centered at $A$ and of radius $R$, $B^{[n]}[A,R]\subset\R n$, is its center.
\eprop 
\proof We can suppose without loss of generality that $A$ is the origin of $\R n$. Let $S=B^{[n]}[R]$. From
$$\int_Sx_n\dvol n=\int_{B^{n-1}[R]}\left(\int_{-\sqrt{R^2-\sum_{i=1}^{n-1}x_i^2}}^{\sqrt{R^2-\sum_{i=1}^{n-1}x_i^2}}x_n\dx n\right)\dvol{n-1}=0,$$
we get that ${\rm Center}(S,n)=0$. Now, for $1\le j\le n-1$,
$$\arraycolsep=0pt\def\arraystretch{3}
\begin{array}{rcl}
\** \int_Sx_j\dvol n&\;=\;&\**\int_{-R}^R\left(\int_{B^{n-1}[\sqrt{R^2-x_n^2}]}x_j\dvol{n-1}\right)\dx n\\
\**&\;=\;&\**\int_{-R}^R\left(\vol (B^{n-1}[\sqrt{R^2-x_n^2}])\,{\rm Center}(B^{n-1}[\sqrt{R^2-x_n^2}],j)\right)\dx n\\
\**&\;=\;&\**\int_{-R}^R0\dx n=0,
\end{array}$$
by using induction on $n$. Hence,  ${\rm Center}(S)$ equals the origin of $\R n$.\qed
\vskip10pt

The second Pappus theorem, as the original version as its generalization in ~\cite{Good}, deals with the calculation of volumes of tubes with constant cross sections and axes passing trough the barycenters of these sections. In our case this means that $S_t=S$, $\Omega_I(S_t)=S\times I$, $\til P={\rm Center}(S)$ and $\Gamma=G(S\times I)$, as in the picture bellow, where we are taking account all discussion just before the equation~(\ref{voltuboform}). For such a genuine tube, we have the following theorem that generalizes the second Pappus theorem for tubes in $\R{n+1}$.

\includegraphics[scale=0.6]{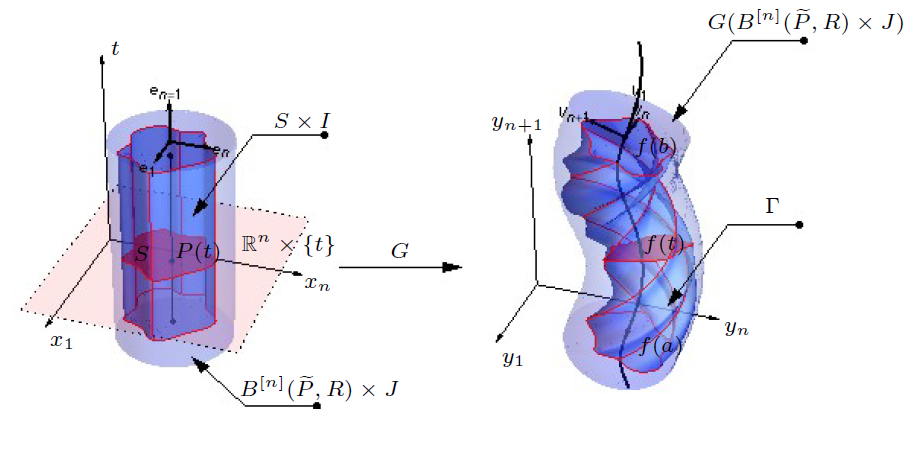}

\bteo \label{pappusteo}$\vol\,\Gamma=l(f,I)\,\vol\, S$, where $l(f,I)$ is the length of the curve $f$ from $t=a$ to $t=b$. In other words, the volume of the tube is the product of the length of the curve traveled by the barycenter of $S$ and the volume of the cross section $S$.
\eteo 
\proof Initially, note that $p_1={\rm Center}(S,1)$. By replacing this fact in~(\ref{voltuboform2}), we obtain easily that
$
\vol\,\Gamma=(\vol\,S)\,l(f,I),
$
as we wanted.\qed
\vskip10pt

Another application of the equation~\refeq{voltuboform2} is the calculation of volumes of the  traditional and genuine disk tubes  with variable radius, say $R(t)$, that is, for each $t\in J$, the cross section $S_t$ equals disk $B^{[n]}[R(t)]$ and  $P(t)=(0,\ldots,0,t)$ or $\til P$ is the center of the disk $B^{[n]}[R(t)]$. In this situation, we have 
$p_1=0$ and ${\rm Center}(S_t,1)=0$ which implies that
\beq\label{voltubedisk}\vol\,\Gamma=\int_a^b\nu(t)\,\vol\,(B^{[n]}[R(t)])\dt=\vol\,(B^{[n]})\int_a^b\nu(t)(R(t))^n\dt,\eeq
where $\vol\,(B^{[n]})$ is the volume of the compact unit disk which can be recursively obtained from
$$\vol\,(B^{[n]})=\vol\,(B^{[n-1]})\int_{-1}^1(1-x^2)^{\frac{n-1}{2}}\dx{},$$
where $\vol\,(B^{[1]})=\vol\,([-1,1])=2$,
as can be seen in~\cite{Flanders} (page~74). In particular, we get the volume of the revolution solid $\Gamma_\alfa$ generated by the plane curve $\alfa(t)=(R(t),0,0,\ldots,0,t)\in\R{n+1}$, $t\in[a,b]$, around the $y_{n+1}$-axis, namely,
$$\vol\, \Gamma_\alfa=\vol\,(B^{[n]})\int_a^b(R(t))^n\dt,$$
because, in this case, $f(t)=(0,0,\ldots,0,t)$ an thus $\nu=1$. This
result  appears in~\cite{Aberra} (equation~(2.2)).
\bexe \label{volumeexe}Consider the tubes $\Gamma_1$ and $\Gamma_2$ in the Example~\ref{mainexe}. $\Gamma_1$ is a disk tube of constant radius, namely $R(t)=2$. The speed $\nu=\sqrt 2$. Thus its volume is
$$\vol\,(\Gamma_1)=\area(B^{[2]}[2])\,l(f,[0,2\pi])=4\pi \int_0^{2\pi}\sqrt 2\dt=8 \sqrt{2} \pi ^2.$$
For $\Gamma_2$, at first, note that ${\rm Center} (S)=(0,0)$ and that $\area\, S=2 \int_{-1}^1 \left(\frac{x^2}{4}+1\right)  \dx{}=\frac{13}{3}$. Hence 
$$\vol\,(\Gamma_2)=\area(S)\,l(f,[0,2\pi])=\frac{26 \sqrt{2} \pi }{3}.$$
In both cases,  we have used the Theorem~\ref{pappusteo}, because the axis of the tube (the circular helix)  passes trough the barycenter of the cross section, that is, $\til P={\rm Center}(S)$. In the examples that we have discussed so far, $\til P$ has been chosen to be the barycenter of the cross section $S$. We now make a more instructive example. Indeed, consider the same cross section $S$ as that of the tube $\Gamma_2$ (see~\refeq{cross}) and $\til P=(\frac12,0)$. Now, it is exactly at this  point that the circular helix $f(t)=(\cos t,\sin t,t)$, $t\in[0,2\pi]$, will meet $S$ to form the tube~$\Gamma_3$. We have that  $p_1=\frac12$ and ${\rm Center}(S,1)=0$. A full geometric description of this case is in the picture bellow, where $\Gamma_3=G(S\times[0,2\pi])$ and the default map $G$ is given by
$$G(x,y,t)=f(t)+\parl x-\frac12\parr\,\Nb(t)+y\,\Bb(t),$$
that is injective in $B^{[2]}[\til P,2]\times \R{}$ by the result in the Example~\ref{mainexe}. Since $S\subset B^{[2]}[\til P,2]$, it follows that $\Gamma_3$ also is  genuine and its volume is
$$\arraycolsep=0pt\def\arraystretch{3}
\begin{array}{rcl}
\** \vol\,(\Gamma_3)&\;=\;&\**(\area\,S)\,l(f,I)+(\area\,S)\left(p_1-{\rm Center}(S,1)\right)\int_a^b\nu(t)\kapa 1(t)\dt\\
\**&\;=\;&\**\frac{13}{3}\,2\pi\sqrt2+\frac{13}{3}\left(\frac12-0\right)2\pi\,\sqrt2\,\frac12
=\frac{26}{3}\,\pi\sqrt2+\frac{13 \pi \sqrt{2} }{6}=\frac{65 \pi \sqrt{2} }{6}.
\end{array}$$

\includegraphics[scale=0.42]{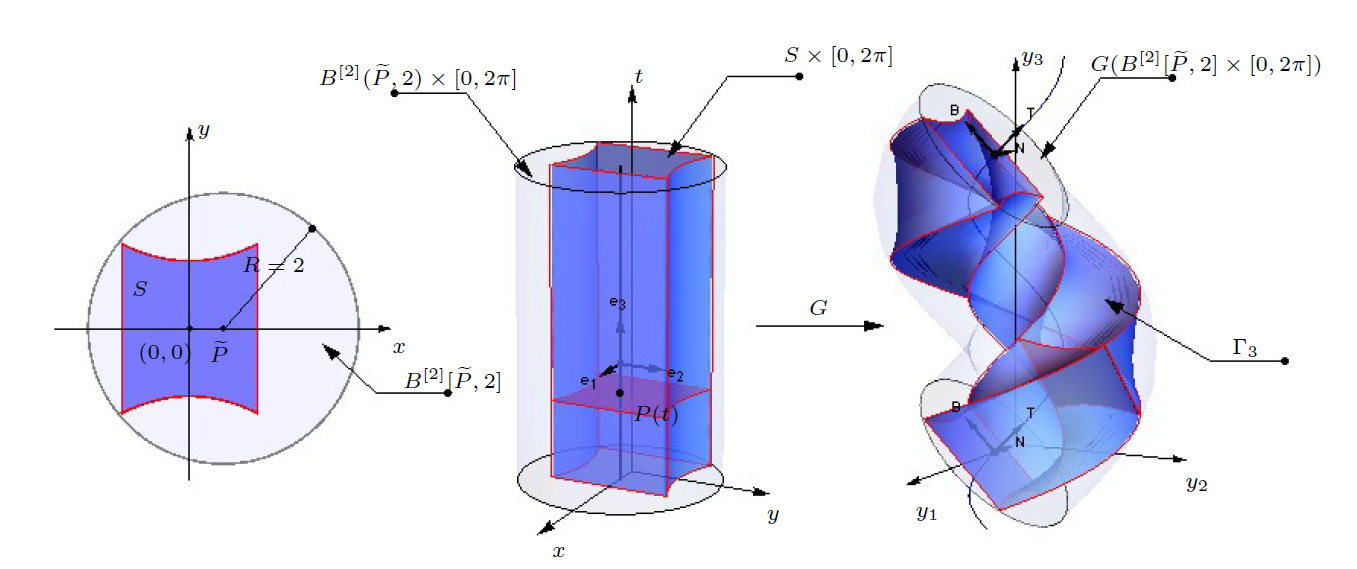}
Observe that the disk tube $G(B^{[2]}[\til P,2]\times[0,2\pi])$ coincides with $\Gamma_1$, that was obtained with another map in the Example~\ref{mainexe}.
For completeness, let us include one more example, namely, a disk tube of variable radius around the same arc of the circular helix $f$. For this, consider $R(t)=1+\frac12\sin t$ and $\Gamma_4$ the disk tube of radius $R(t)$, $0\le t\le2\pi$. This is to be to say that the cross sections of $\Gamma_4$ are the disks $S_t=B^{[2]}[R(t)]$. Since $0<R(t)<2$, $\Gamma_4=G(\Omega_{[0.2\pi]}S_t)$ is genuine, where, of course, $G$ must be changed to
$$G(x,y,t)=f(t)+x\,\Nb(t)+y\,\Bb(t).$$ By using~\refeq{voltubedisk}, we get that
$$\vol\,(\Gamma_4)=\area\,(B^{[2]})\int_a^b\nu(t)(R(t))^2\dt=\pi  \int_0^{2 \pi } \sqrt{2} \left(\frac{\sin (t)}{2}+1\right)^2 \dt=\frac{9 \pi
	^2}{2 \sqrt{2}}.$$

\begin{center}
	\includegraphics[scale=0.5]{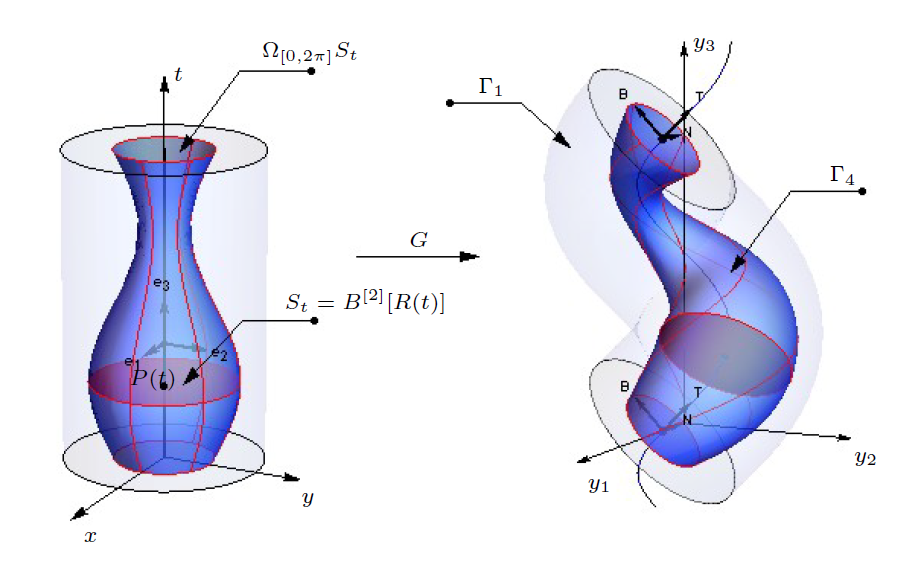}
\end{center}

\eexe
\subsection{On the Volume of Sphere Tubes} We start with a genuine compact disk tube $\Gamma\subset\R{n+1}$ of variable radius, say $R(t)>0$, and centered at the centers of its cross sections (see the picture above). So, we have a compact solid of revolution around the $t$-axis
$$\Omega_I=\Omega_I(B^{[n]}[R(t)])=\{(X,t);\;\norma{X}\leq R(t),\; X=\pupla xn,\; t\in I=[a,b]\},$$
a curve $f$, with Frenet frame  $\{\reffrenetnew{n+1}\}$, defined in $J\supset I$ and
the default map $G$,
$$G(X,t)=f(t)+\sum_{j=1}^n x_j\V{j+1}(t).$$
Next, we suppose that $G$ is injective, satisfies the fundamental regularity condition~\refeq{regularitycond} in some neighborhood of the $\Omega_I$
and that  $\Gamma=G(\Omega_I)$. Denoting by $\esfera{n-1}$ the unit sphere of  $\R n$, we see that the boundary of $\Omega_I$ is given by
$$\partial\omegar=(B^{[n]}[R(a)],a)\cup \partial_l\omegar\cup (B^{[n]}[R(b)],b),$$
where $$\partial_l\omegar=\{(R(t)\,X,t);\;X\in\,\esfera{n-1},\; X=\pupla xn,\; t\in I=[a,b]\}$$
is what we call the {\bf lateral hypersurface} of $\omegar$. The image of this hypersurface of revolution under $G$ will be called {\sl lateral hypersurface} of $\Gamma$ and will be denoted by $\partial_l\Gamma$.  It is a {\sl sphere tube}. Thus, if 
$$\fifi\pupla u {n-1}=(\fifi_1\pupla u{n-1},\fifi_2\pupla u{n-1},\ldots,\fifi_n\pupla u{n-1})$$
is a parametrization for $\esfera{n-1}$, which we may take as an orthogonal parametrization,  then
\beq\label{bordopar}
H(U,t)=G(\fifi(U),t)=f(t)+R(t)\sum_{j=1}^n \fifi_j(U)\V{j+1}(t),
\eeq
$U=\pupla u {n-1}$ ranging over some region of $\R{n-1}$, describes  the hypersurface $\bgamma$.
Our main goal in this section is to get the volume (or area) element of $\bgamma$ associated to $H$, that is, 
$$\dbgamma=\sqrt{\det (g_{ij})}\,\du1\du2\ldots\du{n-1}\dt,$$
where $(g_{ij})$ is the symmetric matrix of the first fundamental form of $\bgamma$ with respect to $H$. More precisely,
$g_{ij}=\dparcial H{u_i}\cdot\dparcial H{u_j}$, $1\le i, j \le n-1$, $g_{in}=\dparcial H{u_i}\cdot\dparcial Y{t}$, $1\le i\le n-1$, and $g_{nn}=\dparcial H{t}\cdot\dparcial H{t}$. We remember that
$$\sqrt{\det (g_{ij})}=\norma{\dparcial H{u_1}\wedge\dparcial H{u_2}\wedge\cdots\wedge\dparcial H{u_{n-1}}\wedge\dparcial H{t}}=\norma{\dparcial H{u_1}\times\dparcial H{u_2}\times\cdots\times\dparcial H{u_{n-1}}\times\dparcial H{t}}$$
and that $$N(U,t)=N(\puplasp u{n-1},t)=\frac{\dparcial H{u_1}\times\dparcial H{u_2}\times\cdots\times\dparcial H{u_{n-1}}\times\dparcial H{t}}{\sqrt{\det (g_{ij})}}$$ is the Gauss map of $\bgamma$.
The next Lemma will be useful for our objective.\def\dHd#1{\dparcial H{#1}}
\blema \label{matrideHlem}Let $D$ be the following $(n\times n)$ ``almost'' diagonal symmetric matrix
$$D=\left(
\begin{array}{ccccc}
a_1 & 0 & \ldots & 0 & b_1 \\
0 & a_2 & 0 & \ldots & b_2 \\
0 & 0 & \ddots & 0 &\vdots \\
0 & 0 & 0 & a_{n-1} & b_{n-1} \\
b_1 & b_2 &\ldots & b_{n-1} & a_n \\
\end{array}
\right),$$
where $a_j\neq0$, for all $j$. Then 
$$\det\,D=a_1a_2\cdots a_{n-1}\left(a_n-\left(\frac{b_1^2}{a_1}+\frac{b_2^2}{a_2}+\cdots +\frac{b_{n-1}^2}{a_{n-1}}\right)\right).$$
\elema
\proof After some elementary operations on the rows of $D$, we obtain that
$$\det D=a_1a_2\cdots a_{n-1}\det\left(
\begin{array}{ccccc}
1 & 0 & \ldots & 0 & \frac{ b_1}{a_1} \\[10pt]
0 & 1 & 0 & \ldots & \frac{ b_2}{a_2} \\[10pt]
0 & 0 & \ddots & 0 &\vdots \\[10pt]
0 & 0 & 0 & 1 &  \frac{ b_{n-1}}{a_{n-1}} \\[10pt]
0 & 0 &\ldots &0 & a_n-\frac{b_1^2}{a_1}-\frac{b_2^2}{a_2}-\cdots -\frac{b_{n-1}^2}{a_{n-1}} \\
\end{array}
\right),$$
whence the result follows easily.\qed
\vskip10pt

Now, let us return to the function $H$ in ~\refeq{bordopar}, where we consider $\fifi$ as an  orthogonal parametrization of $\esfera{n-1}$, that is, $\dparcial{\fifi}{u_i}\cdot\dparcial{\fifi}{u_j}=0$, whenever $i\neq j$. In other words,  given $U=\puplacomdois u{n-1}$, the set 
$$\left\{\dparcial{\fifi}{u_1},\dparcial{\fifi}{u_2},\ldots,\dparcial{\fifi}{u_{n-1}}\right\}$$
is an orthogonal basis of the tangent space of $\esfera{n-1}$ at $\fifi(U)$. Moreover, since $\fifi\cdot\fifi=1$, it follows that $\fifi\cdot\dfi{}{u_j}=0$, $1\le j\le {n-1}$. Now, we will study the properties of  $H$.\vskip10pt
\bprop{} \label{Hprop}The parametrization $H$ of $\bgamma$ has the following properties, where the  $H_0$ stands  for the vector $H(U,t)-f(t)$.\mbox{}\\
\begin{tasks}[style=myenumerate,label-format={\normalfont},label-width={5mm}, item-indent={25pt}, item-format={\rm},label-align=right, after-item-skip={10pt}, column-sep={-30pt}](2)
	\task $\norma{H_0}^2=(R(t))^2$;
	\task $\dparcial H{u_i}\cdot\dparcial H{u_j}=0$, $i\neq j$, $1\le i, j \le n-1$;
	\task $H_0\cdot\dparcial H{u_j}=0$, $i\neq j$, $1\le i, j \le n-1$;
	\task $H_0\cdot\dparcial H{t}=R(t)R'(t)$;
	\task*(2) 
	${\rm span}\left\{H_0,\dparcial H{u_1},\ldots,\dparcial H{u_{n-1}}\right\}= {\rm span}\{\V2(t),\V3(t),\ldots,\V{n+1}(t)\};$
	\task*(2) $\dparcial H{t}(U,t)= \nu(t)\left(1- R(t)\varphi _1(U)\kapa 1(t)\right)\V1(t)+\til{ \dHd t}(U,t)$, where $$\til{ \dHd t}(U,t)\in {\rm span}\{\V2(t),\V3(t),\ldots,\V{n+1}(t)\}.$$
\end{tasks}
\eprop
\proof The properties~(i), (ii) and (iii) follow easily from those of $\fifi$, which we cited above. Differentiating (i) with respect to $t$, we get $H_0\cdot \dparcial {H_0}{t}=R(t)R'(t)$. But $\dparcial {H_0}{t}=\dparcial {H}{t}-\nu(t)\V1(t)$ and $H_0$ is perpendicular to $\V1$. Hence (iv) holds true. For (v), observe that
$${\rm span}\left\{H_0,\dparcial H{u_1},\ldots,\dparcial H{u_{n-1}}\right\}\subset{\rm span}\{\V2(t),\V3(t),\ldots,\V{n+1}(t)\}$$
and that both subspaces have the same dimension. The alternative (vi) is obtained by a direct calculation, by using the Frenet equations.\qed
\vskip 10pt

Now, we can establish the main result of this section.
\bteo \label{elementovolteo}\mbox{}\\[-5pt]
\benum
\item $\**\dbgamma=(R(t))^{n-1}\,\desfera{n-1}\sqrt{(\nu(t))^2\left(1-\varphi _1(U) R(t) \kappa_1(t)\right){}^2+(R'(t))^2}\,\dt,$
where $\desfera{n-1}$ is the volume element of the unit sphere $\esfera{n-1}$. 
\item $N(U,t)=\pm \frac{\nu(t) \left(\fifi _1 R(t) \kapa 1 (t)-1\right)}{R(t) \sqrt{\nu(t)^2 \left(\fifi _1 R(t) \kapa 1 (t)-1\right)^2+R'(t)^2}} \, \left(H_0+\frac{R(t) R'(t)}{\nu(t) \left(\fifi _1 R(t) \kapa 1 (t)-1\right)} \,\V 1(t)\right).$
\eenum
In particular, if $\bgamma$ is of revolution, then
$$\**\dbgamma=(R(t))^{n-1}\,\desfera{n-1}\sqrt{1+(R'(t))^2}\,\dt,$$
where $\desfera{n-1}$ is the volume element of the unit sphere $\esfera{n-1}$, and hence
$$\vol(\bgamma)=\vol(\esfera{n-1})\int_a^b(R(t))^{n-1}\sqrt{1+(R'(t))^2}\,\dt.$$
Also, in this case, the Gauss map is
$$N(U,t)=\pm \frac{1}{\sqrt{1+(R'(t))^2}}\parl \fifi\pupla u{n-1},-R'(t)\parr.$$

\eteo
\proof We start by setting  $U_0=H_0/R(t)$,  $a_j=\norma{\dparcial H{u_j}}^2$, $U_j=\dparcial H{u_j}/\sqrt{a_j}$, $b_j=\dparcial H{u_j}\cdot\dparcial H{t}$,  for $j$  running from 1 to $n-1$, and $a_n=\norma{\dparcial H{t}}^2$. The vectors $U_j$, $0\le j\le n-1$, form an orthonormal basis for the subspace in~(v) of the Proposition~\ref{Hprop}. We have that $(g_{ij})$, the matrix of the first fundamental form of $\bgamma$ with respect to $H$, is 
$$(g_{ij})=\left(
\begin{array}{ccccc}
a_1 & 0 & \ldots & 0 & b_1 \\
0 & a_2 & 0 & \ldots & b_2 \\
0 & 0 & \ddots & 0 &\vdots \\
0 & 0 & 0 & a_{n-1} & b_{n-1} \\
b_1 & b_2 &\ldots & b_{n-1} & a_n \\
\end{array}
\right).$$
It follows from the Lemma~\ref{matrideHlem} that 
\beq \label{detequation}\arraycolsep=0pt\def\arraystretch{3}
\begin{array}{rcl}
	\**\det(g_{ij}) &\;=\;&\**\norma{\dparcial H{u_1}}^2\norma{\dparcial H{u_2}}^2\cdots \norma{\dparcial H{u_{n-1}}}^2\,\,\left(\norma{\dparcial H{t}}^2- \sum _{j=1}^{n-1} \left(\dparcial H{t}\cdot U_j\right)^2\right)\\
	\**&\;=\;&\**\norma{\dparcial H{u_1}}^2\norma{\dparcial H{u_2}}^2\cdots \norma{\dparcial H{u_{n-1}}}^2\,\left(\norma{\dparcial H{t}}^2- \sum _{j=1}^{n-1} \left(\til{\dHd t}\cdot U_j\right)^2\right)\\
\end{array}\eeq
Note that,
$$\norma{\dparcial H{t}}^2- \sum _{j=1}^{n-1}\left( \til{\dHd t}\cdot U_j\right)^2=\nu(t)^2\left(1- R(t)\varphi _1(U)\kapa 1(t)\right)^2+\norma{\til{\dHd t}}^2- \sum _{j=1}^{n-1} \left(\til{\dHd t}\cdot U_j\right)^2,$$
by~(vi) above,
and 
$$\norma{\til{\dHd t}}^2=\left(\til{\dHd t}\cdot U_0\right)^2+\sum _{j=1}^{n-1} \left(\til{\dHd t}\cdot U_j\right)^2,$$
because $\til{\dHd t}\in{\rm span}\{U_0,U_1,\ldots,U_{n-1}\}$. Moreover,
$$\left(\til{\dHd t}\cdot U_0\right)^2=\left(\dparcial H{t}\cdot U_0\right)^2=\left(\dparcial H{t}\cdot \frac{H_0}{R(t)}\right)^2=(R'(t))^2,$$
where we use~(iv) of the Proposition~\ref{Hprop}. Putting this all together in~\refeq{detequation}, we obtain that
$$\arraycolsep=0pt\def\arraystretch{3}
\begin{array}{rcl}
\**\det(g_{ij}) &\;=\;&\**\norma{\dparcial H{u_1}}^2\norma{\dparcial H{u_2}}^2\cdots \norma{\dparcial H{u_{n-1}}}^2\left(\nu(t)^2\left(1- R(t)\varphi _1(U)\kapa 1(t)\right)^2+(R'(t))^2\right)\\
\**&\;=\;&\**\norma{R(t)\dparcial \fifi{u_1}}^2\norma{R(t)\dparcial \fifi{u_2}}^2\cdots \norma{R(t)\dparcial \fifi{u_{n-1}}}^2\left(\nu(t)^2\left(1- R(t)\varphi _1(U)\kapa 1(t)\right)^2+(R'(t))^2\right)\\
\** &\;=\;&\**(R(t))^{2(n-1)}\norma{\dparcial \fifi{u_1}}^2\norma{\dparcial \fifi{u_2}}^2\cdots \norma{\dparcial \fifi{u_{n-1}}}^2\left(\nu(t)^2\left(1- R(t)\varphi _1(U)\kapa 1(t)\right)^2+(R'(t))^2\right),\\
\end{array}$$ 
from which $\dbgamma$ can be easily seen. For the Gauss map, we start by remembering that the vectors $H_0$ and  $\V1$ are perpendicular to $\dHd{u_j}$, for $1\le j\le n-1$, and that  $H_0\cdot \dHd t=RR'$.  Thus, it is natural try to find some $\lambda\in\R{}$ such that $H_0+\lambda \V1$ also is  perpendicular to $\dHd t$. By solving the equation $(H_0+\lambda \V1)\cdot \dHd t=0$, we get that $$\lambda=\frac{R(t)\, R'(t)}{\nu(t)\, \left(\fifi _1(U) R(t)\, \kapa 1(t)-1\right)}.$$
Hence $N$ equals $\pm\frac{H_0+\lambda \V1}{\norma{H_0+\lambda \V1}}$. For the case when the hypersurface $\bgamma$ is of revolution, we just do  $f(t)=(0,0,\ldots,t)$, $\nu=1$, $\V1=(0,0,\ldots,1)$, $\kapa 1=0$ and $\V j=e_{j-1}$, for $2\le j\le n+1$, in (i) and (ii), where $e_j$ is the $j$-th element of the canonical basis of $\R {n+1}$.  To end this proof, we remark that the Gauss map $N$ is well defined due to the regularity condition~\refeq{regularitycond}, which implies that
$$\left(1-\fifi _1(U) R(t)\, \kapa 1(t)\right)>0,$$
because $|\fifi _1(U)|\le 1$. Thus $\bgamma$ is, in fact, a piece of a regular hypersurface of $\R{n+1}$. 
We are done.\qed
\vskip 10pt
We have the following corollary which shows that the first theorem of Pappus holds true for sphere tubes  of radius constant in $\R{n+1}$.
\bcoro If $R(t)=R$ is constant, then $\vol(\bgamma)=\vol(\esfera{n-1}(R))\,l(f,I)$, where $\esfera{n-1}(R)$ is the sphere of radius $R$ centered at origin.
\ecoro
\proof From (i) of the Theorem~\ref{elementovolteo}, we get that 
$$\arraycolsep=0pt\def\arraystretch{3}
\begin{array}{rcl}
\** \vol(\bgamma)&\;=\;&\**\int_{\bgamma}\dbgamma=R^{n-1}\int_a^b\nu(t)\parl
\int_{\esfera{n-1}}\desfera{n-1}
\left(1-\varphi _1(U) R\,\kappa_1(t)\right)\parr\dt\\
\**&\;=\;&\** R^{n-1}\vol(\esfera{n-1})\,l(f,I)-R^{n}\int_a^b\nu(t)\kappa_1(t)\parl
\int_{\esfera{n-1}}\varphi _1(U)\,\desfera{n-1}
\parr\dt\\[-10pt]
\** &\;=\;&\**\**\vol(\esfera{n-1}(R))\,l(f,I),
\end{array}$$
because $\int_{\esfera{n-1}}\varphi _1(U)\,\desfera{n-1}=0$, according to the next lemma.
\blema  $\int_{\esfera{n-1}}x_j\,\desfera{n-1}=0$, for $1\le j\le n$, that is, the barycenter of the sphere is its center.
\elema
\proof Consider the $(n-1)$-form $\omega=\sum_{i=1}^{n}(-1)^{n-1}x_i\,\dx1\wedge\dx2\wedge\cdots\wedge\widehat{\dx i}\wedge\cdots\wedge\dx n$, 
where the hat $\;\widehat{\;\;}\;$ over a 1-form means that the correspondent  1-form is excluded from the wedge product. As we know, the restriction of $\omega$ to $\esfera{n-1}$ equals $\desfera{n-1}$ and ${\rm d}\omega=n\dx1\wedge\dx2\wedge\ldots\wedge\dx n$. Hence,
$$\arraycolsep=0pt\def\arraystretch{2}
\begin{array}{rcl}
\** {\rm d}(x_j\,\desfera{n-1})&\;=\;&\**{\rm d}(x_j\,\omega)=\dx j\wedge\omega+n\,x_j\,\dx1\wedge\dx2\wedge\ldots\wedge\dx n\\
\**&\;=\;&\**(-1)^{j-1}\,\dx j\wedge(\dx1\wedge\dx2\wedge\cdots\wedge\widehat{\dx j}\wedge\cdots\wedge\dx n)+n\,x_j\,\dx1\wedge\dx2\wedge\ldots\wedge\dx n\\
\** &\;=\;&\**(n+1)\,x_j\,\dx1\wedge\dx2\wedge\ldots\wedge\dx n.
\end{array}$$
Now, the theorem of Stokes gives us that 
$$\int_{\esfera{n-1}}x_j\,\desfera{n-1}=(n+1)\int_{B^{[n]}}\,x_j\,\dx1\dx2\ldots\dx n=0,$$
by using the Proposition~\ref{barydisk}. Here $B^{[n]}\subset\R n$ is the compact unit disk.\qed
\brem The same arguments as in the proof above show that $\vol(\esfera{n-1}(R))=n\,R^{n-1}\,\vol(B^{[n]})$.
\erem

\end{document}

%% file: matrizes.tex
\def\dongamum{
\left(
\begin{array}{cccccccc}
 0 & -a_1^2 & 0 & a_1^4 & \ldots & 0 & (-1)^m a_1^{2 m} & 0 \\
 a_1 & 0 & -a_1^3 & 0 & \ldots & (-1)^{m-1} a_1^{2 m-1} & 0 & (-1)^m a_1^{2 m+1} \\
 0 & -a_2^2 & 0 & a_2^4 & \ldots & 0 & (-1)^m a_2^{2 m} & 0 \\
 a_2 & 0 & -a_2^3 & 0 & \ldots & (-1)^{m-1} a_2^{2 m-1} & 0 & (-1)^m a_2^{2 m+1} \\
 0 & -a_3^2 & 0 & a_3^4 & \ldots & 0 & (-1)^m a_3^{2 m} & 0 \\
 \vdots & \vdots & \vdots & \vdots & \vdots & \ddots & \vdots & \vdots \\
 0 & -a_m^2 & 0 & a_m^4 & \ldots & 0 & (-1)^m a_m^{2 m} & 0 \\
 a_m & 0 & -a_m^3 & 0 & \ldots & (-1)^{m-1} a_m^{2 m-1} & 0 & (-1)^m a_m^{2 m+1} \\
 b & 0 & 0 & 0 & \ldots & 0 & 0 & 0 \\
\end{array}
\right)}
\def\dongacinco{\left(
\begin{array}{ccccc}
 0 & -a_1^2 r_1 & 0 & a_1^4 r_1 & 0 \\
 a_1 r_1 & 0 & -a_1^3 r_1 & 0 & a_1^5 r_1 \\
 0 & -a_2^2 r_2 & 0 & a_2^4 r_2 & 0 \\
 a_2 r_2 & 0 & -a_2^3 r_2 & 0 & a_2^5 r_2 \\
 b & 0 & 0 & 0 & 0 \\
\end{array}
\right)}
\def\Sdetcinco{\left(
\begin{array}{ccccc}
 \cos \left(t a_1\right) & -\sin \left(t a_1\right) & 0 & 0 &0\\
 \sin \left(t a_1\right) & \cos \left(t a_1\right) & 0 & 0 &0 \\
 0 & 0 & \cos \left(t a_2\right) & -\sin \left(t a_2\right)  &0\\
 0 & 0 & \sin \left(t a_2\right) & \cos \left(t a_2\right) &0 \\
 0 & 0 &0 &0 &1 \\
\end{array}
\right)}
\def\colcan{\left(
\begin{array}{c}
 e_1 \\
 e_2 \\
 e_3 \\
 e_4 \\
\end{array}
\right)}
\def\sencos{\left(
\begin{array}{cccc}
 \cos \left(a_1 t\right) & \sin \left(a_1 t\right) & 0 & 0 \\
 -\sin \left(a_1 t\right) & \cos \left(a_1 t\right) & 0 & 0 \\
 0 & 0 & \cos \left(a_2 t\right) & \sin \left(a_2 t\right) \\
 0 & 0 & -\sin \left(a_2 t\right) & \cos \left(a_2 t\right) \\
\end{array}
\right)}
\def\transmatrizQ{\left(
\begin{array}{cccc}
 q_{11} & q_{21} & q_{31} & q_{41} \\
 q_{12} & q_{22} & q_{32} & q_{42} \\
 q_{13} & q_{23} & q_{33} & q_{43} \\
 q_{14} & q_{24} & q_{34} & q_{44} \\
\end{array}
\right)}
\def\sys4{\left(
\begin{array}{cccc}
 0 & \nu \kappa _1 & 0 & 0 \\
 -\nu \kappa _1 & 0 & \nu \kappa _2 & 0 \\
 0 & -\nu \kappa _2 & 0 & \nu \kappa _3 \\
 0 & 0 & -\nu \kappa _3 & 0 \\
\end{array}
\right)\left(
\begin{array}{c}
 V_1 \\
 V_2 \\
 V_3 \\
 V_4 \\
\end{array}
\right)=\left(
\begin{array}{c}
 V'_1 \\
 V'_2 \\
 V'_3 \\
 V'_4 \\
\end{array}
\right)}
\def\dongaminor{
\left(
\begin{array}{ccccccc}
 1 & 0 & -a_1^2 & 0 &\ldots & (-1)^{m-1} a_1^{2 m-2}  \\
 0 & -a_2^2 & 0 & a_2^4 &\ldots & 0 \\
 a_2 & 0 & -a_2^3 & 0 &\ldots & (-1)^{m-1} a_2^{2 m-1}  \\
 0 & -a_3^2 & 0 & a_3^4 &\ldots & 0  \\
 \vdots & \vdots & \vdots & \vdots & \vdots & \ddots  \\
 0 & -a_m^2 & 0 & a_m^4 &\ldots & 0  \\
 a_m & 0 & -a_m^3 & 0 &\ldots & (-1)^{m-1} a_m^{2 m-1} \\
\end{array}
\right)}

\def\dongam{
\left(
\begin{array}{ccccccc}
 0 & -a_1^2 & 0 & a_1^4 &\ldots & 0 & (-1)^m a_1^{2 m} \\
 a_1 & 0 & -a_1^3 & 0 &\ldots & (-1)^{m-1} a_1^{2 m-1} & 0 \\
 0 & -a_2^2 & 0 & a_2^4 &\ldots & 0 & (-1)^m a_2^{2 m} \\
 a_2 & 0 & -a_2^3 & 0 &\ldots & (-1)^{m-1} a_2^{2 m-1} & 0 \\
 0 & -a_3^2 & 0 & a_3^4 &\ldots & 0 & (-1)^m a_3^{2 m} \\
 \vdots & \vdots & \vdots & \vdots & \vdots & \ddots & \vdots \\
 0 & -a_m^2 & 0 & a_m^4 &\ldots & 0 & (-1)^m a_m^{2 m} \\
 a_m & 0 & -a_m^3 & 0 &\ldots & (-1)^{m-1} a_m^{2 m-1} & 0 \\
\end{array}
\right)}

\def\donga4{\left(
\begin{array}{cccc}
 0 & -a_1^2 r_1 & 0 & a_1^4 r_1 \\
 a_1 r_1 & 0 & -a_1^3 r_1 & 0 \\
 0 & -a_2^2 r_2 & 0 & a_2^4 r_2 \\
 a_2 r_2 & 0 & -a_2^3 r_2 & 0 \\
\end{array}
\right)}
\def\Sdet{\left(
\begin{array}{cccc}
 \cos \left(t a_1\right) & -\sin \left(t a_1\right) & 0 & 0 \\
 \sin \left(t a_1\right) & \cos \left(t a_1\right) & 0 & 0 \\
 0 & 0 & \cos \left(t a_2\right) & -\sin \left(t a_2\right) \\
 0 & 0 & \sin \left(t a_2\right) & \cos \left(t a_2\right) \\
\end{array}
\right)}
\def\usual3{\left(
\begin{array}{ccccccccc}
 0 & 0 & 0 & \frac{1}{\sqrt{2}} & 0 & 0 & 0 & 0 & 0 \\
 0 & 0 & 0 & 0 & \frac{1}{\sqrt{2}} & 0 & 0 & 0 & 0 \\
 0 & 0 & 0 & 0 & 0 & \frac{1}{\sqrt{2}} & 0 & 0 & 0 \\
 -\frac{1}{\sqrt{2}} & 0 & 0 & 0 & 0 & 0 & \frac{1}{\sqrt{2}} & 0 & 0 \\
 0 & -\frac{1}{\sqrt{2}} & 0 & 0 & 0 & 0 & 0 & \frac{1}{\sqrt{2}} & 0 \\
 0 & 0 & -\frac{1}{\sqrt{2}} & 0 & 0 & 0 & 0 & 0 & \frac{1}{\sqrt{2}} \\
 0 & 0 & 0 & -\frac{1}{\sqrt{2}} & 0 & 0 & 0 & 0 & 0 \\
 0 & 0 & 0 & 0 & -\frac{1}{\sqrt{2}} & 0 & 0 & 0 & 0 \\
 0 & 0 & 0 & 0 & 0 & -\frac{1}{\sqrt{2}} & 0 & 0 & 0 \\
\end{array}
\right) \left(
\begin{array}{c}
 V_{11} \\
 V_{12} \\
 V_{13} \\
 V_{21} \\
 V_{22} \\
 V_{23} \\
 V_{31} \\
 V_{32} \\
 V_{33} \\
\end{array}
\right)=\left(
\begin{array}{c}
 V_ {11}' \\
 V_ {12}' \\
 V_ {13}' \\
 V_ {21}' \\
 V_ {22}' \\
 V_ {23}' \\
 V_ {31}' \\
 V_ {32}' \\
 V_ {33}' \\
\end{array}
\right)}
\def\usualsode{\left(
\begin{array}{ccccccccc}
 0 & \nu\kapa 1 & 0 & \nu\kapa 1 & 0 & 0 & 0 & 0 & 0 \\
 -\nu\kapa 1 & 0 & 0 & 0 & \nu\kapa 1 & 0 & \nu\kapa 2 & 0 & 0 \\
 0 & 0 & 0 & -\nu\kapa 2 & 0 & 0 & 0 & \nu\kapa 1 & 0 \\
 -\nu\kapa 1 & 0 & \nu\kapa 2 & 0 & \nu\kapa 1 & 0 & 0 & 0 & 0 \\
 0 & -\nu\kapa 1 & 0 & -\nu\kapa 1 & 0 & \nu\kapa 2 & 0 & \nu\kapa 2 & 0 \\
 0 & 0 & 0 & 0 & -\nu\kapa 2 & 0 & -\nu\kapa 1 & 0 & \nu\kapa 2 \\
 0 & -\nu\kapa 2 & 0 & 0 & 0 & \nu\kapa 1 & 0 & 0 & 0 \\
 0 & 0 & -\nu\kapa 1 & 0 & -\nu\kapa 2 & 0 & 0 & 0 & \nu\kapa 2 \\
 0 & 0 & 0 & 0 & 0 & -\nu\kapa 2 & 0 & -\nu\kapa 2 & 0 \\
\end{array}
\right)  \left(
\begin{array}{c}
 V_{11} \\
 V_{21} \\
 V_{31} \\
 V_{12} \\
 V_{22} \\
 V_{32} \\
 V_{13} \\
 V_{23} \\
 V_{33} \\
\end{array}
\right)=\left(
\begin{array}{c}
 V_{11}' \\
 V_{21}' \\
 V_{31}' \\
 V_{12}' \\
 V_{22}' \\
 V_{32}' \\
 V_{13}' \\
 V_{23}' \\
 V_{33}' \\
\end{array}
\right)}
\def\Vv(#1,#2){V_{#1#2}}
\def\Wv(#1,#2){W_{#1#2}}
\def\WW{\left(
\begin{array}{c}
 \Wv(1,1) \\
 \Wv(1,2) \\
 \Wv(1,3) \\
 \Wv(1,4) \\
 \Wv(2,1) \\
 \Wv(2,2) \\
 \Wv(2,3) \\
 \Wv(2,4) \\
 \Wv(3,1) \\
 \Wv(3,2) \\
 \Wv(3,3) \\
 \Wv(3,4) \\
\end{array}
\right)}
\def\AA{\left(
\begin{array}{ccc}
 a_{11} & a_{12} & a_{13} \\
 a_{21} & a_{22} & a_{23} \\
\end{array}
\right)}
\def\Vs{\left(
\begin{array}{c}
 V_1 \\
 V_2 \\
 V_3 \\
\end{array}
\right)}
\def\VV{\left(
\begin{array}{c}
 \Vv(1,1) \\
 \Vv(1,2) \\
 \Vv(1,3) \\
 \Vv(1,4) \\
 \Vv(2,1) \\
 \Vv(2,2) \\
 \Vv(2,3) \\
 \Vv(2,4) \\
 \Vv(3,1) \\
 \Vv(3,2) \\
 \Vv(3,3) \\
 \Vv(3,4) \\
\end{array}
\right)}
\def\Mbidexp{\left(
\begin{array}{ccc}
 a_{11} \left(
\begin{array}{cccc}
 1 & 0 & 0 & 0 \\
 0 & 1 & 0 & 0 \\
 0 & 0 & 1 & 0 \\
 0 & 0 & 0 & 1 \\
\end{array}
\right) & a_{12} \left(
\begin{array}{cccc}
 1 & 0 & 0 & 0 \\
 0 & 1 & 0 & 0 \\
 0 & 0 & 1 & 0 \\
 0 & 0 & 0 & 1 \\
\end{array}
\right) & a_{13} \left(
\begin{array}{cccc}
 1 & 0 & 0 & 0 \\
 0 & 1 & 0 & 0 \\
 0 & 0 & 1 & 0 \\
 0 & 0 & 0 & 1 \\
\end{array}
\right) \\
 a_{21} \left(
\begin{array}{cccc}
 1 & 0 & 0 & 0 \\
 0 & 1 & 0 & 0 \\
 0 & 0 & 1 & 0 \\
 0 & 0 & 0 & 1 \\
\end{array}
\right) & a_{22} \left(
\begin{array}{cccc}
 1 & 0 & 0 & 0 \\
 0 & 1 & 0 & 0 \\
 0 & 0 & 1 & 0 \\
 0 & 0 & 0 & 1 \\
\end{array}
\right) & a_{23} \left(
\begin{array}{cccc}
 1 & 0 & 0 & 0 \\
 0 & 1 & 0 & 0 \\
 0 & 0 & 1 & 0 \\
 0 & 0 & 0 & 1 \\
\end{array}
\right) \\
\end{array}
\right)}
\def\Mbid{\left(
\begin{array}{cccccccccccc}
 a_{11} & 0 & 0 & 0 & a_{12} & 0 & 0 & 0 & a_{13} & 0 & 0 & 0 \\
 0 & a_{11} & 0 & 0 & 0 & a_{12} & 0 & 0 & 0 & a_{13} & 0 & 0 \\
 0 & 0 & a_{11} & 0 & 0 & 0 & a_{12} & 0 & 0 & 0 & a_{13} & 0 \\
 0 & 0 & 0 & a_{11} & 0 & 0 & 0 & a_{12} & 0 & 0 & 0 & a_{13} \\
 a_{21} & 0 & 0 & 0 & a_{22} & 0 & 0 & 0 & a_{23} & 0 & 0 & 0 \\
 0 & a_{21} & 0 & 0 & 0 & a_{22} & 0 & 0 & 0 & a_{23} & 0 & 0 \\
 0 & 0 & a_{21} & 0 & 0 & 0 & a_{22} & 0 & 0 & 0 & a_{23} & 0 \\
 0 & 0 & 0 & a_{21} & 0 & 0 & 0 & a_{22} & 0 & 0 & 0 & a_{23} \\
\end{array}
\right)}
\def\Mav{\left(
\begin{array}{c}
 a_{11} \left(
\begin{array}{c}
 \Vv(1,1) \\
 \Vv(1,2) \\
 \Vv(1,3) \\
 \Vv(1,4) \\
\end{array}
\right)+a_{12} \left(
\begin{array}{c}
 \Vv(2,1) \\
 \Vv(2,2) \\
 \Vv(2,3) \\
 \Vv(2,4) \\
\end{array}
\right)+a_{13} \left(
\begin{array}{c}
 \Vv(3,1) \\
 \Vv(3,2) \\
 \Vv(3,3) \\
 \Vv(3,4) \\
\end{array}
\right) \\
 a_{21} \left(
\begin{array}{c}
 \Vv(1,1) \\
 \Vv(1,2) \\
 \Vv(1,3) \\
 \Vv(1,4) \\
\end{array}
\right)+a_{22} \left(
\begin{array}{c}
 \Vv(2,1) \\
 \Vv(2,2) \\
 \Vv(2,3) \\
 \Vv(2,4) \\
\end{array}
\right)+a_{23} \left(
\begin{array}{c}
 \Vv(3,1) \\
 \Vv(3,2) \\
 \Vv(3,3) \\
 \Vv(3,4) \\
\end{array}
\right) \\
\end{array}
\right)}